\documentclass[11pt]{article}
\usepackage{amssymb}
\usepackage{latexsym}

\newtheorem{theorem}{Theorem}[section]
\newtheorem{lemma}[theorem]{Lemma}
\newtheorem{conjecture}[theorem]{Conjecture}
\newtheorem{remark}[theorem]{Remark}
\newtheorem{corollary}[theorem]{Corollary}
\newtheorem{example}[theorem]{Example}

\newcommand {\lls}{\lessapprox}
\newcommand {\ggs}{\gtrapprox}
\newcommand {\nat}{{\mathbb N}}
\newcommand {\ap}{{\em a.p.}}
\newcommand {\bp}{{\em b.p.}}
\newcommand {\mod}{\mbox{mod }}

\renewcommand {\thefootnote}{\fnsymbol{footnote}}

\begin{document}

\baselineskip=18pt

\begin{center}
\bigskip
\bigskip
   {\Large\bf Characterizing the structure of $A$\\
when the ratio $|2A|/|A|$ is bounded by $3+\epsilon$}
  \end{center}
\bigskip
\bigskip
  \begin{center}
Renling Jin\footnotetext{The author was supported in part by the
NSF grant DMS--\#0070407 when he wrote the first draft of this paper two years ago.}\\
Department of Mathematics\\
College of Charleston\\
Charleston, SC 29424\\
{\tt jinr@cofc.edu}
  \end{center}

\renewcommand{\thefootnote}{\arabic{footnote}}\addtocounter{footnote}{0}

  \bigskip

  \begin{quote}

    \centerline{Abstract}
\bigskip
\small
Let $\nat$ be the set all of non-negative integers, let
$A\subseteq\nat$ be a finite set, and let $2A$ be the set of all
numbers of form $a+b$ for each $a$ and $b$ in $A$. In \cite{freiman2} 
the arithmetic structure of $A$ was accurately characterized
when (i) $|2A|\leqslant 3|A|-4$, (ii) $|2A|=3|A|-3$, or (iii)
$|2A|=3|A|-2$. It is also suggested in \cite{freiman2} that for
characterizing the arithmetic structure of $A$ when $|2A|\geqslant
3|A|-1$, analytic methods need to be used. However, 
the interesting and more general results in
\cite{freiman2}, which use analytic methods, no longer give the
arithmetic structure of $A$ as precise as the results mentioned
above. In this paper we characterize, with the help of nonstandard 
analysis, the arithmetic structure of $A$ along the same 
lines as Freiman's results mentioned above when $|2A|=3|A|-3+b$, 
where $b$ is positive but not too large. Precisely,
we prove that there is a real number $\epsilon >0$ and
there is a $K\in\nat$ such that if $|A|>K$ 
and $|2A|=3|A|-3+b$ for $0\leqslant b\leqslant
\epsilon |A|$, then $A$ is either a subset of an arithmetic progression
of length at most $2|A|-1+2b$ or a subset of a bi-arithmetic
progression\footnote{See the definition in the beginning of
Introduction.} of length at most $|A|+b$. 

\end{quote}
\newpage
\begin{quote}
{\center\tableofcontents}
\end{quote}
\bigskip

\section{Introduction}

Inverse problems study the structural properties of the sets $A_i$
when the sum of the sets $\sum_{i=1}^{n}A_i= \{\sum_{i=1}^k
a_i:a_i\in A_i\}$ satisfies certain conditions. When $A_i=A$ for
every $i$, we write $nA$ for $\sum_{i=1}^{n}A_i$. Note that the term
$nA$ should not be confused with the term $n\ast A=\{an:a\in A\}$,
which will also be used later in this paper. For a number
$x$ we write $x\pm A$ for the set $\{x\}\pm A$ and
write $A\pm x$ for the set $A\pm\{x\}$. G. A. Freiman
and many others have studied inverse problems for the addition of
finite sets and have obtained many results showing that if $A+B$ 
is small relative to the size of $A$ and the size
of $B$, then $A$ and $B$ must have some arithmetic structure (cf.
\cite{nathanson, DLY}). In this paper we consider the addition of two
copies of the same finite set $A$ of natural numbers. Let
$a,d,k\in\nat$ with $d,k\geqslant 1$. A set of the form 
$\{a,a+d,a+2d,\ldots,a+(k-1)d\}$
is called an {\em arithmetic progression} of length $k$ with
difference $d$. A set of the form $I\cup J$
is called a {\em bi-arithmetic progression} of length $k$ with
difference $d$ if both $I$ and $J$ are arithmetic progressions of
difference $d$, $|I|+|J|=k$, and $I+I$, $I+J$, $J+J$ are pairwise
disjoint. We will write \ap\ and \bp\ as an abbreviation for
``arithmetic progression'' and ``bi-arithmetic progression'', respectively.
For two integers $m,n$ the term $[m,n]$
represents exclusively the interval of integers. For a set
of integers $A$, we write $A[m,n]$ for the set $A\cap [m,n]$ and
$A(m,n)$ for the number $|A[m,n]|$. The reader needs to be able to
distinguish $2A(a,b)$, which is $2$ times the number $A(a,b)$, from 
$(2A)(a,b)$, which is the number of elements in the set $(2A)\cap [a,b]$.

Suppose $|A|=k$. It is well known that if $|2A|=2k-1$, then $A$
must be an \ap\, Note that it is always true that $|2A|\geqslant
2k-1$. In \cite{freiman2} Freiman obtained the interesting
generalizations of these facts by showing that

(1) if $k>3$ and $|2A|=2k-1+b<3k-3$, then $A$ is a subset of an
\ap\ of length at most $k+b$;

(2) if $k>6$ and $|2A|=3k-3$, then either $A$ is a subset of an
\ap\ of length at most $2k-1$ or $A$ is a \bp

In \cite{freiman2} the structure of $A$ was also characterized when
$|A|>10$ and $|2A|=3k-2$. The proof of the $3k-3$ theorem above in
\cite{freiman2} was not short while the proof of the $3k-2$
theorem was omitted there because, commented by Freiman, it was too
tedious\footnote{The conclusion of Freiman's $3k-2$ Theorem in
\cite{freiman2} seems missing at least one case. For example, if
$A=[0,k-3]\cup\{4k,4k+2\}$, then $|2A|=3k-2$. This case of $A$ was
not covered by the theorem.}. There has
been no further accurate characterization, until now, of the structure
of $A$ when, for example, $|2A|=3k-1$. In fact, Freiman
made the following conjecture a few years ago in \cite{freiman1}.

\begin{conjecture}[G. A. Freiman]\label{conjecture}

There exists a natural number $K$ such that for any finite set of
natural numbers $A$ with $|A|=k>K$ and $|2A|=
3k-3+b$ for $0\leqslant b<\frac{1}{3}k-2$, $A$ is either a
subset of an \ap\ of length at most $2k-1+2b$ or a subset of a
\bp\ of length at most $k+b$.

\end{conjecture}

Note that Conjecture \ref{conjecture} is clearly false if $b=\frac{1}{3}k-2$
as shown in the following easy example.

\begin{example} 

Suppose $k=3a$ and $c>2k$. Let
$A=[0,a-1]\cup [c,c+a-1]\cup [2c,2c+a-1]$. Then $|A|=k$ and
$|2A|=3|A|-3+\frac{1}{3}k-2$. But $A$ is neither a
subset of an \ap\ of length $2k-1+2(\frac{1}{3}k-2)$ nor
a subset of a \bp\ of length $k+(\frac{1}{3}k-2)$.

\end{example}

It is easy to prove Freiman's conjecture if one
adds an extra condition that the set $A$ is a subset of a \bp\, 
We prove this in Theorem \ref{subsetofbp} as a simple 
consequence of Theorem \ref{twolines}.

Let $A$ and $B$ be two subsets of two torsion--free groups,
respectively. A bijection $\phi:A\mapsto B$ is called a 
{\em $F_2$-isomorphism}\footnote{$F$ is the initial of Freiman.} 
if for all $a_1,a_2,a_3,a_4\in A$, $a_1+a_2=a_3+a_4$
if and only if $\phi(a_1)+\phi(a_2)=\phi(a_3)+\phi(a_4)$. A set
\[P=P(x_0;x_1,x_2;b_1,b_2)=\{x_0+ix_1+jx_2: 0\leqslant i<b_1\mbox{
and }0\leqslant j<b_2\}\] with $b_1\geqslant b_2>0$ is called a 
{\em $F_2$-progression} if the map $\phi:[0,b_1-1]\times [0,b_2-1] \mapsto P$ 
with $\phi(i,j)=x_0+ix_1+jx_2$ is a $F_2$-isomorphism.
$P$ is called to have rank 2 if $b_2>1$ and rank 1 if $b_2=1$.

\begin{theorem}\label{subsetofbp}

Suppose $A$ is a subset of a \bp\ $I\cup J$ such that $|A|=k>10$ 
and both $A\cap I$ and $A\cap J$ are non-empty. If
$|2A|=3k-3+b$ for $0\leqslant b<\frac{1}{3}k-2$, then
$I$ and $J$ can be chosen so that $|I|+|J|\leqslant k+b$.

\end{theorem}

\noindent {\bf Proof}:\quad Without loss of generality
we can assume that $I\cup J$ has the shortest length. Clearly, $I\cup J$
is $F_2$-isomorphic to the set 
\begin{equation}\label{l1l2}
M=\{(0,0),(1,0),\ldots,(l_1-1,0)\}\cup\{(0,1),(1,1),\ldots,(l_2-1,1)\}
\end{equation}
in ${\mathbb Z}^2$ where 
$l_1$ is the length of $I$ and $l_2$ is the length of $J$. 
Let $\phi$ be the $F_2$-isomorphism from $I\cup J$ to $M$. Then 
$|\phi(A)|=k$ and $|2\phi(A)|=|2A|=3|A|=3k-3+b<\frac{10}{3}k-5$. 
By Theorem \ref{twolines} we have that $l_1+l_2\leqslant k+b$. 
Hence $A$ is a subset of a \bp\ of length at most $k+b$.
\quad $\Box$(Theorem \ref{subsetofbp})

\medskip

It is worth to mention another interesting generalization of
Freiman's $3k-3$ Theorem in \cite{HP}, where the condition
$|2A|=3k-3$ is replaced by $|A+t\ast A|=3k-3$ for an integer
$t$. The most interesting case of this generalization is when $t=-1$. 
However, this generalization does not concern the case
when $|2A|=3k-3+b$ with $b>0$.
Recently, we developed some ideas with the help of  
nonstandard analysis in the research of the inverse problem 
for upper asymptotic density \cite{jin3} and  
found that these ideas can be applied to the case
when $|2A|=3k-3+b$ with some relatively small $b>0$.
The following is the main result of this paper.

\begin{theorem}\label{main}

There exists a positive real number $\epsilon$ and 
a natural number $K$ such that for every finite set of
natural numbers $A$ with $|A|=k$, if $k>K$ and $|2A|=
3k-3+b$ for $0\leqslant b\leqslant\epsilon k$, then $A$ is either
a subset of an \ap\ of length at most $2k-1+2b$ or a subset of
a \bp\ of length at most $k+b$.

\end{theorem}

Theorem \ref{main} gives an affirmative answer to 
Conjecture \ref{conjecture} when $0\leqslant b\leqslant\epsilon |A|$. 
Note that we have a new result even when $b=2$ . Note also that the
upper bound $2k-1+2b$ of the length of the \ap\ and the upper bound 
$k+b$ of the length of the \bp\
in Theorem \ref{main} are optimal as shown in the following two 
easy examples.

\begin{example}

For $k>15$ let $A=[0,k-3]\cup\{k+10,2k+20\}$. Then $|A|=k$ and
$|2A|=3k-3+11$. The shortest \ap\ containing $A$ has
length $2k-1+2\times 11$ and $A$ is not a subset of a \bp\
of length $k+11$.

\end{example}

\begin{example}

For $k>14$ let $A=[0,k-3]\cup\{3k,3k+12\}$. Then $|A|=k$ and
$|2A|=3k-3+11$. The shortest \bp\ containing $A$ has length $k+11$ 
and $A$ is not a subset of an \ap\ of length $2k-1+2\times 11$.

\end{example}

The proofs in this paper use methods from nonstandard analysis.
The reader is assumed to have some basic knowledge of nonstandard
analysis such as the existence of infinitesimals, differences among
standard sets, internal sets, and external sets, the transfer principle,
countable saturation, etc. For details we recommend the reader to consult
\cite{lindstrom}, \cite{henson}, \cite{jin1}, or other introductory
nonstandard analysis textbooks.

Notations involved in nonstandard methods need to be introduced. 
Some of these notations may not be common in other literature. 
We work within a fixed countably saturated nonstandard universe 
$^*\!V$. For each standard set $A\subseteq\nat$, 
we write $^*\!A$ for the nonstandard version of
$A$ in $^*\!V$. For example, $^*\!\nat$ is the set of all natural
numbers in $^*\!V$, and if $A$ is the set of all even numbers in
$\nat$, then $^*\!A$ is the set of all even numbers in $^*\!\nat$.
If we do not specify that $A,B,C,\ldots$ are sets of standard
natural numbers, then they are always assumed to be {\em internal}
subsets of $^*\!\nat$. The lower case letters are used for
integers. The integers in $^*\!\nat\smallsetminus \nat$ are called
hyperfinite integers. From now on, the letters $H$ and $N$ are
exclusively reserved for hyperfinite integers. The Greek letters
$\alpha$, $\beta$, $\gamma$, $\delta$, and $\epsilon$ are reserved
exclusively for standard real numbers.

For the convenience of handling nonstandard arguments, we introduce 
some notations of comparison (quasi-order). For real
numbers $r,s$ in $^*\!V$, by $r\approx s$ we mean that $r-s$ is an
infinitesimal, i.e. the absolute value of $r-s$ is less
than any standard positive real numbers; by $r\ll s$ ($r\gg s$) we
mean $r<s$ ($r>s$) and $r\not\approx s$; by $r\lls s$
($r\ggs s$) we mean $r<s$ ($r>s$) or $r\approx s$. Given a
hyperfinite integer $H$ and two real numbers $r,s$, by $r\sim_H s$
we mean $\frac{s-r}{H}\approx 0$; by $r\prec_H s$ ($r\succ_H
s$) we mean $r<s$ ($r>s$) and $r\not\sim_H s$; by
$r\preceq_H s$ ($r\succeq_H s$) we mean $r\prec_H s$
($r\succ_H s$) or $r\sim_H s$. It is often said that a quantity
$a$ is insignificant with respect to $H$ if $a\sim_H 0$. When using
$\sim$, $\prec$, $\preceq$, etc.\ insignificant quantities can
often be neglected. For example, instead of writing
$A(0,H)\sim_H\alpha (H+1)$, we can write its equivalent form
$A(0,H)\sim_H\alpha H$. For another example, when $a\leqslant
c\leqslant b$, we often write $A(a,c)\sim_H A(a,b)+A(b,c)$ instead
of $A(a,c)=A(a,b)+A(b+1,c)$. We will omit the subscript $H$
when it is clearly given. For a real number
$r\in\,^*\!{\mathbb R}$ bounded by a standard real number, let
$st(r)$, the {\em standard part of $r$}, be the unique standard real
number $\alpha$ such that $r\approx\alpha$. Note that $\approx$,
$\ll$, $\lls$, $\sim_H$, $\prec_H$, and $\preceq_H$ are external 
relations. 
If $A\subseteq [0,H]$ is a hyperfinite set
with $a=\min A$ and $b=\max A$, then $A$ is said to be {\em full} (in
$I$) if $A$ is a subset of an \ap\ $I$ such that
$|A|\sim I(a,b)$. We say that $A$ is {\em full in a} \bp\ $I_0\cup I_1$ 
if $A\subseteq I_0\cup I_1$ and $|A|\sim I_0(l_0,u_0)+I_1(l_1,u_1)$ where
$u_i=\max (A\cap I_i)$ and $l_i=\min (A\cap I_i)$ for $i=0,1$. 
Note that if $A\subseteq [0,H]$ be a subset of an \ap\ $I$ and $|I|\sim 0$,
then $A$ is always full.
We always assume that $A\cap I_0$ and $A\cap I_1$
are non-empty when we say that $A$ is a subset of the \bp\ $I_0\cup I_1$.

In order to apply nonstandard methods, we need to translate Theorem
\ref{main} into the following nonstandard version of it. Then we proof the
nonstandard version in the rest of the paper.

\begin{theorem}\label{nsamain}

Let $H$ be a hyperfinite integer and $A\subseteq [0,H]$ be an
internal set. Suppose $0=\min A$, $H=\max A$, $|A|\succ 0$, 
$\gcd(A)=1$, and $|2A|=3|A|-3+b$ for $0\leqslant b\sim 0$. 
Then either $H+1\leqslant 2|A|-1+2b$ or $A$ is a subset of a \bp\ 
of length at most $|A|+b$.

\end{theorem}

\noindent {\bf Proof of Theorem \ref{main} from Theorem
\ref{nsamain}}: \quad Suppose Theorem \ref{main} is not true. Then
for $\epsilon_k=\frac{1}{k}$ and $K_k=k$ for each $k\in\nat$, 
there is a finite set $A_k\subseteq [0,h_k]$ 
satisfying the following: $0=\min A_k$, $h_k=\max A_k$, $|A_k|>k$,
$\gcd(A_k)=1$, $|2A_k|=3|A_k|-3+b_k$ for 
$0\leqslant b_k\leqslant \frac{|A|}{k}$, 
$h_k+1>2|A_k|-1+2b_k$, and $A_k$ is not 
a subset of a \bp\ of length at most $|A_k|+b_k$.

Let $K$ be a hyperfinite integer and let $A=A_K$ be the term in 
the internal sequence $\langle A_k:k\in\,^*\!\nat\rangle$.
Then we have the following: $0=\min A$, $H=h_K=\max A$, $|A|>K$,
$\gcd(A)=1$, $|2A|=3|A|-3+b$ for some $b\geqslant 0$ with
$\frac{b}{|A|}\leqslant\frac{1}{K}\approx 0$, 
$H+1>2k-1+2b$, and $A$ is not a subset of
a \bp\ of length at most $|A|+b$. If in addition we have $\frac{|A|}{H}\gg
0$, then the set $A$ contradicts Theorem \ref{nsamain}. Hence it
suffices to prove $\frac{|A|}{H}\gg 0$ or equivalently $|A|\succ 0$.

Suppose $|A|\sim 0$. By Theorem \ref{bilufreiman} the set
$A$ is a subset of a $F_2$-sequence $P=P(x_0;x_1,x_2;b_1,b_2)$ such that
$\frac{|A|}{|P|}\gg 0$. If $P$ has rank 1, then $P$ is an \ap\,  Since
$\gcd(A)=1$, then $[0,H]\subseteq P$. This contradicts 
$|A|\sim 0$. Hence we can assume that $P$ has rank 2.
Let $\phi:P\mapsto [0,b_1-1]\times [0,b_2-1]$ be a $F_2$-isomorphism
and $B=\phi(A)$. Then $B$ is not a subset of a straight line.
Since $B$ is a $F_2$-isomorphic image of $A$, we have 
$|2B|=|2A|$. Hence by Theorem \ref{twolines}, $B$ is $F_2$-isomorphic
to a subset of $M$ in (\ref{l1l2}) such that $l_1+l_2\leqslant |B|+b$. This
shows that $A$ is a subset of a \bp\ of length at most $|A|+b$, 
which contradicts the assumption that $A$ is not a subset of a \bp\ of
length at most $|A|+b$.
\quad $\Box$(Theorem \ref{main})

\medskip

The approach of eliminating the possibility of 
$|A|\sim 0$ in the proof above is from \cite{bordes}. 
In fact the same approach can be used to prove that there exists 
a small positive number $\delta$ such that Conjecture \ref{conjecture} 
is true when an extra condition $|A|<\delta(\max A -\min A)$ is added.
It is possible but much more tedious to prove $|A|\succ 0$ in the proof
above directly without citing Theorem \ref{bilufreiman}.

We prove Theorem \ref{nsamain} in the next several sections. 
The proof is done in two steps. In the first step we deal with the case
when $A\subseteq [0,H]$ contains significantly less than half of the 
elements in $[0,H]$. In the second step we deal with the case when
$A\subseteq [0,H]$ contains roughly half of the elements in $[0,H]$.
The main theorem in each step is preceded by a list of lemmas, which
prove the theorem under various circumstances.
Before these two steps we present a list of general lemmas. 
For convenience we include some existing theorems
in Appendix for quick references. In this paper, theorems, lemmas,
cases, and claims are numbered in such a way that the reader should
be able to see how they are nested.

\section{General Lemmas}

In this section we state some lemmas from the author's previous
papers without proof and state some other new lemmas with proof. 
The first lemma in this section will play an important role in the 
proof of Theorem \ref{nsamain}. It uses
a concept called cut from nonstandard analysis.

An infinite initial segment $U$ of $^*\!\nat$ is called a {\em cut} if
$U+U\subseteq U$. Clearly $U=\nat$ and $U=\,^*\!\nat$ are cuts. A cut
$U\not=\,^*\!\nat$ is external because it has no maximum element. 
For example, $\nat$ is external. For a hyperfinite integer $H$, the set
\begin{equation}\label{Hcut}
U_H=\bigcap_{n\in \nat}[0,[H/n]]
\end{equation} 
is an external cut. We often write $x>U$ for 
$x\in\,^*\!\nat\smallsetminus U$ and write $x<U$ for $x\in U$. 
Note that if $x<U$ and $y>U$, then $\frac{x}{y}\approx 0$.

Let $U$ be a cut. A \bp\ $B=I\cup J$ is called a {\em $U$--unbounded} 
\bp\ if both $I\cap U$ and $J\cap U$ are upper unbounded in $U$.
Note that a $U$--unbounded \bp\ always has the difference 
greater than $2$.

Suppose $U$ is a cut. Given a function $f:U\mapsto\,^*\!{\mathbb
R}$ (not necessarily internal) bounded by a standard real number,
the {\em lower $U$--density} of $f$ is defined by
\[\underline{d}_U(f)=\sup\{\inf\{st(f(n)):n\in
U\smallsetminus [0,m]\}:m\in U\}.\] Given a set $A\subseteq
[0,H]$, let $f_A(x)=\frac{A(0,x)}{x+1}$ for each $x\in [0,H]$. The
{\em lower $U$--density} of $A$ is defined by
\[\underline{d}_U(A)=\underline{d}_U(f_A).\]
For any $x\in\,^*\!\nat$, we define the {\em lower $(x+U)$--density} and
{\em lower $(x-U)$--density} of $A$ by
\[\underline{d}_{x+U}(A)=\underline{d}_U((A-x)\cap\,^*\nat)\] and
\[\underline{d}_{x-U}(A)=\underline{d}_U((x-A)\cap\,^*\nat).\]

\begin{remark}\label{remark1}

(1) For any $A\subseteq \nat$, $\underline{d}(A)=
\underline{d}_{\nat}(^*\!A)$, where
$\underline{d}(A)=\liminf_{n\rightarrow\infty} \frac{A(0,n-1)}{n}$ is
the standard definition of the lower asymptotic density of $A$.

(2) It is easy to check that for each $a\in U$,
\[\underline{d}_U(A+a)=\underline{d}_U(A)\] and
\[\underline{d}_U(A\smallsetminus [0,a])=
\underline{d}_U(A).\]

(3) Let $H$ be hyperfinite and $A\subseteq [0,H]$. If
$\underline{d}_U(A)>\gamma$, then there exist $x\in U$ and 
$y\in [0,H]\smallsetminus U$ such that for
any $x\leqslant z\leqslant y$, $\frac{A(0,z)}{z+1}>\gamma$.
Clearly one can find a $x\in U$ such that for every 
$z\geqslant x$ in $U$, $\frac{A(0,z)}{z+1}>\gamma$.
Now the set of all $z\in [x,H]$ such that 
$\frac{A(0,z)}{z+1}>\gamma$ is internal and contains 
all elements in $U\cap [x,H]$, hence contains all elements
in $[x,y]$ for some $y>U$.

(4) If $\underline{d}_U(A)=\alpha$, then there is a $x\in U$ such
that for every $y\in U$ with $y>x$, one has $\frac{A(0,y)}{y+1}
\ggs\alpha$. This can be proven by first choosing a $x_n\in U$
for each $n\in\nat$ such that for all $z\geqslant x_n$ in $U$, 
$\frac{A(0,z)}{z+1}>\alpha -\frac{1}{n}$ and then
choosing a $x\in U$ such that $x\geqslant x_n$ for every $n\in\nat$.
The element $x$ exists because, by countable saturation, the
cofinality of $U$ is uncountable, i.e. any countable
increasing sequence in $U$ is upper bounded in $U$. 

(5) If $\underline{d}_U(A)>\frac{1}{2}$, then there exists an $a\in U$
such that $A(0,a-1)=\frac{1}{2}a$ and $A(a,c)>\frac{1}{2}(a-c+1)$ for 
every $c\in U$ with $c\geqslant a$. As a by-product we have 
$a,a+1\in A$ and $A(a,a+3)\geqslant 3$. The existence of $a$ 
is guaranteed by $\underline{d}_U(A)>\frac{1}{2}$.

\end{remark}

From now on, the only cut we need is $U_H$ defined by (\ref{Hcut})
for a given $H$. Hence when $H$ is clearly given, the letter $U$ always 
represents the cut $U_H$. Note that with $H$ fixed we have that 
$x<U$ iff $x\sim 0$ or equivalently $x>U$ iff $x\succ 0$.

The first lemma of this section bellow is \cite[Lemma 2.12]{jin3}. 

\begin{lemma}\label{kneser}

Let $H$ be hyperfinite, $U=U_H$, and $A\subseteq\,^*\!\nat$ be
such that $0<\underline{d}_U(A)=\alpha<\frac{2}{3}$. If $A\cap U$
is neither a subset of an \ap\ of difference greater than $1$
nor a subset of a $U$--unbounded \bp, then there is a $\gamma>0$
such that for every $N>U$, there is a $K\in A\smallsetminus U$ 
with $K<N$ such that
\begin{equation}\label{3timesmore}
\frac{(2A)(0,2K)}{2K+1}\geqslant \frac{3}{2}\frac{A(0,K)}{K+1}+\gamma.
\end{equation}
\end{lemma}

The following is \cite[Lemma 2.4]{jin3}.

\begin{lemma}\label{smallbig}

Let $A\subseteq [0,H]$. Suppose $0,H\in A$. If
$0\leqslant x_1\prec x_2\leqslant H$ satisfy the following

(1) $(2A)(2x_1,2x_2)\succ 3A(x_1,x_2)$,

(2) if $0\prec x_1$, then $\gcd(A[0,x])=1$ and 
$A(0,x)\preceq\frac{1}{2} (x+1)$ for some $x\sim x_1$ in $A$,

(3) if $x_2\prec H$, then $\gcd(A[x,H]-x)=1$ and 
$A(x,H)\preceq\frac{1}{2} (H-x+1)$ for some $x\sim x_2$ in $A$,

\noindent then $|2A|\succ 3|A|$.

\end{lemma}

\begin{lemma}\label{smallbig2}

Let $A\subseteq [0,H]$. Suppose $0,H\in A$, $|A|\sim\frac{1}{2}$,
$\underline{d}_U(A)=\frac{1}{2}$, there is an $a\succ 0$ in $A$ 
such that $\gcd(A[a,H]-a)=1$, and for every $N\succ 0$ there is
a $K\in A$ with $0\prec K\leqslant N$ such that (\ref{3timesmore})
is true. Then $|2A|\succ 3|A|$.

\end{lemma}

\noindent {\bf Proof}:\quad Let $0<\epsilon<1$ be such
that for any $N\succ 0$ there is a $K\in A$ with $0\prec K\leqslant
N$ such that
\[\frac{(2A)(0,2K)}{2K+1}> \frac{(3+\epsilon)A(0,K)}{2(K+1)}.\]
Let $\delta>0$ be such that $\delta<\frac{\epsilon}{6}$ and let
$y\in A$ be such that $y\succ 0$, $y\leqslant a$, $A(0,y)\geqslant
(\frac{1}{2}-\delta)(y+1)$, and $(2A)(0,2y)>(3+\epsilon)A(0,y)$.

If $A(y,H)\preceq\frac{1}{2}(H-y)$, then the lemma follows from 
Lemma \ref{smallbig} and Theorem \ref{2k-1+b}.
So we can assume $A(y,H)\succ\frac{1}{2}(H-y)$. By Theorem \ref{A+B}
we have $|A[y,H]+A[y,H]|\succeq H-y+A(y,H)$. Hence
\begin{eqnarray*}
\lefteqn{|2A|\succeq (2A)(0,2y)+|A[y,H]+A[y,H]|}\\
& &\succeq (3+\epsilon)A(0,y)+H-y+A(y,H)\\
& &\succeq 3A(0,y)+\epsilon A(0,y)+2|A|-y+A(y,H)\\
& &\succeq 3|A|+\epsilon A(0,y)+2A(0,y)-y\\
& &\succeq 3|A|+(\epsilon +2)(\frac{1}{2}-\delta)y-y\\
& &\succeq 3|A|+(\frac{\epsilon}{2}-\epsilon\delta-2\delta)y\\
& &\succ 3|A|+(\frac{\epsilon}{2}-\frac{\epsilon}{6}
   -\frac{2\epsilon}{6})y\\
& &=3|A|.
\end{eqnarray*}
\quad $\Box$(Lemma \ref{smallbig2})

It is worth to mention here that Lemma \ref{kneser}, 
Lemma \ref{smallbig}, and Lemma \ref{smallbig2}, combined together,
will frequently be used to show $|2A|\succ 3|A|$ in various situations. 
For example, if $|A|\preceq\frac{1}{2}H$ and $A\cap U$ does not 
have ``nice arithmetic structures'', then one can find an
arbitrarily small $y\succ 0$ in $A$ such that 
$(2A)(0,2y)\succ 3A(0,y)$. By Lemma \ref{smallbig} or 
Lemma \ref{smallbig2} one needs only to make sure that 
$A[x,H]$ is not a subset of an \ap\ of difference $>1$ for
some $x\succ 0$ in order to conclude that $|2A|\succ 3|A|$. 

The next lemma is trivial and will be frequently referred as 
the pigeonhole principle.

\begin{lemma}\label{pigeon}

Let $d\geqslant 1$. Suppose $a,b\in [0,d-1]$, 
$A\subseteq a+(d\ast\,^*\nat)$, $B\subseteq b+(d\ast\,^*\nat)$,
$x\in a+(d\ast\,^*\nat)$, $y\in b+(d\ast\,^*\nat)$, and
$t\in (d\ast\,^*\nat)$. If $A(x,x+t)+A(y-t,y)>\frac{t}{d}+1$, 
then $x+y\in (2A)$.

\end{lemma}

For convenience we give a name for each of the following two
sets with special structural properties.
Let $a\prec b$ in $[0,H]$. A set $F\subseteq [a,b]$ is called a {\em
forward triangle} from $a$ to $b$ if $|F|\sim\frac{1}{2}(b-a)$ 
and for every $x$ with $a\prec x\prec b$,
$F(a,x)\succ\frac{1}{2}(x-a)$. A set $B\subseteq [a,b]$ is
called a {\em backward triangle} from $a$ to $b$ if the set
$(b+a)-B$ is a forward triangle from $a$ to $b$. By the symmetry
of the forward triangle and the backward triangle, we often 
prove a result about forward (backward) triangle and assume
the symmetric result about backward (forward) triangle without proof. 

Note that if $F$ is a forward triangle from $a$ to $b$, then
there is a $z\sim a$ such that $z,z+1\in A$ and
$A(z,z+3)\geqslant 3$. The number $z$ can be obtained by
letting $z-1$ be the greatest number in $a+U$ such that
$A(a,z-1)\leqslant\frac{1}{2}(z-a)$.

\begin{lemma}\label{triangle}

Let $A\subseteq [0,H]$ be such that $0,H\in A$ and 
$0\prec |A|\preceq\frac{1}{2}H$. Let $0<\alpha\leqslant\frac{1}{2}$
and $0\prec x\prec H$.

(1) If $A(0,x)\preceq\alpha x$ and $|A|\succeq\alpha H$,
then there exists a $y\succeq x$ such that
$A(0,y)\sim\alpha y$ and either $y\sim H$ or for any $z\succ y$ in $[0,H]$,
$A(0,z)\succ\alpha z$.

(2) If $A(0,x)\succ\frac{1}{2}x$, then there are
$0\leqslant y\prec x\prec y'\leqslant H$ such that $A(0,y)\sim\frac{1}{2}y$ and
$A[y,y']$ is a forward triangle.

(3) If $\underline{d}_U(A)>\frac{1}{2}$, then there is a $y\succ 0$ such that
$A[0,y]$ is a forward triangle.

(4) If $\underline{d}_U(A)<\frac{1}{2}$ and $|A|\sim\frac{1}{2}H$,
then there are $0\leqslant y\prec y'\leqslant H$ such that
$A(y',H)\sim\frac{1}{2}(H-y')$ and $A[y,y']$ is a backward triangle.

\end{lemma}

\noindent {\bf Proof}:\quad (1) Let \[\beta=\sup\{st(\frac{z}{H+1}):
z\in [0,H]\mbox{ and }A(0,z)\preceq\alpha z\},\]
where $st$ is the standard part map.
By the completeness of the standard real line, $\beta$ is well defined.
Let $y\in [0,H]$ be such that $\frac{y}{H+1}\approx\beta$. Clearly
$y\succeq x$.

It is easy to see that if $y\prec H$, then $A(0,z)\succ\alpha z$
for any $y\prec z\leqslant H$ by the supremality of $\beta$. It is also easy 
to see that both $A(0,y)\succ\alpha y$ and $A(0,y)\prec\alpha y$
are impossible by the fact that $\beta$ is the least upper bound.

(2) By the same idea as in (1) we can find $y'\succ x$ 
such that $A(0,y')\sim\frac{1}{2}y'$ and $A(0,z)
\succ\frac{1}{2}z$ for any $x\leqslant z\prec y'$. Let $0\leqslant y\prec x$
be such that $A(0,y)\sim\frac{1}{2}y$ and $A(0,z)\succ\frac{1}{2}z$
for any $y\prec z\leqslant x$. It is easy to see that $A(y,y')$ is a forward triangle.

(3) By the definition of $\underline{d}_U$ and (1) above there exists $y\succ 0$
such that $A(0,y)\sim\frac{1}{2}y$ and
$A(0,z)\succ\frac{1}{2}z$ for every $0\prec z\prec y$. Clearly $A[0,y]$
is a forward triangle.

(4) Choose a $x\succ 0$ such that $A(0,x)\prec\frac{1}{2}x$. Hence
$A(x,H)\succ\frac{1}{2}(H-x)$. Now the conclusion follows from (2) above
with the order of $[0,H]$ reversed.
\quad $\Box$(Lemma \ref{triangle})

\medskip

The following lemma in nonstandard analysis, which is already used
in (3) of Remark \ref{remark1}, will be frequently--
sometimes implicitly--used.

\begin{lemma}\label{overspill}

Let $X\subseteq\,^*\nat$ be a proper external initial segment 
of non-negative integers and let $A\subseteq\,^*\nat$ be an 
internal set. (a) If $A\cap X$ is upper unbounded
in $X$, then $A\smallsetminus X\not=\emptyset$. (b) If 
$A\smallsetminus X$ is lower unbounded in $^*\nat\smallsetminus X$, 
then $A\cap X\not=\emptyset$. 

\end{lemma}

\noindent {\bf Proof}:\quad If (a) of the lemma is not true, then
\[X=\{v\in\,^*\nat:(\exists x\in A)\,(v\leqslant x)\},\]
which means that $X$ is internal. The proof of (b) is similar.
\quad $\Box$(Lemma \ref{overspill})

\begin{lemma}\label{ftbt}

Suppose $a\prec b$ in $[0,H]$. 

(1) If $T$ is a forward triangle from
$a$ to $b$, then $[a',b']\subseteq (2T)$ for some $a'\sim 2a$ and
$b'\sim a+b$.

(2) If $B$ is a backward triangle from $a$ to $b$,
then $[a',b']\subseteq (2B)$ for some $a'\sim a+b$ and $b'\sim 2b$.

\end{lemma}

\noindent {\bf Proof}:\quad Given each $x$ with $a\prec x\prec b$,
since $T(a,x)\succ\frac{1}{2}(x-a)$, then by the pigeonhole
principle, $T[a,x]\cap (x+a-T[a,x])\not=\emptyset$. This implies
$x+a\in (2T)$. Since $2T$ is an internal set, then by 
Lemma \ref{overspill} there are
$a'\sim 2a$ and $b'\sim a+b$ such that $[a',b']\subseteq (2T)$.
The proof of the second part follows from the symmetry.
\quad $\Box$(Lemma \ref{ftbt})

\medskip

The following is \cite[Lemma 2.5]{jin3}

\begin{lemma}\label{block}

Let $A\subseteq [0,H]$ for a hyperfinite integer $H$. If $A[0,a]$ is
a forward triangle from $0$ to $a$ and $(2A)(a,c)\sim 0$ for some
$0\prec a\prec c$, then there is a $b\sim \frac{a}{2}$ such that 
$A[0,a]\subseteq [0,b]$.

\end{lemma}

The following is a technical lemma, which will be used in the next two sections.

\begin{lemma}\label{bpleft}

Suppose $0\prec a\prec H$, $A(0,a)\succ 0$, $\gcd(A[0,a])=1$, $A[0,a]$ is 
a subset of a \bp\ $I\cup J$ of difference $d\geqslant 3$,  
$A[0,a+1]$ is not a subset of a \bp\ of difference $d$, $|2A|\sim 3|A|$, and
$|A[a+1,H]+A[a+1,H]|\sim 3A(a+1,H)$. Then $A[0,a]$ is full in $I\cup J$
and $\max (A\cap I)\sim\max (A\cap J)\sim a$.

\end{lemma}

\noindent {\bf Proof}:\quad Let $A_0=A[0,a]\cap I$, $A_1=A[0,a]\cap J$, 
$l_i=\min A_i$, and $u_i=\max A_i$ for $i=0,1$. Since 
\begin{eqnarray*}
\lefteqn{|2A|\succeq |2A_0|+|2A_1|+|A_0+A_1|+|A[a+1,H]+A[a+1,H]|}\\
& &\succeq 3|A_0|+3|A_1|+3A(a+1,H)\sim 3|A|,
\end{eqnarray*}
then $|2A|\sim 3|A|$ implies that $|2A_i|\sim 2|A_i|$ for $i=0,1$.
Hence by Theorem \ref{2k-1+b} we have that $A_0$ is full in $I$ and
$A_1$ is full in $J$. 

Without loss of generality we assume $u_0<u_1$. If $u_1\prec a$,
then \[|2A|\succeq 3A(0,a)+3A(a+1,H)+|a+1+A[2u_1-a,u_1]|\succ 3|A|.\]
If $u_1\sim a$ and $u_0\prec a$, then 
\[|2A|\succeq 3A(0,a)+3A(a+1,H)+|a+1+A_1[2u_0-a,a]|\succ 3|A|\]
because $(a+1+A_1[2u_0-a,a])\cap (A[0,a]+A[0,a])=\emptyset$. 
Hence we have $u_i\sim a$ for $i=0,1$.
\quad $\Box$(Lemma \ref{bpleft})

\section{First Step: When $\frac{|A|}{H}$ is significantly less than $\frac{1}{2}$.}

In this section we always assume that $H$ is a hyperfinite integer,
$A\subseteq [0,H]$, $0,H\in A$, and $\gcd(A)=1$. We will prove 
Theorem \ref{nsamain} under one extra condition
\begin{equation}\label{sparse}
|A|\prec\frac{1}{2}H.
\end{equation}
We will prove that if $|2A|\sim 3|A|$ and (\ref{sparse}) is true, 
then $A$ must be a subset
of a \bp, which, by Theorem \ref{subsetofbp}, implies
Theorem \ref{nsamain}. In this section the condition $|2A|=3|A|-3+b$ 
is not explicitly used. Hence the letter $b$ is not reserved. In order to
make the lemmas in this section available for the other sections, we
will not automatically assume (\ref{sparse}). The condition (\ref{sparse}) 
will be explicitly stated when it is needed.

We will first prove various versions of the main theorem of the section 
as lemmas when some additional structural properties of $A$ are assumed.
After all needed versions are proven we combine them into the main theorem.

\begin{lemma}\label{ftap}

If there are $0\prec a<b\prec H$ in $A$ such that $A=T\cup P$
where $T=A[0,a]$ is a forward triangle from $0$ to $a$ and
$P=A[b,H]$ is a subset of an \ap\ of difference $d>1$
with $|P|\succ 0$, then either $A$ is a subset of
a \bp\ of difference $3$ or $|2A|\succ 3|A|$.

\end{lemma}

\noindent {\bf Proof}:\quad Let $P$ be a subset of an \ap\ $I$
of difference $d>1$ such that $b\in P$
is the least element of $I$, and $H\in P$ is the largest element
of $I$. Suppose $T$ is not a subset of a \bp\ of difference $3$.
Since $T$ is a forward triangle, there exist $z,z+1\in A\cap U$
such that for every $x$ with $z\leqslant x\prec a$,
\[\frac{T(z,x)}{x-z+1}>\frac{1}{2}.\] 
Without loss of generality (except in Case \ref{ftap}.1.2)
we can assume $z=0$. Under this assumption
we have $0,1\in A$ and $A(0,3)\geqslant 3$.

\medskip

{\bf Claim \ref{ftap}.1}\quad If $P$ is not full in $I$, then
$(T+P)(b,H)\succ 2|P|$.

Proof of Claim \ref{ftap}.1:\quad Since $P$ is not full, then
$|I\smallsetminus P|\succ 0$. Let $\cal I$ be the collection
of all intervals $[x,y]\subseteq [b,H]$ such that $y-x\geqslant 2d-2$,
$[x,y]\cap P=\emptyset$, and $x-1,y+1\in P$. Then
\[I\smallsetminus P=\bigcup_{[x,y]\in
{\cal I}}([x-1+d,y+1-d]\cap I).\] Hence 
\begin{eqnarray*}
\lefteqn{|I\smallsetminus P|= \sum_{[x,y]\in{\cal I}}\frac{1}{d}(y-x+2-d)}\\
& &\preceq\sum_{[x,y]\in{\cal I}}\frac{1}{d}(y-x)\preceq
\sum_{[x,y]\in{\cal I}}\frac{1}{2}(y-x)\preceq\sum_{[x,y]\in{\cal I}}(y-x-1).
\end{eqnarray*}
If there is an interval $[x,y]\in{\cal I}$ such that 
$y-x\geqslant\frac{a}{2}$, then
\begin{eqnarray*}
\lefteqn{(P+T)(b,H)}\\
& &\succeq |P|+|1+P|+|(x-1+T)(x+1,y)|\\
& &\succeq 2|P| +T(2,y-x+1)\\
& &\succeq 2|P|+\min\{|T|-2,\frac{1}{2}(y-x+2)-2\}\\
& &\succeq 2|P|+\frac{1}{4}a\succ 2|P|.
\end{eqnarray*}
So we can assume $y-x<\frac{a}{2}$ for every $[x,y]\in {\cal I}$.
Since for each interval $[x,y]\in{\cal I}$, we have
\begin{eqnarray*}
\lefteqn{(P+T)(x+1,y)\geqslant |x-1+T[2,y-x+1]|}\\
& &=T(2,y-x+1)>\frac{1}{2}(y-x+2)-2\\
& &=\frac{1}{2}(y-x-1)-\frac{1}{2}.
\end{eqnarray*}
Hence $(P+T)(x+1,y)\geqslant\frac{1}{2}(y-x-1)$ because the
left-side is an integer. So
\begin{eqnarray*}
\lefteqn{(P+T)(b,H)}\\
& &\succeq |P|+|1+P|+\sum_{[x,y]\in{\cal I}}(P+T)(x+1,y)\\
& &\succeq |P|+|1+P|+\sum_{[x,y]\in{\cal I}}\frac{1}{2}(y-x-1)\\
& &\succeq 2|P|+\frac{1}{2}|I\smallsetminus P|\succ 2|P|.
\end{eqnarray*}
\quad $\Box$(Claim \ref{ftap}.1)

\medskip 

By the claim above we can assume that $P$ is full in $I$ because
otherwise
\[|2A|\succeq a+(P+T)(b,H)+|H+A|\succ 2|T|+2|P|+|A|=3|A|.\]
Next we divide the proof of the lemma into three cases with $d=2$,
$d=3$, and $d\geqslant 4$.

\medskip

{\bf Case \ref{ftap}.1}\quad $d=2$.

Let $c\prec H$ in $P$ be such that $H-c<\frac{a}{2}$. Then
\begin{eqnarray*}
\lefteqn{|2A|\succeq a+|\{0,1\}+P|+(H+T)(H,c+a)}\\
& &\quad +(c+P)(c+b,H+b)+(H+P)(H+b,2H)\\
& &\sim 2|T|+2|P|+T(0,c+a-H)+\frac{1}{2}(H-c)+|P|\\
& &\succ 3|T|+3|P|=3|A|
\end{eqnarray*}
because $T(c+a-H+1,a)\prec\frac{1}{2}(H-c)$.
\quad $\Box$(Case \ref{ftap}.1)

\medskip

{\bf Case \ref{ftap}.2}\quad $d=3$. Let $A$ be the original set
with $z,z+1\in A\cap U$.
If $A$ is a subset of a \bp\ of difference $3$, then the lemma 
is trivially true. Suppose $A$ is not a subset of a \bp\ of difference
$3$. Let $c=\min\{x\in A:x\equiv z+2\,(\mod 3)\}$. Note that either 
$c\in T$ or $c=b$.

Suppose $c\succ 0$. Let $b\leqslant x\prec H$ be such that $x\in
A$ and $H-x<c$. Note that
$A[0,c-1]\subseteq (z+(3\ast\,^*\!\nat))\cup (z+1+(3\ast\,^*\!\nat))$ 
and $A[x,H]+c\subseteq b+z+2+(3\ast\,^*\!\nat)$. Hence
$(A[x,H]+c)\cap (H+A[0,c-1])=\emptyset$. So we have
\begin{eqnarray*}
\lefteqn{|2A|\succeq a+2|P|+|H+A|+|c+A[x,H-1]|}\\
& &\sim 2|T|+2|P|+|A|+A(x,H)\succ 3|A|.
\end{eqnarray*}
Suppose $c\sim 0$. Then
\[|2A|\succeq a+|\{z,z+1,c\}+P|+|H+A|\sim 3|A|+|P|\succ 3|A|.\]
\quad $\Box$(Case \ref{ftap}.2)

\medskip

{\bf Case \ref{ftap}.3}\quad $d\geqslant 4$.

Since $T(0,3)\geqslant 3$, then
\[|2A|\succeq a+|T[0,3]+P|+|H+A|\sim 3|A|+|P|\succ 3|A|.\]
This ends the proof of Lemma \ref{ftap}.
\quad $\Box$(Lemma \ref{ftap})

\begin{lemma}\label{ftbt2}

Suppose there are $0\prec a<b\prec H$ such that $A=F\cup B$, where
$F$ is a forward triangle from $0$ to $a$ and $B$ is a backward
triangle from $b$ to $H$. If $|2A|\sim 3|A|$, then $\bar{a}=\max
F\sim\frac{a}{2}$ and $\bar{b}=\min B\sim\frac{b+H}{2}$. Hence $F$
is full in $[0,\bar{a}]$ and $B$ is full in $[\bar{b},H]$. So $A$ is
a full subset of the \bp\ $[0,\bar{a}]\cup [\bar{b},H]$ of difference $1$.

\end{lemma}

\noindent {\bf Proof}:\quad Suppose $\bar{b}=\min B\sim b$. Let $0\prec
x\prec \min\{a,\frac{H-b}{2}\}$. Then by Lemma \ref{ftbt}
\begin{eqnarray*}
\lefteqn{|2A|\succeq a+|\bar{b}+F[0,x]|+B(\bar{b}+x,H)}\\
& &\quad +|H+F|+|[H+b,2H]|\\
& &\sim 2|F|+F(0,x)+B(\bar{b}+x,H)+|F|+2|B|\\
& &\succ 3|F|+\frac{1}{2}(x+1)+B(\bar{b}+x,H)+2|B|\\
& &\succ 3|F|+2|B|+B(\bar{b},H)\sim 3|A|.
\end{eqnarray*}
Hence we can assume that $\bar{b}\succ b$. But this implies
\begin{eqnarray*}
\lefteqn{|2A|\succeq a+(2A)(a,\bar{b})+A(\bar{b},H)+
    |H+F|+|[H+b,2H]|}\\
& &\sim 2|F|+(2A)(a,\bar{b})+|B|+|F|+2|B|\sim
3|A|+(2A)(a,\bar{b}).
\end{eqnarray*}
Hence $|2A|\sim 3|A|$ implies $(2A)(a,\bar{b})\sim 0$.
By Theorem \ref{block},
$\bar{a}\sim\frac{a}{2}$. By a symmetric argument, we can also
show that $\bar{b}\sim\frac{b+H}{2}$.\quad $\Box$(Lemma \ref{ftbt2})

\begin{lemma}\label{ftmore}

Suppose there are $0\prec a<b\prec H$ such that $A=F\cup C$,
where $F\subseteq [0,a]$ is a forward triangle from $0$ to $a$ and
$C\subseteq [b,H]$ with $b\in C$,
$|C|\preceq\frac{1}{2}(H-b+1)$, and $\gcd(C-b)=1$. Then
$|2A|\succ 3|A|$.

\end{lemma}

\noindent {\bf Proof}:\quad First we assume that there is an
$x\in C$ such that $0\prec x-b<\frac{a}{2}$. Then
\begin{eqnarray*}
\lefteqn{|2A|\succeq a+|b+F[0,x-b]|}\\
& &\quad +|x+F[0,a+b-x]|+|C[b,H]+C[b,H]|\\
& &\succeq 2|F|+F(0,x-b)+F(0,a+b-x)+3|C|\\
& &\succ 2|F|+3|C|+F(a+b-x,a)+F(0,a+b-x)\\
& &\sim 3|F|+3|C|=3|A|.
\end{eqnarray*}
If the assumption above is not true, let
$x=\min\{z\in C:z\geqslant b+\frac{a}{2}\}$ and 
$y=\max\{z\in C:z< b+\frac{a}{2}\}$. Then $y\sim b$,
$(2C)(2b,b+x)\sim 0$, and $(2C)(2b,2H)\sim (2C)(b+x,2H)$.
Hence
\begin{eqnarray*}
\lefteqn{|2A|\succeq a+|b+F[0,x-b]|+|x+F|+(2C)(b+x,2H)}\\
& &\succeq 3|F|+F(0,x-b)+3|C|\sim 3|A|+F(0,x-b)\succ 3|A|.
\end{eqnarray*}
This ends the proof. \quad $\Box$(Lemma \ref{ftmore})

\begin{lemma}\label{ftbtmore}

Suppose there are $0\prec a<b\prec c\prec H$ such that $A=F\cup
B\cup C$, where $F$ is a forward triangle from $0$ to $a$, $B$ is
a backward triangle from $b$ to $c$, and $C\subseteq [c+1,H]$ with
$|C|\preceq\frac{1}{2}(H-c)$. Then $|2A|\succ 3|A|$.

\end{lemma}

\noindent {\bf Proof}:\quad Without loss of generality we can
assume $c,c-1\in B$. By Lemma \ref{smallbig} we have that
$|A[0,c]+A[0,c]|\succ 3A(0,c)$ implies $|2A|\succ 3|A|$.
So we can now assume $|A[0,c]+A[0,c]|\sim 3A(0,c)$.
Let $\bar{a}=\max F$ and $\bar{b}=\min B$. By Lemma \ref{ftbt2},
we have $\bar{a}\sim\frac{a}{2}$, $\bar{b}\sim\frac{b+c}{2}$,
$F$ is full in $[0,\bar{a}]$, and $B$ is full in $[\bar{b},c]$.

\medskip

{\bf Case \ref{ftbtmore}.1}\quad There is a $x\in C$ with $x\succ
c$ such that $C(c,x)\sim 0$.

We have
\begin{eqnarray*}
\lefteqn{|2A|\succeq 3|F|+3|B|+|x+B[2c-x,c]|}\\
& &\quad +|(\{c-1,c\}\cup C[x,H])+(\{c-1,c\}\cup C[x,H])|\\
& &\succeq 3|F|+3|B|+B(2c-x,c)+3|\{c-1,c\}\cup C[x,H])|\\
& &\sim 3|F|+3|B|+3|C|+B(2c-x,c)\succ 3|A|.
\end{eqnarray*}
\quad $\Box$(Case \ref{ftbtmore}.1)

\medskip

{\bf Case \ref{ftbtmore}.2}\quad For every $x\in C$, if $x\succ c$,
then $C(c,x)\succ 0$.

The assumption implies that for every $y\succ c$, there is a $x\in
C$ with $c\prec x\prec y$. Let $x\in C$ be such that
$0\prec x-c<\bar{a}$. Then
\begin{eqnarray*}
\lefteqn{|2A|\succeq a+B(b,c)+|c+F[0,x-c]|+|x+F[0,\bar{a}]|}\\
& &\quad +|[c+b,2c]|+|(\{c-1,c\}\cup C)+(\{c-1,c\}\cup C)|\\
& &\succeq 2|F|+|B|+F(0,x-c)+|F|+2|B|+3|C|\succ 3|A|.
\end{eqnarray*}
This ends the proof. \quad $\Box$(Lemma \ref{ftbtmore})

\begin{lemma}\label{ftthm}

Suppose there is an $a$ with $0\prec a\prec H$ such that $F=A[0,a]$
is a forward triangle from $0$ to $a$ and $A(a,H)\preceq\frac{1}{2}(H-a)$.
If $|2A|\sim 3|A|$, then $A$ is a full subset of a \bp\ of difference $3$ or 
a full subset of a \bp\ of difference $1$.

\end{lemma}

\noindent {\bf Proof}:\quad Note that if $A$ is a subset of a \bp\,
then $A$ must be a full subset of that \bp\ when $|2A|\sim 3|A|$.
Let $b=\min A[a+1,H]$. If $b\sim H$, then
\[|2A|\succeq a+(2A)(a+1,H)+|H+F|\sim 3|A|+(2A)(a+1,H).\]
Hence $(2A)(a+1,H)\sim 0$. By Lemma \ref{block}, $\bar{a}=\max
F\sim \frac{1}{2}(a+1)$. This shows $2\bar{a}\prec b$ and
$\bar{a}+H\prec 2b$. Hence $[0,\bar{a}]\cup [b,H]$ is the desired 
\bp\ of difference $1$. So we
can assume $b\prec H$. If $A(b,H)\sim 0$, then
\[|2A|\succeq a+|b+A[0,H-b]|+|H+A[0,a]|\sim 3|A|+A(0,H-b)\succ 3|A|.\]
Hence we can assume $A(b,H)\succ 0$.

Suppose $A$ is not a full subset of a \bp\ of difference $3$.
By Lemma \ref{ftap} we can assume $\gcd(A[b,H]-b)=1$.
If $A(b,H)\preceq\frac{1}{2}(H-b)$,
then the lemma follows from Lemma \ref{ftmore}.
So now we can assume that $A(b,H)\succ\frac{1}{2}(H-b)$.

By (2) of Lemma \ref{triangle} there are $a<c\prec b\prec c'
\leqslant H$ such that $A[c,c']$ is a backward triangle and
$A(c',H)\sim\frac{1}{2}(H-c')$. If $c'\prec H$, then
by Lemma \ref{ftbtmore} we have $|2A|\succ 3|A|$. Hence we can assume
$c'\sim H$. But now $A$ becomes the
union of a forward triangle $A[0,a]$ and a backward triangle
$A[c,H]$. Now the lemma follows from Lemma \ref{ftbt2}.
\quad $\Box$(Lemma \ref{ftthm})

\begin{lemma}\label{moreftmore}

Suppose there are $0\prec a\prec a'\prec H$ such that 
$A(0,a)\preceq\frac{1}{2}a$, $A(a',H)\preceq\frac{1}{2}(H-a')$,
$A[a,a']$ is a forward triangle from $a$ to $a'$, and
$A[a,H]$ is not a subset of a \bp\ of difference $3$. 
Then $|2A|\succ 3|A|$.

\end{lemma}

\noindent {\bf Proof}:\quad By Lemma \ref{smallbig} we can assume
$|A[a,H]+A[a,H]|\sim 3A(a,H)$. By Lemma \ref{ftthm} we can assume
that $A[a,H]$ is full in a \bp\ $[a,c]\cup [c',H]$ for some $c<c'$ 
in $[a,H]$. If $c'\prec H$, then the lemma follows from Lemma \ref{ftbtmore}.
So we can assume $c'\sim H$. Note that $c\sim\frac{a+a'}{2}$. Hence
we have $2c\prec a+H$. If there is a $x\prec a$ in $A$ with $x\geqslant
2c-H$ such that $A(x,a)\succ 0$, then
\begin{eqnarray*}
\lefteqn{|2A|\succeq (2A)(0,2a)+(2A)(2a,2c)}\\
& &\quad +(2A)(H+x,H+a)+|H+A[a,c]|\\
& &\succeq 3A(0,a)+3A(a,c)+A(x,a)\succ 3|A|.
\end{eqnarray*}
Otherwise choose $x\prec a$ in $A$ such that $A(x,a)\sim 0$. Without
loss of generality let $a,a+1\in A$. Then
\[(2A)(0,x+a)\succeq |(A[0,x]\cup\{a,a+1\})+
(A[0,x]\cup\{a,a+1\})|\succeq 3A(0,x)\sim 3A(0,a).\]
So we have
\begin{eqnarray*}
\lefteqn{|2A|\succeq (2A)(0,x+a)+|x+A[a,2a-x]|}\\
& &\quad +(2A)(2a,2c)+|H+A[a,c]|\\
& &\succeq 3A(0,a)+3A(a,c)+A(a,2a-x)\succ 3|A|.
\end{eqnarray*}
\quad $\Box$(Lemma \ref{moreftmore})

\begin{lemma}\label{moreftmore2}

Suppose there are $0\prec a\prec a'\prec H$ such that 
$A(0,a)\preceq\frac{1}{2}a$, $A(a',H)\preceq\frac{1}{2}(H-a')$,
$A[a,a']$ is a forward triangle from $a$ to $a'$, and
$A[a',H]$ is not a subset of an \ap\ of difference $3$. 
Then $|2A|\succ 3|A|$.

\end{lemma}

\noindent {\bf proof}:\quad By Lemma \ref{moreftmore} it suffices to 
show that $A[a,H]$ is not a subset of a \bp\ of difference $3$.
If $A[a,H]$ is a subset of a \bp\ of difference $3$, then $|2A|\sim 3|A|$
implies that $A[a,H]$ is a full subset of the \bp\, This implies that
$A[a,a']$ is a subset of the union of an \ap\ of difference $3$ of
length $\sim \frac{1}{3}(a'-a)$ and an \ap\ of difference $3$ of
length $\sim\frac{1}{6}(a'-a)$-both have the left-end points $\sim a$, 
and $A[a',H]$ is a full subset of an \ap\ of difference $3$,
which contradicts the assumption of the lemma.
\quad $\Box$(Lemma \ref{moreftmore2})

\medskip

Starting from the next lemma to the end of this section, the condition
(\ref{sparse}) is assumed.

\begin{lemma}\label{apmore}

Assume that $|A|\prec\frac{1}{2}H$ and $A$ is neither a
subset of an \ap\ of difference $>1$ nor a subset of a \bp\,
Suppose that there is a $x\succ 0$ such that
$A(0,x)\succ 0$ and $A[0,x]$ is a subset of an \ap\ of difference
$>1$. Then $|2A|\succ 3|A|$.

\end{lemma}

\noindent {\bf Proof}:\quad Let $a=\min\{y\in A:\gcd(A[0,y])=1\}$
and $c=\max A[0,a-1]$. Let $d=\gcd(A[0,c])$. Note that $d>1$ is
a standard natural number because $A(0,x)\succ 0$.

First, we can assume that $a\prec H$ by the following reason:
Suppose $a\sim H$. If there is a $b\in A$ such that $b\not\equiv
0\,(\mod d)$ and $b\not\equiv a\,(\mod d)$, then
\[|2A|\succeq |A[0,c]+A[0,c]|+|a+A[0,c]|+|b+A[0,c]|\sim 4|A|.\]
If for any $b\in A$, $b\equiv 0\,(\mod d)$ or $b\equiv a\,(\mod
d)$, then $A$ is a subset of a \bp\ unless $d=2$. Assume $d=2$.
Hence $A[0,c]$ is a set of even numbers and $a$ is odd. If
$A[a,H]$ is a set of odd numbers, then $A=A[0,c]\cup A[a,H]$ is a
subset of a \bp\, So we can assume that there is an even number
$b\in A$ with $b>a$. Clearly $b\sim H$.
Let $A_e$ be the set of all even number in $A$. Then
\[|2A|\succeq |2A_e|+|a+A_e|\succeq 3|A_e|\sim 3|A|.\]
Hence if $|2A|\sim 3|A|$, then $|2A_e|\sim 2|A_e|$.
By Theorem \ref{2k-1+b} and the fact that $b\sim H$ the set 
$A_e$ is full in the set of all even numbers in $[0,H]$,
which contradicts (\ref{sparse}). 

Second, we can assume that $A(a,H)\succ 0$ by the following reason:
Suppose $A(a,H)\sim 0$. If $|A[0,c]+A[0,c]|\succ 2A(0,c)$, then
\[|2A|\succeq |A[0,c]+A[0,c]|+|a+A[0,c]|\succ 3|A|.\]
So we can assume $|A[0,c]+A[0,c]|\sim 2A(0,c)$.
By Theorem \ref{2k-1+b} $A[0,c]$ is full.
This implies $A(a+c-H,c)\succ 0$. Hence we have
\begin{eqnarray*}
\lefteqn{|2A|\succeq |A[0,c]+A[0,c]|}\\
& &\quad +|a+A[0,c]|+|H+A[a+c-H,c]|\\
& &\succeq 3A(0,c)+A(a+c-H,c)\\
& &\sim 3|A|+A(a+c-H,c)\succ 3|A|.
\end{eqnarray*}

Now we are ready to prove the lemma. The proof is divided into
five cases.

\medskip

{\bf Case \ref{apmore}.1}\quad $d=2$ and
$A(a,H)\preceq\frac{1}{2}(H-a)$.

Clearly $a$ is odd. Since $A$ is not a subset of a \bp, 
then $\gcd(A[a,H]-a)=d'$ is not an even number.

Suppose $d'>2$. Let $c'=\max\{x\in A:\gcd(A[x,H]-x)=1\}$. Then
$c'\leqslant c$ and $A[c'+1,H]\subseteq (H-(d''\ast\,^*\!\nat))$
for some $d''>1$ and $d''|d'$.
By a symmetric argument of showing $A(a,H)\succ 0$ above,
we can assume $A(0,c')\succ 0$. With a little more effort we can
show that
\[(A[a,H]+A[a,H])\cap (c'+A[a+1,H])=\emptyset,\]
\[(A[0,c]+A[0,c])\cap (c'+A[a+1,H])=\emptyset,\mbox{ and }\]
\[(a+A[0,c])\cap (c'+A[a+1,H])=\emptyset.\]
The second equality above is due to the fact that if
$x_1,x_2\in A[0,c]$ and $a'\in A[a+1,H]$ having 
$x_1+x_2=c'+a'\geqslant c'+a+1$, then $x_1,x_2>c'$,
which implies $x_1,x_2\in (H-(d''\ast\,^*\!\nat))$. 
This implies $c'=x_1+x_2-a'\in
(H-(d''\ast\,^*\!\nat))$, which contradicts the choice
of $c'$. The reason for the third equality above is similar.
Hence
\begin{eqnarray*}
\lefteqn{|2A|\succeq |A[0,c]+A[0,c]|+|a+A[0,c]|}\\
& &\quad +|A[a,H]+A[a,H]|+|c'+A[a+1,H]|\\
& &\succeq 3A(0,c)+3A(a,H)\sim 3|A|.
\end{eqnarray*}
If $|2A|\sim 3|A|$, then
$|A[0,c]+A[0,c]|\sim 2A(0,c)$ and
$|A[a,H]+A[a,H]|\sim 2A(a,H)$. Hence $A[0,c]$ is full in the set
of all even numbers in $[0,c]$ and $A[a,H]$ is full in
$(a+(d'\ast\,^*\!\nat))\cap [a,H]$. Without loss of generality,
we can assume $c,c-2,c-4\in A$ and $c'=c$. Note that $c+A[a,H]$,
$c-2+A[a,H]$, and $c-4+A[a,H]$ are pairwise disjoint because
$d'$ is odd and $d'>2$. Hence we have
\begin{eqnarray*}
\lefteqn{|2A|\succeq |A[0,c]+A[0,c]|+|a+A[0,c]|}\\
& &\quad +|\{c,c-2,c-4\}+A[a,H]|+|H+A[a,H]|\\
& &\succeq 3A(0,c)+4A(a,H)\succ 3|A|.
\end{eqnarray*}
So $|2A|\succ 3|A|$ must be true.
This ends the proof of the case for $d'>2$.

Now assume that $\gcd(A[a,H]-a)=d'=1$. This implies
$|A[a,H]+A[a,H]|\succeq 3A(a,H)$.
Hence
\[|2A|\succeq (2A)(0,2c)+|a+A[0,c]|+(2A)(2a,2H)
\succeq 3|A|.\] 
We now derive a contradiction by assuming $|2A|\sim 3|A|$.
By the inequality above we have that
$A[0,c]$ is full in the set of all
even numbers in $[0,c]$. Suppose $c\prec a$.
If there is a $x\succ a$ in $A$ such that
$x-a<a-c$. Then we have
\begin{eqnarray*}
\lefteqn{|2A|\succeq 2A(0,c)+|a+A[0,c]|}\\
& &\quad +|x+A[a+c-x,c]|+|A[a,H]+A[a,H]|\\
& &\succeq 3A(0,c)+3A(a,H)+A(a+c-x,c)\succ 3|A|,
\end{eqnarray*}
which contradicts $|2A|\sim 3|A|$. Otherwise we
can find a $x\succ a$ in $A$ such that $A(a,x)\sim 0$. 
Let $F\subseteq A[a,H]$ be finite such that $a\in F$
and $\gcd((F\cup A[x,H])-a)=1$. Then 
\[(2A)(x+a,2H)\succeq |(F\cup A[x,H])+(F\cup A[x,H])|\succeq 3A(x,H)\sim 3A(a,H).\]
Hence we have
\[|2A|\succeq 3A(0,c)+3A(a,H)+|x+A[2c-x,c]|\succ 3|A|.\] 
So we can assume $c\sim a$. Recall
that we have $A(0,c)\succ 0$, $A(a,H)\succ 0$,
$\gcd(A[0,c])=2$, $A[0,c]$ is full, $\gcd(A[a,H]-a)=1$, 
and $A(a,H)\preceq\frac{1}{2}(H-a)$.
Note that since $A(0,c)\sim\frac{1}{2}(c+1)$, then
$|A|\prec\frac{1}{2}H$ implies
$A(a,H)\prec\frac{1}{2}(H-a)$. Since $A[0,c]$ is full,
we can, without loss of generality, assume that
$c,c-2,c-4\in A$.

\medskip

{\bf Subcase \ref{apmore}.1.1}\quad $\underline{d}_{a+U}(A)=0$.

Choose a $x\in A$ with $x\succ a$ such that
$A(a,x)<\frac{1}{8}(x-a+1)$. Let $F\subseteq A[a,H]$ be finite
such that $a\in F$ and $\gcd(F\cup A[x,H])=1$. Then
\begin{eqnarray*}
\lefteqn{|2A|\succeq |A[0,c]+A[0,c]|+|a+A[0,c]|+|x+A[a+c-x,c]|}\\
& &\quad +|(F\cup A[x,H])+(F\cup A[x,H])|\\
& &\succeq 3A(0,c)+\frac{1}{2}(x-a+1)+3A(x,H)\\
& &\succeq 3|A|+\frac{1}{2}(x-a+1)-3A(a,x)\succ 3|A|,
\end{eqnarray*}
which is again a contradiction.
\quad $\Box$(Subcase \ref{apmore}.1.1)

\medskip

{\bf Subcase \ref{apmore}.1.2}\quad
$\underline{d}_{a+U}(A)>\frac{1}{2}$.

By (3) of Lemma \ref{triangle} there exists a $b\succ a$ 
such that $A[a,b]$ is a forward triangle
from $a$ to $b$. Since $A(a,H)\prec\frac{1}{2}(H-a)$, then 
$A(b,H)\prec\frac{1}{2}(H-b)$ and $b\prec H$. If $|2A|\sim
3|A|$, then $|A[c-4,H]+A[c-4,H]|\sim 3A(c-4,H)$. Note that $A[c-4,H]$
is not a subset of a \bp\ of difference $3$ because $c,c-2,c-4\in A$.
Hence by Lemma \ref{ftthm}, $A[c-4,H]$ is a full subset
of a \bp\ $[c-4,a']\cup [b',H]$ for some $a',b'\in A$. If $b'\sim H$,
then by the fact that $2a'\prec a+H$ we have
\[|2A|\succeq 3A(0,c)+2a'-2a+|H+A[2a'-H,H]|\succeq 3|A|+A(2a'-H,a)\succ 3|A|.\]
If $b'\prec H$, then the lemma follows from Lemma \ref{ftbtmore}.
\quad $\Box$(Subcase \ref{apmore}.1.2)

\medskip

{\bf Subcase \ref{apmore}.1.3}\quad
$0<\underline{d}_{a+U}(A)\leqslant\frac{1}{2}$.

Suppose for any $x\succ a$ in $A$ we have $\gcd(A[x,H]-x)>1$.
Choose a $x\sim a$ in $A$ such that $\gcd(A[x,H]-x)=d'>1$. Since
$\gcd(A[a,H]-a)=1$, then $|A[a,H]+A[a,H]|\sim 3A(a,H)$ implies
that $A[x,H]$ is full.

If $d'=2$, then $|A|\sim\frac{1}{2}H$, which contradicts
the condition (\ref{sparse}).

Suppose $d'=4$. Let $c'=c$ and $c''=c-2$ when $x$ is odd, or let 
$c'\in\{c,c-2\}$ such that $c'+x\equiv 2x+2\,(\mod d')$ and 
$c''=a$ when $x$ is even. Then $c'+A[x,H]$, $c''+A[x,H]$, and
$A[x,H]+A[x,H]$ are pairwise disjoint. Hence
\begin{eqnarray*}
\lefteqn{|2A|\succeq |A[0,a]+A[0,a]|+|c'+A[x,H]|}\\
& &\quad +|c''+A[x,H]|+|A[x,H]+A[x,H]|\\
& &\succeq 3|A|+A(x,H)\succ 3|A|.
\end{eqnarray*}

Suppose $d'=3$ or $d'>4$. Then there are $c',c''$ in $\{c,c-2,c-4\}$
such that $c'+A[x,H]$, $c''+A[x,H]$, and $A[x,H]+A[x,H]$ are
pairwise disjoint. Hence $|2A|\succ 3|A|$ by the same reason
above.

Therefore, we can now assume that there is a $x\succ a$ in $A$
such that $\gcd(A[x,H]-x)=1$. Since $\gcd((A[c-4,H]-c-4)\cap U)=1$
and $(A[c-4,H]-c-4)\cap U$ is not a subset of a $U$--unbounded
\bp\ of difference $d>1$ because $a,c,c-2,c-4\in A$, then by 
Lemma \ref{kneser} there exists a $y\succ
a$ in $A$ with $c\prec y\leqslant x$ and
$A(y,H)\preceq\frac{1}{2}(H-y+1)$ such that
$(2A)(2(c-4),2y)\succ 3A(c-4,y)$. Hence by Lemma \ref{smallbig}
$|A[c-4,H]+A[c-4,H]|\succ 3A(c-4,H)$, which implies $|2A|\succ
3|A|$. This ends the proof.\quad $\Box$(Case \ref{apmore}.1)

\medskip

{\bf Case \ref{apmore}.2}\quad $d=2$ and
$A(a,H)\succ\frac{1}{2}(H-a)$.

By Lemma \ref{triangle} we can find $0\leqslant a'\prec a\prec a''
\leqslant H$ such that $A[a',a'']$
is a backward triangle from $a'$ to $a''$ and
$A(a'',H)\sim\frac{1}{2}(H-a'')$. Without loss of generality we
can assume $\gcd(A[a'',H]-a'')=1$. Then by Lemma \ref{ftap} we have
that $|A[0,a'']+A[0,a'']|\succ 3A(0,a'')$ or $A[0,a'']$ is a full 
subset of a \bp\ of difference $3$. However, the former
implies $|2A|\succ 3|A|$ by Lemma \ref{smallbig} and the latter
is impossible because $d=2$.  \quad $\Box$(Case \ref{apmore}.2)

\medskip

{\bf Case \ref{apmore}.3}\quad $d=3$ and
$A(a,H)\preceq\frac{1}{2}(H-a)$.

(Note that this case does not occur when $|A|\sim\frac{1}{2}H$.)
Since $A$ is not a subset of a \bp, we can define
\[b=\min\{x\in A:x\not\in\{0,a\}\,(\mod 3)\}.\] Let
$A_0=A\cap (3\ast\,^*\!\nat)$, $A_a=A\cap (a+(3\ast\,^*\!\nat))$,
and $A_b=A\cap (b+(3\ast\,^*\!\nat))$. Let $l_0,l_a,l_b$ be the
least element of $A_0,A_a,A_b$, respectively. Let $u_0,u_a,u_b$ be
the largest element of $A_0,A_a,A_b$, respectively.
Note that the rest of the proof does not use the fact that $a\succ 0$.

\medskip

{\bf Subcase \ref{apmore}.3.1}\quad $b\sim H$.

We have $|A|\sim |A_0|+|A_a|$. We can also assume
$|A_a|\succ 0$ because otherwise
\[|2A|\succeq |2A_0|+|a+A_0|+|b+A_0|=4|A_0|=4|A|.\]
Since $A_0\cup A_a$ is a subset of
a \bp, then by Theorem \ref{2k-1+b}, $A_0$ is full and $A_a$ is full. 
This implies $u_a\prec H$ or $u_0\prec H$ because
$A(a,H)\preceq\frac{1}{2}(H-a)$.
Suppose $u_a\prec H$ and $u_a\leqslant u_0$. Then
\begin{eqnarray*}
\lefteqn{|2A|\succeq|2A_0|+|2A_a|}\\
& &\quad +|A_0+A_a|+|b+A_0[u_a+u_0-b,u_0]|\\
& &\succeq 3|A|+A_0(u_a+u_0-b,u_0)\succ 3|A|.
\end{eqnarray*}
By the same reason, if $u_0\prec H$ and $u_0\leqslant u_a$, then
$|2A|\succ 3|A|$. Note that if both $u_0\prec H$ and $u_a\prec H$
are true, then either $u_0\leqslant u_a$ or $u_a\leqslant u_0$.
\quad $\Box$(Subcase \ref{apmore}.3.1)

\medskip

{\bf Subcase \ref{apmore}.3.2}\quad $b\prec H$.

Suppose $d'=\gcd(A[b,H]-b)>1$. 
If $d'=2$, then the proof of this case is same as the proof
in Case \ref{apmore}.1 and Case \ref{apmore}.2 by considering
$H-A$ in the place of $A$. So we can assume that $d'>2$.

If $d'=3$, then $u_0,u_a<b$. Note that $b\not\in\{0,a\}\,(\mod 3)$. 
We can assume $|A_a|\succ 0$ because if $|A_a|\sim 0$, then
\[|2A|\succeq |2A_0|+|2A_b|+|A_0+A_b|+|u_a+A_0|\succ 3|A|.\]
We can also assume $|A_b|\succ 0$ because otherwise let
$c\leqslant b$ such that $A(c,b)\sim 0$ and for every $x\prec c$,
$A(x,c)\succ 0$. Then
\[|2A|\succeq 3|A_0|+3|A_a|+|H+A[c+b-H,c]|\succ 3|A|.\]
Let $u=\max\{u_0,u_a\}$. If $|2A|\sim 3|A|$, then
\begin{eqnarray*}
\lefteqn{|2A|\succeq |2A_0|+|2A_a|}\\
& &\quad +|A_0+A_a|+|2A_b|+|u+A_b|\\
& &\succeq 3|A_0|+3|A_a|+3|A_b|=3|A|
\end{eqnarray*}
implies that $A_0$, $A_a$, and $A_b$ are full. If $u_0\prec b$ and
$u_0<u_a$, then
\begin{eqnarray*}
\lefteqn{|2A|\succeq |2A_0|+|2A_a|+|A_0+A_a|+|2A_b|}\\
& &\quad +|b+A_a[u_a+u_0-b,u_a]|+|u_a+A_b|\\
& &\succeq 3|A_0|+3|A_a|+3|A_b|+A_a(u_a+u_0-b,u_a)\succ 3|A|.
\end{eqnarray*}
So we can assume $u_0\sim b$. By a similar argument we can also
assume $u_a\sim b$. However, above assumptions imply that
\begin{eqnarray*}
\lefteqn{|2A|\succeq |2A_0|+|2A_a|+|A_0+A_a|}\\
& &\quad +|2A_b|+|u_a+A_b|+|u_0+A_b|\\
& &\succeq 3|A_0|+3|A_a|+4|A_b|\succ 3|A|.
\end{eqnarray*}

Suppose $d'\geqslant 4$. We re-define $A_0$ to be $A_0[0,b-1]$,
$A_a$ to be $A_a[0,b-1]$, $u_0=\max A_0$, and $u_a=\max A_a$.
Let $u=\max\{u_0,u_a\}$. Then
\[|2A|\succeq |2A_0|+|2A_a|+|A_0+A_a|+|2A_b|+|u+A_b|\succeq 3|A|\]
together with $|2A|\sim 3|A|$ imply that $A_0$, $A_a$, and $A_b$
are all full. Note that $|A_0|\succ 0$ is always true. We can also 
assume $|A_a|\succ 0$ because otherwise we have
\[|2A|\succeq |2A_0|+|a+A_0|+|b+A_0|+|u+A_b|+|2A_b|\succeq
4|A_0|+3|A_b|\succ 3|A|.\]
Hence we can assume $u,u-3,u-6\in A_0\cup A_a$. Since there are
$u',u''\in\{u,u-3,u-6\}$ such that $u'+A_b$, $u''+A_b$, and
$2A_b$ are pairwise disjoint, we have
\begin{eqnarray*}
\lefteqn{|2A|\succeq |2A_0|+|2A_a|+|A_0+A_a|}\\
& &\quad +|u'+A_b|+|u''+A_b|+|2A_b|\\
& &\succeq 3|A_0|+3|A_a|+4|A_b|\succ 3|A|.
\end{eqnarray*}

Therefore, we can now assume that $d'=1$.
If $A(b,H)\succ\frac{1}{2}(H-b)$, then by Lemma \ref{triangle}
there exist $b'\prec b\prec b''\leqslant H$ such that
$A(b'',H)\sim\frac{1}{2}(H-b'')$ and $A[b',b'']$ is a backward triangle.
Since $0,a,b\in A[0,b'']$, then $A[0,b'']$ is not a subset of a \bp\
of difference $3$. Clearly $A[0,b'']$ is not a subset of a \bp\ of
difference $1$ because $d>1$. Hence we have $|2A|\succ 3|A|$ by
Lemma \ref{ftthm} and Lemma \ref{smallbig}. 
So we can now assume that $A(b,H)\preceq\frac{1}{2}(H-b)$.
let's re-define $A_0$ to be $A_0[0,b-1]$,
$A_a$ to be $A_a[0,b-1]$, $u_0=\max A_0$, and $u_a=\max A_a$.
Then by Lemma \ref{bpleft} we have that $A_0$ and $A_a$ are full and
$u_0,u_a\sim b$. We can also assume $A_a(l_a,u_a)\succ 0$ because
otherwise $(2A)(0,2b)\succeq 4A(0,b)$, which implies $|2A|\succ 3|A|$.

If $A(b,H)\sim\frac{1}{2}(H-b)$, then $A(0,b)\prec\frac{1}{2}b$. Since
$u_0\sim b$, $u_a\sim b$, and $d=3$, then $u_a-l_a\prec\frac{1}{2}b$,
which implies $l_a\succ\frac{1}{2}b$. Hence
\begin{eqnarray*}
\lefteqn{|2A|\succeq |2A_0|+|2A_a|+|A_0+A_a|}\\
& &\quad +|A[b,H]+A[b,H]|+|b+A_0[0,2l_a-b]|\\
& &\succeq 3A(0,b)+3A(b,H)+A_0(0,2l_a-b)\succ 3|A|.
\end{eqnarray*}
So we can assume $A(b,H)\prec\frac{1}{2}(H-b)$.

If $\underline{d}_{b+U}(A)=0$, then there is a $x\in A$, $x\succ b$
such that either $x-b<u_0$ and
$A(b,x)\preceq\frac{1}{10}(x-b+1)$, or $A(b,x)\sim 0$. Let
$F\subseteq A[b,H]$ be a finite set such that $b\in F$ and
$\gcd((F\cup A[x,H])-b)=1$. Then
\begin{eqnarray*}
\lefteqn{|2A|\succeq |2A_0|+|2A_a|+|A_0+A_a|}\\
& &\quad +|(F\cup A[x,H])+(F\cup A[x,H])|+|x+A_0[2u_0-x,u_0]|\\
& &\succeq 3A(0,b)+3A(x,H)+A_0(2u_0-x,u_0)\\
& &\succeq 3A(0,b)+3A(x,H)+\frac{1}{3}(x-u_0+1)\\
& &\succeq 3|A|+\frac{1}{3}(x-u_0+1)-\frac{3}{10}(x-b+1)\succ 3|A|.
\end{eqnarray*}

If $\underline{d}_{b+U}(A)>\frac{1}{2}$, then there is a $x\succ
b$ such that $A[b,x]$ is a forward triangle. 
Clearly $x\prec H$ and $A(x,H)\prec\frac{1}{2}(H-x)$. 
Let $u'=\min\{u_0,u_a\}$. Note that $u'\sim b$.
By Lemma \ref{ftthm} and Lemma \ref{smallbig}, $|2A|\sim 3|A|$ implies
that $A[u',H]$ is either a full subset of a \bp\ of difference $3$ or
a full subset of a \bp\ of difference $1$. Since $u_0,u_a,b\in A[u',H]$,
$A[u',H]$ cannot be a subset of a \bp\ of difference $3$.
Let $A[b,H]$ be a full subset of the \bp\ $[b,z]\cup [z',H]$ for some 
$z<z'$ in $A[b,H]$. Note that $2z\prec b+z'$. Then
\begin{eqnarray*}
\lefteqn{|2A|\succeq |2A_0|+|2A_a|+|A_0+A_a|+|A[b,z]+A[b,z]|}\\
& &\quad +|A[b,z]+A[z',H]|+|A[z',H]+A[z',H]|+|z'+A[2z-z',b]|\\
& &\succeq 3A(0,b)+3A(b,H)+A(2z-z',b)\succ 3|A|.
\end{eqnarray*}

Now we can assume $0<\underline{d}_{b+U}(A)\leqslant\frac{1}{2}$. 
Suppose there is a $b'\sim b$ in $A$ such that $\gcd(A[b',H]-b')=d''>1$, 
If $d''=2$, then there is a $b'''\sim b$ such that $d'''-d''$ is odd.
Hence $|A[b''',H]+A[b''',H]|\sim 3A(b''',H)$ implies that $A[b'',H]$
is full, which contradicts $A(b,H)\prec\frac{1}{2}(H-b)$.
If $d''\geqslant 3$, then $|A[u',H]+A[u',H]|\succeq 4A(u',H)\succ 3A(u',H)$,
which contradicts $|2A|\sim 3|A|$ by Lemma \ref{smallbig}.
So we can assume that there is an
$x\succ b$ in $A$ such that $\gcd(A[x,H]-x)=1$. 
Since $A_0$ and $A_a$ are full, we can assume $u'-3\in A$. 
Hence $A\cap (u'-3+U)$ is
neither a subset of an \ap\ of difference $>1$ nor a subset
of a $(u'-3+U)$--unbounded \bp\,
By Lemma \ref{kneser} there exists a $y\in A$ with $b\prec y<x$ 
such that $A(y,H)\preceq\frac{1}{2}(H-y)$ and $(2A)(2(u'-2),2y)\succ
3A(u'-2,y)$. By Lemma \ref{smallbig} we have
$|A[b,H]+A[b,H]|\succ 3A(b,H)$, which implies $|2A|\succ 3|A|$
again by Lemma \ref{smallbig}.
\quad $\Box$(Case \ref{apmore}.3)

\medskip

{\bf Case \ref{apmore}.4}\quad $d=3$ and $A(a,H)\succ\frac{1}{2}(H-a)$.

By Lemma \ref{ftbt} and (\ref{sparse}) we can find 
$0\prec a'\prec a\prec a''\leqslant H$
such that $A[a',a'']$ is a backward triangle and 
$A(a'',H)\sim\frac{1}{2}(H-a'')$. By (\ref{sparse}) we have
$A(0,a'')\prec\frac{1}{2}a''$. If $a''\sim H$, then the lemma follows
from Lemma \ref{ftthm}. So we can assume $a''\prec H$. Now the lemma
follows from Lemma \ref{moreftmore} unless $A[0,a'']$ is a subset
of a \bp\ of difference $3$. Without loss of generality, let
$A[0,a'']$ be a subset of a \bp\ $I_0\cup I_1$ of difference $3$, 
where $I_i=i+(3\ast\,^*\!\nat)$ for $i\in [0,2]$. 
Let $A_i=A\cap I_i$ and $b=\min (A\cap I_2)$. Then $b\succeq a''$. 
Let $l_1=\min A_1[0,a'']$. Then $2l_1\succ a''$. 
If $b\sim a''$, then
\[(2A)(0,2a'')\succeq 3A(0,a'')+|b+A_0[0,2l_1-b]|\succ 3A(0,a''),\]
which implies $|2A|\succ 3|A|$ by Lemma \ref{smallbig}.
So we can assume $b\succ a''$.

\medskip

{\bf Subcase \ref{apmore}.4.1}\quad $A(a'',b)\prec\frac{1}{2}(b-a'')$.

Then we have $A(b,H)\succ\frac{1}{2}(H-b)$. Hence we can find 
$a''\leqslant c\prec b\prec c'\leqslant H$ such that $A[c,c']$ is a
backward triangle from $c$ to $c'$ and $A(c',H)\sim\frac{1}{2}(H-c')$.
Since $A[0,c']$ contains two backward triangles, it cannot be
a subset of a \bp\ of difference $d'$ for $d'=1$ or $d'=3$. Hence
by Lemma \ref{ftthm} we have $|A[0,c']+A[0,c']|\succ 3A(0,c')$, which
implies $|2A|\succ 3|A|$. \quad $\Box$(Subcase \ref{apmore}.4.1)

\medskip

{\bf Subcase \ref{apmore}.4.2}\quad $A(a'',b)\succeq\frac{1}{2}(b-a'')$.

Let $c=\max\{x\in [a'',b-1]:x,x-1\in A\}$. It is easy to see that
$c\succ a''$ and $A(c,b)\preceq\frac{1}{3}(b-c)$. Since $|2A|\sim 3|A|$
implies $|A[0,c]+A[0,c]|\sim 3A(0,c)$, then we can assume that 
$A[0,c]$ is full in the \bp\ $I_0\cup I_1$, which implies 
$A(a'',c)\sim\frac{2}{3}(c-a'')$. Hence $A(c,H)\prec\frac{1}{2}(H-c)$.
So we can find a $m$ with $a''\prec m\prec c$ such that 
$A(m,H)\prec\frac{1}{2}(H-m)$. It is easy to show that there is
a $m'\prec H$ such that $A[m,m']$ is a forward triangle.
Since $A[m,H]$ cannot be a full subset of a \bp\ of difference $3$,
because it contains $b$, or $1$ because $A[m,b-1]\prec b-m$, then
by Lemma \ref{ftthm} we have $|A[m,H]+A[m,H]|\succ 3A(m,H)$.
Since $A[0,m]$ is a subset of a \bp\ of difference $3$ we have
$|A[0,m]+A[0,m]|\succeq 3A(0,m)$. By Lemma \ref{smallbig} we have
$|2A|\succ 3|A|$. 
\quad $\Box$(Case \ref{apmore}.4)

\medskip

{\bf Case \ref{apmore}.5}\quad $d\geqslant 4$.

Since $A$ is not a subset of a \bp, the number
$b=\min\{x\in A:x\not\in\{0,a\}\,(\mod d)\}$ is well defined.
Let $I_i=i+(d\ast\,^*\!\nat)$ and $A_i=A\cap I_i$. Let
$u_i=\max A_i[0,b-1]$ and $l_i=\min A_i[0,b-1]$ for $i=0,a$.

If $b\equiv 2a\,(\mod d)$, then $a+b\not\equiv 0\, (\mod d)$
because otherwise $A[0,b]$ is a subset of an \ap\ with difference
$\frac{d}{3}>1$. Hence we have either $b\not\equiv 2a\,(\mod d)$ or
$a+b\not\equiv 0\, (\mod d)$. This implies
\[(b+A_0[0,u_0])\cap (A[0,b-1]+A[0,b-1])=\emptyset\]
or \[(b+A_a[l_a,u_a])\cap (A[0,b-1]+A[0,b-1])=\emptyset.\] If
$A(b,H)\succ\frac{1}{2}(H-b)$, then there are 
$0\prec b'\prec b\prec b''\leqslant H$ such that 
$A(b'',H)\sim\frac{1}{2}(H-b'')$ and $A[b',b'']$ is a backward
triangle. If $|2A|\sim 3|A|$, then $A[0,b'']$ is a full subset
of either a \bp\ of difference $3$ or a \bp\ of difference $1$.
But both contradict $d\geqslant 4$.
So we can assume $A(b,H)\preceq\frac{1}{2}(H-b)$. 

Suppose $\gcd(A[b,H]-b)=1$. By Lemma \ref{bpleft}
we can show that $u_a\sim b$, $u_0\sim b$, $A_a(l_a,u_a)\succ 0$,
$A_0(0,u_0)\succ 0$, $A_0[0,u_0]$ is full in $[0,u_0]$, and
$A_a[l_a,u_a]$ is full in $[l_a,u_a]$. Hence
\begin{eqnarray*}
\lefteqn{|2A|\succeq |A_0[0,u_0]+A_0[0,u_0]|+|A_a[l_a,u_a]+A_a[l_a,u_a]|}\\
& &\quad +|A_0[0,u_0]+A_a[l_a,u_a]|+|A[b,H]+A[b,H]|\\
& &\quad +\min\{|b+A_0[0,u_0]|,|b+A_a[l_a,u_a]|\}\\
& &\succeq 3A(0,b)+3A(b,H)+\min\{A_0(0,u_0),A_a(l_a,u_a)\}\succ 3|A|.
\end{eqnarray*}

Suppose $\gcd(A[b,H]-b)=d'>1$. Let $c'=\max\{x\in A:\gcd(A[x,H]-x)<d'\}$.
Then we have \[(c'+A[b,H])\cap (A[b,H]+A[b,H])=\emptyset,\]
\[(c'+A[b,H])\cap (A[0,b-1]+A[0,b-1])=\emptyset,\,\mbox{ and}\]
\[(c'+A[b+1,H])\cap (b+A[0,b-1])=\emptyset.\] If $A_a(l_a,u_a)\succ 0$, 
then we have
\begin{eqnarray*}
\lefteqn{|2A|\succeq |A_0[0,u_0]+A_0[0,u_0]|+|A_0[0,u_0]+A_a[l_a,u_a]|}\\
& &\quad +|A_a[l_a,u_a]+A_a[l_a,u_a]|+\min\{|b+A_0[0,u_0]|,|a+A_a[l_a,u_a]|\}\\
& &\quad +|A[b,H]+A[b,H]|+|c'+A[b,H]|\succ 3|A|.
\end{eqnarray*}
If $A_a(l_a,u_a)\sim 0$, then we have 
\begin{eqnarray*}
\lefteqn{|2A|\succeq |A_0[0,u_0]+A_0[0,u_0]|+|a+A_0[0,u_0]|}\\
& &\quad +|b+A_0[0,u_0]|+|A[b,H]+A[b,H]|+|c'+A[b,H]|\\
& &\succeq 4A_0(0,u_0)+3A(b,H)\succ 3|A|.
\end{eqnarray*}
This ends the proof of Case \ref{apmore}.5 as well as the
lemma.\quad $\Box$(Lemma \ref{apmore})

\begin{lemma}\label{bpmore}

Assume $|A|\prec\frac{1}{2}H$,
$0<\underline{d}_U(A)\leqslant\frac{1}{2}$, and $A\cap U$ is a subset
of a $U$--unbounded \bp\ of difference $d$. If $A$ is not a subset of
a \bp, then $|2A|\succ 3|A|$.

\end{lemma}

\noindent {\bf Proof}:\quad Suppose $A\cap U$ is a subset of 
the $U$--unbounded \bp\ $(d\ast U)\cup(a+(d\ast U))$ for some 
$a\in A\cap U$. Clearly $d\not=2$. 
Let $b=\min\{x\in A:x\not\in\{0,a\}\,(\mod d)\}$. 
Note that $A(0,b-1)\succ 0$. By Lemma \ref{apmore} we can assume 
$\gcd(A[0,b-1])=1$. If $A(b,H)\succ\frac{1}{2}(H-b)$, 
then we can find $0\prec b'\prec b\prec b''\leqslant H$
such that $A(b'',H)\sim\frac{1}{2}(H-b'')$ and $A[b',b'']$ is a backward
triangle from $b'$ to $b''$. Note that $A[0,b'']$ cannot be a subset
of a \bp\ of difference $3$ because otherwise $A\cap U$ is a subset
of an \ap\ of difference $3$. By Lemma \ref{smallbig} we can assume
$|A[0,b'']+A[0,b'']|\sim 3A(0,b'')$. By Lemma \ref{ftthm} $A[0,b'']$
is a full subset of a \bp\ $[0,c]\cup [c',b'']$, which contradicts
the assumption $0<\underline{d}_U(A)\leqslant\frac{1}{2}$.
So we can assume $A[b,H]\preceq\frac{1}{2}(H-b)$.
If $d>3$, then the proof of the lemma is the same as the proof of 
Case \ref{apmore}.5. Suppose $d=3$. If $b\prec H$ and $\gcd(A[b,H]-b)>1$, 
then the lemma follows from Lemma \ref{apmore}. 
If $\gcd(A[b,H]-b)=1$ or $b\sim H$, 
then $|2A|\sim 3|A|$ implies $|A[0,b-1]+A[0,b-1]|\sim 3A(0,b-1)$,
which implies $A[0,b-1]$ is a full subset of a \bp\ of difference $3$.
Since $A\cap U$ is already a subset of the \bp, then
$\underline{d}_U(A)=\frac{2}{3}$, which contradicts
$\underline{d}_U(A)\leqslant\frac{1}{2}$.
\quad $\Box$(Lemma \ref{bpmore})

\medskip

Now we summarize all the proofs in this section into a theorem,
which takes care of the case in Theorem \ref{nsamain} under the condition
(\ref{sparse}).

\begin{theorem}\label{lesshalf}

Assume $A\subseteq [0,H]$ and $0,H\in A$. Suppose $\gcd(A)=1$ and
$0\prec |A|\prec\frac{1}{2}H$. If $A$ is not a subset of a \bp,
then $|2A|\succ 3|A|$.

\end{theorem}

\noindent {\bf Proof}:\quad By Lemma \ref{apmore} we can assume
that for every $x\succ 0$, if $A(0,x)\succ 0$, then
$\gcd(A[0,x])=1$ for every $x\succ 0$. If there is a $x\prec
H$ in $A$ such that $A(x,H)\succ 0$ and $\gcd(A[x,H]-x)>1$,
then the theorem is true
again by Lemma \ref{apmore} with $A$ replaced by $H-A$.
So we can assume that for every $x\prec H$ in $A$, if
$A(x,H)\succ 0$, then $\gcd(A[x,H]-x)=1$. We now
divide the proof into four cases according to the value of
$\underline{d}_{U}(A)$.

\medskip

{\bf Case \ref{lesshalf}.1}\quad
$\underline{d}_{U}(A)>\frac{1}{2}$.

Then there is a $c\succ 0$ such that $A[0,c]$ is a forward
triangle from $0$ to $c$. Since $|A|\prec\frac{1}{2}(H+1)$,
then $c\prec H$. Now the theorem follows from Lemma
\ref{ftthm}.

\medskip

{\bf Case \ref{lesshalf}.2}\quad $0<\underline{d}_{U}(A)\leqslant
\frac{1}{2}$.

If $A\cap U$ is a subset of an \ap\ of difference $>1$, then
the theorem follows from Lemma \ref{apmore}. If $A\cap U$ is a
subset of a $U$--unbounded \bp, then the theorem follows from
Lemma \ref{bpmore}. Otherwise by Lemma \ref{kneser} we can
find a $y\in A$ with $0\prec y\prec H$ such
that $A(y,H)\preceq\frac{1}{2}(H-y)$ and $(2A)(0,2y)\succ
3A(0,y)$. If $A(y,H)\sim 0$, then the theorem is already true
because $|A|\sim A(0,y)$. So we can assume $A(y,H)\succ 0$.
If $\gcd(A[y,H]-y)>1$, then the theorem
follows from Lemma \ref{apmore}. If $\gcd(A[y,H]-y)=1$,
then the theorem follows from Lemma \ref{smallbig}.
\quad $\Box$(Case \ref{lesshalf}.2)

\medskip

{\bf Case \ref{lesshalf}.3}\quad $\underline{d}_{U}(A)=0$ and
there is a $x\succ 0$ such that $A(0,x)\sim 0$.

By Lemma \ref{triangle} we can find such $x\in A$ such that
for any $y\succ x$, $A(x,y)\succ 0$.

If $A(x,H)\succ\frac{1}{2}(H-x)$, then we can find
$0\prec c'\prec x\prec c\leqslant H$ such that, 
$A(c,H)\sim\frac{1}{2}(H-c)$, and $A[c',c]$ is a
backward triangle. If $c\sim H$, then the theorem follows from
Lemma \ref{ftthm}. Suppose $c\prec H$. Note that $A[0,c]$ cannot
be a full subset of a \bp\ of difference $3$ by the 
condition of the case. Hence the lemma follows from 
Lemma \ref{moreftmore}.

If $A(x,H)\preceq\frac{1}{2}(H-x)$ and $A[x,H]$ is a subset of
an \ap\ of difference $>1$, then the theorem follows from Lemma \ref{apmore}.
If $A(x,H)\preceq\frac{1}{2}(H-x)$ and $A[x,H]$ is not a subset of
an \ap\ of difference $>1$, then
\[|2A|\succeq A(x,2x)+|A[x,H]+A[x,H]|\succeq 3|A|+A(x,2x)\succ 3|A|.\]
\quad $\Box$(Case \ref{lesshalf}.3)

\medskip

{\bf Case \ref{lesshalf}.4}\quad $\underline{d}_{U}(A)=0$ and for
every $x\succ 0$, $A(0,x)\succ 0$.

By symmetry we can also assume $\underline{d}_{H-U}(A)=0$ 
and for every $y\prec H$, $A(y,H)\succ 0$.

Let $|A|\sim\alpha H$. Then $0<\alpha<\frac{1}{2}$.
By Lemma \ref{triangle} there is a $b\succ 0$ in $A$ such that
$A(0,b)\sim\alpha b$ and $A(0,x)\prec\alpha x$ for every
$0\prec x\prec b$. By the assumption of this case, we have
$A(0,b)\succ 0$ and $A(b,H)\succ 0$. If there is a $0\prec x\prec H$
such that $A[0,x]$ or $A[x,H]$ is a subset of an \ap\ of
difference $>1$, then the theorem follows from Lemma \ref{apmore}.
Note that $\underline{d}_{b-U}(A)\geqslant\alpha$ by the choice of 
$b$. By Lemma \ref{smallbig} we can assume $|A[0,b]+A[0,b]|\sim 3A(0,b)$.
By Case \ref{lesshalf}.1 and Case \ref{lesshalf}.2 for $A[0,b]$
we can assume that $A[0,b]$ is a subset of a \bp\ of difference $d$.
Clearly $A[0,b]$ is a full subset of the \bp\,
If $d=1$, then $A[0,b]$ is a full subset of
$[0,x]\cup [x',b]$, which implies either $A(0,x')\sim 0$ or
$\underline{d}_U(A)=1$. Each of them contradicts the assumption of
the case. If $d>1$, then
$\underline{d}_U(A)=\frac{2}{d}$, which is again a contradiction
to the assumption of the case. \quad $\Box$(Theorem \ref{lesshalf})

\section{Second Step: When $\frac{|A|}{H}$ is almost $\frac{1}{2}$.}

In this section we again assume $A\subseteq [0,H]$, $0,H\in A$,
and $\gcd(A)=1$. In addition we also assume
\begin{equation}\label{half}
|A|\sim\frac{1}{2}H,
\end{equation}
\begin{equation}\label{notbp}
A\mbox{ is not a subset of a \bp,}
\end{equation}
and 
\begin{equation}\label{cond3k-3+b}
|2A|=3|A|-3+b
\end{equation}
for $0\leqslant b\sim 0$. (\ref{cond3k-3+b}) implies $|2A|\sim 3|A|$. 
Under the condition above we want to prove
\begin{equation}\label{shortap}
H+1\leqslant 2|A|-1+2b.
\end{equation}
Without loss of generality, we can assume that 
\begin{equation}\label{smallA}
|A|\leqslant\frac{1}{2}(H+1)
\end{equation}
because otherwise (\ref{shortap}) is trivially true.
In this section the letter $b$ is reserved only for the
purpose in (\ref{cond3k-3+b}). 

\begin{lemma}\label{addone}

Let $z\in [0,H]\smallsetminus A$ and let $A'=A\cup\{z\}$. Suppose
$|(2A')\smallsetminus (2A)|\leqslant 2$ and $|2A'|=3|A'|-3+b'$. If 
\begin{equation}\label{longap}
H+1>2|A|-1+2b,
\end{equation} 
then $0\leqslant b'\leqslant b-1$,
$|A'|\leqslant\frac{1}{2}(H+1)$, and $H+1>2|A'|-1+2b'$.

\end{lemma}

\noindent {\bf Proof}:\quad If $b=0$, then $|2A|=3|A|-3$.
By Theorem \ref{3k-3} we have $H+1\leqslant 2|A|-1$, which
contradicts $|A|\leqslant\frac{1}{2}(H+1)$. So we can assume
$b>0$. By the assumption of the lemma we have $H+1\geqslant
2|A|+2b$. Hence $|A'|=|A|+1\leqslant\frac{1}{2}(H+1)-b+1
\leqslant\frac{1}{2}(H+1)$. Since
\begin{eqnarray*}
\lefteqn{|2A'|=3|A'|-3+b'}\\
& &\leqslant |2A|+2=3|A|-1+b\\
& &=3|A'|-3+(b-1),
\end{eqnarray*}
then $b'\leqslant b-1$. If $b'<0$, then by Theorem \ref{2k-1+b}
$A'$ is a subset of an \ap\ of length $\leqslant 2|A'|-3=2|A|-1$,
which implies $H+1\leqslant 2|A|-1$, a contradiction to (\ref{longap}). 
Hence $b'\geqslant 0$. Finally $H+1>2|A|-1+2b=2|A'|-1+2(b-1)\geqslant
2|A'|-1+2b'$.\quad $\Box$(Lemma \ref{addone})

\begin{lemma}\label{bplittle}

If there is an $a\sim 0$ such that $A[a+1,H]$ is a subset of a
\bp\ of difference $3$, then $H+1\leqslant 2|A|-1+2b$.

\end{lemma}

\noindent {\bf Proof}:\quad Without loss of generality we can
assume $a\in A$ and $A[a,H]$ is not a subset of a \bp\ of difference
$3$. Fix $j\in [-2,0]$ such that $A[a+1,H]\subseteq A_0\cup A_1$
where $A_i=A\cap J_i$ and $J_i=j+i+(3\ast\,^*\!\nat)$ for $i=0,1,2$. 
For $i=0,1,2$ let $l_i=\min A_i$ and $u_i=\max A_i$. For $i=0,1$ 
let $I_i=J_i[0,u_i]$ and let $I_2=J_2[l_2,u_2]$. Clearly $a=u_2\sim 0$.
Since \[|2A|\succeq |2A_0|+|2A_1|+|A_0+A_1|\succeq 3|A_0|+3|A_1|
\sim 3|A|,\] then we have $|2A_i|\sim 2|A_i|$ for $i=0,1$. By
Theorem \ref{2k-1+b} we have that $A_i$ is full for $i=0,1$.
Note that $|A_i|\succ 0$ for $i=0,1$ by (\ref{half}).
If $l_0\succ 0$ and $l_0\geqslant l_1$, then $|2A|\succeq
3|A|+|a+A_1[l_1,2l_0]|\succ 3|A|$. By symmetry we can also 
prove that $l_1\succ 0$ and $l_1\geqslant l_0$ together are
impossible. So we can assume
$l_0\sim 0$ and $l_1\sim 0$. Without loss of generality we assume
$u_0=H$. Then $|A|\sim\frac{1}{2}(H+1)$ implies $u_1\sim\frac{H}{2}$.

Suppose the lemma is not true. Then we can assume that
(\ref{cond3k-3+b}), (\ref{smallA}), and (\ref{longap}) are true.
Without loss of generality we can assume that $|A|$ is the maximum
among all the sets in $I_0\cup I_1\cup I_2$ containing the original 
set and satisfying (\ref{cond3k-3+b}), (\ref{smallA}), and (\ref{longap}).

\medskip

{\bf Claim \ref{bplittle}.1}:\quad If $l_i\prec z\prec u_i$ and 
$z\equiv j+i\,(\mod 3)$ for $i=0$ or $i=1$ , then $z\in A$.

\medskip

Proof of Claim \ref{bplittle}.1:\quad Suppose not and let $A'=A\cup\{z\}$.
By Lemma \ref{addone} and by the maximality of $|A|$ we need only to show
that $|(2A')\smallsetminus (2A)|\leqslant 2$ for a contradiction. 
First let $z\equiv j\,(\mod 3)$. Let $y\in A'$.

If $y\in A_0\cup\{z\}$ and $y\prec u_0$, then
$A_0[y+1,y+t]\cap (y+z-A_0[z-t,z-1])\not=\emptyset$ for some
$0\prec t\prec\min\{u_0-y,z\}$, which implies $y+z\in (2A_0)\subseteq (2A)$,
by the pigeonhole principle. If $y\in A_0\cup\{z\}$ and $y\sim u_0$, then
$A_0[y-t,y-1]\cap (y+z-A_0[z+1,z+t])\not=\emptyset$ for some
$0\prec t\prec\min\{y,u_0-z\}$, which implies $y+z\in (2A_0)\subseteq (2A)$.
If $y\in A_1$ and $y\prec u_1$, then
$A_1[y+1,y+t]\cap (y+z-A_0[z-t,z-1])\not=\emptyset$ for some
$0\prec t\prec\min\{u_1-y,z\}$, which implies $y+z\in A_1+A_0\subseteq (2A)$.
If $y\in A_1$ and $y\sim u_1$, then
$A_1[y-t,y-1]\cap (y+z-A_0[z+1,z+t])\not=\emptyset$ for some
$0\prec t\prec\min\{y,u_0-z\}$, which implies $y+z\in A_1+A_0\subseteq (2A)$.
If $y\in A_2$, then $0\prec y+z\prec u_0\sim 2u_1$. Since $A_1$ is full,
then there are $x\sim 2l_1$ and $x'\sim 2u_1$ such that 
$(J_0+J_2)[x,x']\subseteq (2A_1)$. Hence $y+z\in (2A_1)\subseteq (2A)$.
By all the arguments above we have $(2A')=(2A)$.

For the case that $z\equiv j+1\,(\mod 3)$ the proof is similar.
\quad $\Box$(Claim \ref{bplittle}.1)

\medskip

{\bf Claim \ref{bplittle}.2}:\quad If $\frac{u_1}{2}<z<u_1$ and
$z\equiv j+1\,(\mod 3)$, then $z\in A$.

Proof of Claim \ref{bplittle}.2:\quad Suppose not and
let $z$ be the least number such that the claim is not true.
By Claim \ref{bplittle}.1 we have $z\sim u_1$.
Let $A'=A\cup\{z\}$. It suffices to show $|(2A')\smallsetminus (2A)|\leqslant 2$.
Let $y\in A'$.

If $y\in A_0$ and $y\prec u_0$, then
$A_0[y+1,y+t]\cap (y+z-A_1[z-t,z-1])\not=\emptyset$ for some
$0\prec t\prec\min\{u_0-y,z\}$, which implies $y+z\in A_0+A_1\subseteq (2A)$.
If $y\in A_0$ and $u_0\sim y<u_0$, then $y+z=u_0+(z-(u_0-y))\in A_0+A_1\subseteq (2A)$
by the minimality of $z$. If $y\in A_1\cup\{z\}$ and $y\prec u_1$, then
$A_1[y+1,y+t]\cap (y+z-A_1[z-t,z-1])\not=\emptyset$ for some
$0\prec t\prec\min\{u_1-y,z\}$, which implies $y+z\in (2A_1)\subseteq (2A)$.
If $y\in A_1\cup\{z\}$ and $u_1\sim y<u_1$, then $y+z=u_1+(z-u_1+y)\in (2A_1)\subseteq (2A)$.
If $y\in A_2$, then $y+z\in (2A_0)$ by the facts that $A_0$ is full and
$2l_0\prec y+z\prec 2u_0$. Hence $((2A')\smallsetminus (2A))\subseteq\{z+u_0,z+u_1\}$.
\quad $\Box$(Claim \ref{bplittle}.2)

\medskip

{\bf Claim \ref{bplittle}.3}:\quad If $\frac{u_0}{2}<z<u_0$ and $z\equiv j\,(\mod 3)$,
then $z\in A$.

Proof of Claim \ref{bplittle}.3:\quad Suppose not and let $z$ be the least
number such that the claim is not true. By Claim \ref{bplittle}.1 we have
$z\sim u_0$. Let $A'=A\cup\{z\}$. Again it suffices to show
$|(2A')\smallsetminus (2A)|\leqslant 2$. Let $y\in A'$.

If $y\in A_0\cup\{z\}$ and $y\prec u_0$, then
$A_0[y+1,y+t]\cap (y+z-A_0[z-t,z-1])\not=\emptyset$ for some
$0\prec t\prec\min\{u_0-y,z\}$, which implies $y+z\in (2A_0)\subseteq (2A)$.
If $y\in A_0\cup\{z\}$ and $u_0\sim y<u_0$, then $y+z=u_0+(z-(u_0-y))\in (2A_0)
\subseteq (2A)$ by the minimality of $z$. If $y\in A_1$ and $y\prec u_1$, then
$A_1[y+1,y+t]\cap (y+z-A_0[z-t,z-1])\not=\emptyset$ for some
$0\prec t\prec\min\{u_1-y,z\}$, which implies $y+z\in A_1+A_0\subseteq (2A)$.
If $y\in A_1$ and $u_1\sim y\leqslant u_1$, then $y+z=(y-(u_0-z))+u_0\in A_1+A_0
\subseteq (2A)$. Note that $y-u_0+z\in A_1$ by Claim \ref{bplittle}.1 and
Claim \ref{bplittle}.2. If $y\in A_2$ and $y<u_2$, then $y+z=u_2+(z-u_2+y)\in A_2+A_0
\subseteq (2A)$. Hence $((2A')\smallsetminus (2A))\subseteq\{u_0+z,u_2+z\}$.
\quad $\Box$(Claim \ref{bplittle}.3)

\medskip

{\bf Claim \ref{bplittle}.4}:\quad There is an $i\in\{0,1\}$ such that
$l_i<z<\frac{u_i}{2}$ and $z\equiv j+i\,(\mod 3)$ imply $z\in A$.

Proof of Claim \ref{bplittle}.4:\quad Suppose not and let
\[z_i=\max\{z\in [0,H]:l_i<z<\frac{u_i}{2},\,z\equiv j+i\,(\mod 3)\mbox{ and }
z\not\in A_i\}\]
for $i=0,1$. By Claim \ref{bplittle}.1 we have $z_i\sim 0$.

\medskip

{\bf Subclaim \ref{bplittle}.4.1}:\quad $z_0-l_0=z_1-l_1$.

Proof of Subclaim \ref{bplittle}.4.1:\quad Suppose the subclaim is not true.
Without loss of generality we assume $z_0-l_0<z_1-l_1$. Let $A'=A\cup\{z_1\}$.
Since $z_0+l_1<z_1+l_0$, then $z_1+l_0=(z_0+t)+l_1\in A_0+A_1\subseteq (2A)$
for $t=(z_1+l_0)-(z_0+l_1)$ by the maximality of $z_0$. By the similar arguments
in the last several claims we have $((2A')\smallsetminus (2A))\subseteq\{z_1+l_1,z_1+l_2\}$.
This contradicts the maximality of $|A|$ by Lemma \ref{addone}. By a symmetric argument
we can show $z_0-l_0>z_1-l_1$ is also impossible.
\quad $\Box$(Subclaim \ref{bplittle}.4.1)

\medskip

{\bf Case \ref{bplittle}.4.1}:\quad $z_0+l_2<z_1+l_1$.

Let $A'=A\cup\{z_1\}$. Note that $z_0+l_2\equiv z_1+l_1\,(\mod 3)$. Then
$z_1+l_1=z_0+t+l_2\in A_0+A_2\subseteq (2A)$ for $t=(z_1+l_1)-(z_0+l_2)>0$.
Hence by the similar arguments as in
Subclaim \ref{bplittle}.4.1 we can show $((2A')\smallsetminus (2A))\subseteq
\{z_1+l_0,z_1+l_2\}$.\quad $\Box$(Case \ref{bplittle}.4.1)

\medskip

{\bf Case \ref{bplittle}.4.2}:\quad $z_0+l_2>z_1+l_1$.

Let $A'=A\cup\{z_0\}$. Then $z_0+l_2=z_1+t+l_1\in (2A_1)$ for
$t=(z_0+l_2)-(z_1+l_1)>0$ by the maximality of $z_1$. Hence
$((2A')\smallsetminus (2A))\subseteq\{z_0+l_0,z_0+l_1\}$.
\quad $\Box$(Case \ref{bplittle}.4.2)

\medskip

Following the two cases above we have $z_0+l_2=z_1+l_1$. By symmetric arguments
we can also show that $z_1+l_2=z_0+l_0$. Subtracting the second equality from the
first we have $z_0-z_1=z_1-z_0+l_1-l_0$. This implies $2(z_0-z_1)=l_1-l_0$.
But by Subclaim \ref{bplittle}.4.1 we have $z_0-z_1=-(l_1-l_0)$. Hence $l_1-l_0=0$,
which is absurd.\quad $\Box$(Claim \ref{bplittle}.4)

\medskip

{\bf Claim \ref{bplittle}.5}:\quad If $l_i<z<\frac{u_i}{2}$ and $z\equiv j+i\,(\mod 3)$,
then $z\in A$ for $i=0,1$.

Proof of Claim \ref{bplittle}.5:\quad Suppose the claim is not true.
By Claim \ref{bplittle}.4 we can assume, without loss of generality,
that $l_0<z<\frac{u_0}{2}$ and $z\equiv j\,(\mod 3)$ imply $z\in A$.

Let $z_1=\max\{z\in [0,H]:l_1<z<\frac{u_1}{2},\,z\equiv j+1\,(\mod 3)\mbox{ and }
z\not\in A_1\}$. Let $A'=A\cup\{z_1\}$. Then $z_1+l_0=l_1+(l_0+z_1-l_1)\in A_1+A_0
\subseteq (2A)$. It is now easy to see that $((2A')\smallsetminus (2A))\subseteq
\{z_1+l_1,z_1+l_2\}$. This contradicts the maximality of $|A|$.
\quad $\Box$(Claim \ref{bplittle}.5)

\medskip

{\bf Claim \ref{bplittle}.6}:\quad $l_0+l_2\geqslant 2l_1-3$ and
$l_1+l_2\geqslant 2l_0-3$.

Proof of Claim \ref{bplittle}.6:\quad By symmetry we need only to show the
first inequality. Assume it is not true and we have $l_0+l_2\leqslant 2l_1-6$.
Then $l_1\not\in [0,2]$. Let $z=l_1-3$ and $A'=A\cup\{z\}$. Let $y\in A'$.

If $y\in A_0$ and $y\not=l_0$, then $y+z=(y-3)+l_1\in A_0+A_1\subseteq (2A)$.
If $y\in A_1$ and $y\succ l_1$, then $A_1[y-t,y-1]\cap (y+z-A_1[z+1,z+t])\not=
\emptyset$ for some $0\prec t\prec\min\{y,u_1\}$, which implies $y+z\in (2A_1)
\subseteq (2A)$. If $y\in A_1\cup\{z\}$ and $y\sim l_1$, then
$y+z=(l_0+t)+l_2\in A_0+A_2$ for $t=(y+z)-(l_0+l_2)\geqslant (2l_1-6)-(l_0+l_2)
\geqslant 0$. If $y\in A_2$ and $y>l_2$, then $y+z=l_2+(z+y-l_2)\in A_2+A_1\in (2A)$.
Hence $((2A')\smallsetminus (2A))\subseteq\{z+l_0,z+l_2\}$, a contradiction to
the maximality of $|A|$ by Lemma \ref{addone}.\quad $\Box$(Claim \ref{bplittle}.6)

\medskip

{\bf Claim \ref{bplittle}.7}:\quad Let $z=u_2+3$ and $A'=A\cup\{z\}$. Then
(\ref{cond3k-3+b}), (\ref{smallA}), and (\ref{longap}) maintain true with
$A$ and $b$ being replaced by $A'$ and $b'$, respectively.

Proof of Claim \ref{bplittle}.7:\quad By Lemma \ref{addone} it suffices to prove
$|(2A')\smallsetminus (2A)|\leqslant 2$. Suppose not. We derive a contradiction.

\medskip

{\bf Subclaim \ref{bplittle}.7.1}:\quad $l_0+l_1-6\geqslant u_2+l_2$.

Prove of Subclaim \ref{bplittle}.7.1:\quad Assume the subclaim is not true. 
So we have $u_2+l_2\geqslant l_0+l_1-3$. Let $y\in A'$

If $y\in A_0$ and $y<u_0$, then $y+z=(y+3)+(z-3)\in A_0+A_2\subseteq (2A)$.
If $y\in A_1$ and $y<u_1$, then $y+z=(y+3)+(z-3)\in A_1+A_2\subseteq (2A)$.
If $y\in A_2\cup\{z\}$, then $y+z=y+u_2+3\geqslant l_2+u_2+3\geqslant l_0+l_1$.
Hence $y+z=(l_0+t)+l_1\in A_0+A_1\subseteq (2A)$ for $t=(y+z)-(l_0+l_1)\geqslant 0$.
Now we have $((2A')\smallsetminus (2A))\subseteq\{z+u_0,z+u_1\}$, which
contradicts the assumption that $|(2A')\smallsetminus (2A)|\leqslant 2$.
\quad $\Box$(Subclaim \ref{bplittle}.7.1)

\medskip

We now ready to derive a contradiction. By 
Claim \ref{bplittle}.6 and Subclaim \ref{bplittle}.7.1
we have $2(l_0+l_1+l_2)\geqslant 2l_0+2l_1+u_2+l_2$. This implies $l_2\geqslant u_2$.
Hence $A_2=\{l_2\}$. So by Subclaim \ref{bplittle}.7.1 again we have
$l_0+l_1-6\geqslant 2l_2$. Since $0=\min A$, then $0\in\{l_0,l_1,l_2\}$.
We want to show $l_2=0$. Suppose $l_0=0$. Then by Claim \ref{bplittle}.6 and
Subclaim \ref{bplittle}.7.1 we have $l_2\geqslant 2l_1-3$ and $l_1-6\geqslant 2l_2$.
So $l_1-6\geqslant 2(2l_1-3)=4l_1-6$ implies $l_1\geqslant 4l_1$, which is absurd
because $l_0=0$ implies $l_1>0$. By symmetry we also have $l_1>0$. Hence $l_2=0$.

By Claim \ref{bplittle}.6 and Subclaim \ref{bplittle}.7.1 again we have
$l_0\geqslant 2l_1-3$ and $l_1\geqslant 2l_0-3$, which imply
$l_0+l_1\geqslant 2(l_0+l_1)-6$ or equivalently $l_0+l_1\leqslant
6$. Hence by Subclaim \ref{bplittle}.7.1 we have $l_0+l_1=6$. Note that
$l_0\not\equiv l_1\,(\mod 3)$. So $(l_0,l_1)\not=(3,3)$. Assume
$l_0<l_1$. The $(l_0,l_1)=(2,4)$ or $(l_0,l_1)=(1,5)$. But each of
the two cases contradicts the inequality $l_0\geqslant 2l_1-3$ in
Claim \ref{bplittle}.6 with $l_2=0$.\quad $\Box$(Claim \ref{bplittle}.7)

\medskip

By Claim \ref{bplittle}.7 we can add $u_2+3,u_2+6,u_2+9,\ldots$ successively
to $A$ to form a set $A'$ so that (\ref{cond3k-3+b}), (\ref{smallA}), 
and (\ref{longap}) maintain true with $A$ and $b$ being replaced by $A'$ 
and $b'$, respectively. However, (\ref{smallA}) will be eventually violated
in this process.
\quad $\Box$(Lemma \ref{bplittle})

\begin{lemma}\label{bplittle2}

Let $A_i=\{z\in A:z\equiv i\,(\mod 3)\}$ for $i=0,1,2$. If 
there is an $i\in [0,2]$ such that $\max A_i-\min A_i\sim 0$, 
then $H+1\leqslant 2|A|-1+2b$.

\end{lemma}

\noindent {\bf Proof}:\quad The ideas are same 
as in the proof of Lemma \ref{bplittle}.
We will describe the steps without too much 
technical details. Let $I_i=i+(3\ast\,^*\!\nat)$, $A_i=A\cap I_i$,
$l_i=\min A_i$, and $u_i=\max A_i$ for $i=0,1,2$.

Without loss of generality let $u_2\sim l_2$. 
By Lemma \ref{bplittle} we can assume
$0\prec l_2\leqslant u_2\prec u_2$. Since 
\[|2A|\succeq |2A_0|+|2A_1|+|A_0+A_1|\succeq 3|A_0|+3|A_1|\sim 3|A|,\]
then $|2A|\sim 3|A|$ implies that $A_i$ is full for $i=0,1$.
Note that $|A_0|\succeq\frac{H}{6}$, $|A_1|\succeq\frac{H}{6}$, and
$|A_0\cup A_1|\sim\frac{1}{2}H$. 

Suppose the lemma is not true. Without loss of generality we can
assume that $|A|$ is the maximum
among all the sets in $\cup_{i=0}^2I_i[l_i,u_i]$ containing the original 
set and satisfying (\ref{cond3k-3+b}), (\ref{smallA}), and (\ref{longap}).
Without loss of generality let's assume $l_0=0$.

\medskip

{\bf Case \ref{bplittle2}.1}\quad $H=u_0$.

If $l_2\prec 2l_1$, then
\[|2A|\succeq |2A_0|+|2A_1|+|A_0+A_1|+|l_2+A_0[0,2l_1-l_2]\succ 3|A|.\]
So we can assume $2l_1\preceq l_2$. By symmetry we can assume
$u_2+H\preceq 2u_1$. Since $u_1-l_1\sim\frac{1}{2}H$ by (\ref{half}),
we have $2l_1\sim l_2$ and $u_2+H\sim 2u_1$. Hence 
\[A_2+A_1\subseteq (2A_0)[l_1+l_2,u_1+u_2]=(2I_0)[l_1+l_2,u_1+u_2].\]
This implies $A_1=I_1[l_1,u_1]$ by Lemma \ref{addone} and the maximality of $|A|$. 
Then we can show $A_0=I_0[l_0,u_0]$ again by Lemma \ref{addone} and the
maximality of $|A|$. Furthermore, we can show $A_2=I_2[l_2,u_2]$ by the fact that
$(2I_2)[2l_2,2u_2]\subseteq A_0+A_1$ and by Lemma \ref{addone}.
Now we add $z=u_2+3$, $z=u_2+6$, $z=u_2+9$, etc.\ successively to $A$ 
so that the set maintains satisfying (\ref{cond3k-3+b}), (\ref{smallA}), 
and (\ref{longap}). However, this process will eventually violate
(\ref{smallA}).\quad $\Box$(Case \ref{bplittle2}.1)

\medskip

{\bf Case \ref{bplittle2}.2}\quad $H=u_1$.

We can again show that $2l_1\preceq l_2$ and $u_2+H\preceq 2u_0$
because otherwise we can show $|2A|\succ 3|A|$. Since 
$H-l_1+u_0\sim\frac{3}{2}H$, then $u_0-l_1\sim\frac{H}{2}$, which
implies $2l_1\sim l_2$ and $u_2+H\sim 2u_0$. Again assume that
$A\subseteq \cup_{i=0}^2I_i[l_i,u_i]$ has the maximum cardinality among
the sets satisfying (\ref{cond3k-3+b}), (\ref{smallA}), and 
(\ref{longap}). Then we can show $A_i=I_i[l_0,u_0]$ for $i=0,1,2$.
Finally we can again add $z=u_2+3$, $z=u_2+6$, $z=u_2+9$, etc.\ 
successively to $A$ so that the set maintains satisfying 
(\ref{cond3k-3+b}), (\ref{smallA}), and (\ref{longap}). Again 
this process will eventually violate (\ref{smallA}).
\quad $\Box$(Lemma \ref{bplittle2})

\begin{lemma}\label{bplittle3}

Suppose there are $0\prec a\sim c\prec H$ such that 
$A[0,a]$ is a backward triangle as well as a subset of a \bp\ of 
difference $3$ and $A[c,H]$ is a forward triangle as well as a
subset of a \bp\ of difference $3$. Then $H+1\leqslant 2|A|-1+2b$.

\end{lemma}

\noindent {\bf Proof}:\quad The ideas are again the same as in the
proof of Lemma \ref{bplittle}. Let $I_i=(i+(3\ast\,^*\!\nat))$ for
$i=0,1,2$. By Lemma \ref{bplittle2} we can
assume that $A[0,a]\subseteq I_0\cup I_1$ and 
$A[c,H]\subseteq I_0\cup I_2$. For $i=0,1,2$ let
$A_i=A\cap I_i$, $l_i=\min A_i$, $u_i=\max A_i$, and 
$J_i=I_i\cap [l_i.u_i]$. Then
we have $u_1\sim l_2\sim a$. Suppose the lemma is not true.
Then $A$ satisfies (\ref{cond3k-3+b}), (\ref{smallA}), and 
(\ref{longap}). We again assume the maximality of $|A|$ for
$A\subseteq J_0\cup J_1\cup J_2$ satisfying (\ref{cond3k-3+b}), 
(\ref{smallA}), and (\ref{longap}). 
By Lemma \ref{addone} we can prove that for each $x$,
$l_i\prec x\prec u_i$ implies $x\in A$. Then we can prove
$A_i=J_i$ by the same ideas as in the proof of Lemma \ref{bplittle}.
Now we can add $u_1+3,u_1+6,u_1+9,\ldots$ successively to $A$
such that the conditions (\ref{cond3k-3+b}), (\ref{smallA}), and 
(\ref{longap}) maintain true. But this process will eventually
violate (\ref{smallA}). \quad $\Box$(Lemma \ref{bplittle3})

\begin{lemma}\label{middleu2}

If there is a $x\sim 0$ in $A$ such that $\gcd(A[x,H]-x)=d>1$, 
then $H+1\leqslant 2|A|-1+2b$.

\end{lemma}

\noindent {\bf Proof}:\quad Since (\ref{half}) is true, then $d=2$. 
Let $c=\min\{x\in A:\gcd(A[x,H]-x)=2\}$.
Then $c\sim 0$, $c>0$, and $A[c,H]$ is full. Let $E$ be the
set of all even numbers and $O$ be the set of all odd numbers. Let
$A_e=A\cap E$ and $A_o=A\cap O$. Let $l_e=\min A_e$, $l_o=\min
A_o$, $u_e=\max A_e$, and $u_o=\max A_o$.

\medskip

{\bf Case \ref{middleu2}.1}\quad $c$ is even.

Then $l_e=0$, $u_e=H$, and $A_e$ is full. We want to show $H+1\leqslant
2|A|-1+2b$.

Let $|2A_e|=2|A_e|-1+b_e$. Then $b_e\sim 0$.
By Theorem \ref{2k-1+b} we have $\frac{H}{2}+1\leqslant |A_e|+b_e$.
On the other hand, by Theorem \ref{A+B},
\[|A_e+A_o|\geqslant\min\{|A_e|+2|A_o|-2,\frac{H}{2}+|A_o|\}.\]

If $\frac{H}{2}+|A_o|\geqslant |A_e|+2|A_o|-2$, then
\begin{eqnarray*}
\lefteqn{3|A|-3+b=|2A|}\\
& &\geqslant 2|A_e|-1+b_e+|A_e|+2|A_o|-2\\
& &\geqslant 3|A|-3+b_e-|A_o|.
\end{eqnarray*}
This implies $b_e\leqslant b+|A_o|$. Hence
$\frac{H}{2}+1\leqslant |A_e|+b_e\leqslant |A|+b$,
which implies $H+1\leqslant 2|A|-1+2b$.

If $\frac{H}{2}+|A_o|<|A_e|+2|A_o|-2=|A|+|A_o|-2$, then $|A|>\frac{H}{2}+2$,
which contradicts (\ref{smallA}). \quad $\Box$(Case \ref{middleu2}.1)

\medskip

{\bf Case \ref{middleu2}.2}\quad $c$ is odd.

Clearly, $0=l_e$ and $H=u_o$. If $l_o>u_e$, then $A$ is a subset
of a \bp\, Hence we can assume $l_o<u_e$ and need to show $H+1\leqslant
2|A|-1+2b$. Suppose $H+1>2|A|-1+2b$.
Let \[S=\{x\in O[l_o,u_o]:\frac{A_o(l_o,x)}{x-l_o+1}
\leqslant\frac{1}{4}\}.\] If $S\not=\emptyset$, let $l'=(\max
S)+1$. Otherwise, let $l'=l_o$. Let \[T=\{x\in
O[l_o,u_o]:\frac{A_o(x,u_o)}{u_o-x+1} \leqslant\frac{1}{4}\}.\] If
$T\not=\emptyset$, let $u'=(\min T)-1$. Otherwise, let $u'=u_o$.
Note that $l',u'\in A_o$. Since $A_o$ is full, then $l'\sim l_o$
and $u'\sim u_o$. For each $x\in O[l',u_o]$, we have
$\frac{A_o(l_o,x)}{x-l_o+1}>\frac{1}{4}$ and for any $x\in
O[l_o,u']$, we have $\frac{A_o(x,u_o)}{u_o-x+1}>\frac{1}{4}$.
By the pigeonhole principle
\[E[l_o+l',u_o+u']\subseteq (2A_o).\]
Let $p=O(l',u')-A_o(l',u')$, $p'=O(l',l'+u_e)-A_o(l',l'+u_e)$, and
$p''=p-p'$.

Let $\bar{A}_o=A_o\cup O[l',u']$. Then $2\bar{A}_o=2A_o$. Hence
\[|2\bar{A}_o|=|2A_o|=2|A_o|-1+b_o=2|\bar{A}_o|-1+b_o-2p.\]
This implies, by Theorem \ref{2k-1+b},
\[\frac{1}{2}(H-l_o)+1\leqslant |\bar{A}_o|+b_o-2p=|A_o|+b_o-p.\]

\medskip

{\bf Subcase \ref{middleu2}.2.1}\quad $p''\geqslant 2|A_e|$.

Since
\begin{eqnarray*}
\lefteqn{3|A|-3+b=|2A|}\\
& &\geqslant |2A_o|+|A_e+A_o|\\
& &\geqslant 2|A_o|-1+b_o+|0+A_o[l_o,l'+u_e-2]|+|u_e+A_o[l',H]|\\
& &\geqslant 2|A_o|-1+b_o+A_o(l_o,l'+u_e-2)+A_o(l',H)\\
& &\geqslant 3|A|-3+b_o-3|A_e|+A_o(l',l'+u_e)+1,
\end{eqnarray*}
then $b_o\leqslant b+3|A_e|-A_o(l',l'+u_e)-1$. Hence
\begin{eqnarray*}
\lefteqn{\frac{1}{2}(H-l_o)+1\leqslant |A_o|+b_o-p}\\
& &\leqslant |A_o|+b+3|A_e|-A_o(l',l'+u_e)-1-p\\
& &\leqslant |A|+b+(2|A_e|-p'')-O(l',l'+u_e)-1\\
& &\leqslant |A|+b-\frac{1}{2}u_e-2<|A|+b-\frac{1}{2}u_e.
\end{eqnarray*}
This implies $H-l_o+2\leqslant 2|A|+2b-u_e$ and \[H+1\leqslant
2|A|-1+2b-(u_e-l_o)<2|A|-1+2b.\]
\quad $\Box$(Subcase \ref{middleu2}.2.1)

\medskip

{\bf Subcase \ref{middleu2}.2.2}\quad $p''<2|A_e|$.

Since $u'-(l'+u_e)\succ 0$ and $u_e\sim 0$, then there exists a
$t,t'\in O[l'+2u_e+1,u'-2u_e-1]$ such that
$t'-t>\max\{u_o-u',l'-l_o\}$ and
$|A_e|+A_o(t-u_e,t'+2u_e)>\frac{3}{2}u_e+\frac{t'-t}{2}+1$
because otherwise we can find three disjoint intervals of length
$t'-t+3u_e+1$ for $t'-t=2+\max\{u_o-u',l'-l_o\}$ in 
$[l'+u_e+2,u']$ such that each contains
at least $|A_e|$ elements from the set $O\smallsetminus A_o$. This
contradicts the assumption $p''<2|A_e|$. Suppose $H+1>2|A|-1+2b$.
We want to derive a contradiction by induction on the size of the
counterexamples $A'\supseteq A$.

Suppose $A'=A_e\cup A'_o$ is the set with the maximum cardinality
$|A'|$ such that $A_o\subseteq A'_o\subseteq O[l_o,u_o]$,
$|2A'|=3|A'|-3+b'$ for $0\leqslant b'\leqslant b$, and $H+1>2|A'|-1+2b'$.

\medskip

{\bf Claim \ref{middleu2}.2.2.1}\quad $A'_o[t,t'+u_e]=O[t,t'+u_e]$.

Proof of Claim \ref{middleu2}.2.2.1:\quad Suppose the claim is not
true and let $g\in O[t,t'+u_e]\smallsetminus A'_o$. Let
$A_o''=A_o'\cup\{g\}$ and $A''=A_e\cup A_o''$. 
Note that if $x\in A_o''[l_o,u_o]$, then
$l_o+l'<x+g<u_o+u'$. Hence $x+g\in (2A_o)\subseteq (2A')$. Let
$x\in A_e[0,u_e]$. Since
\begin{eqnarray*}
\lefteqn{|A_e|+A'_o(g+x-u_e,g+x)}\\
& &=|A_e|+A'_o(t-u_e,t'+2u_e)\\
& &\quad -A'_o(t-u_e,g+x-u_e-2)-A'_o(g+x+2,t'+2u_e)\\
& &>\frac{3}{2}u_e+\frac{t'-t}{2}+1-u_e-\frac{t'-t}{2}\\
& &=\frac{1}{2}u_e+1,
\end{eqnarray*}
then $(g+x-A_e)\cap A'_o[g+x-u_e,g+x]\not=\emptyset$.
Hence $g+x\in A'_o+A_e\subseteq (2A')$. Now we have that
$(2A'')=(2A')$, which contradicts the maximality of $|A'|$ by Lemma \ref{addone}.
\quad $\Box$(Claim \ref{middleu2}.2.2.1)

\medskip

{\bf Claim \ref{middleu2}.2.2.2}\quad
$A'_o[t'+u_e+2,H]=O[t'+u_e+2,H]$.

Suppose the claim is not true and let $g=\min
(O[t'+u_e+2,H]\smallsetminus A'_o)$. Let $A''_o=A'_o\cup\{g\}$, and
$A''=A_e\cup A''_o$. Then $g>t'+u_e$. Note that if $x\in
A''_o[l',u']$, then $l_o+l'\leqslant x+g\leqslant u_o+u'$. Hence
$x+g\in (2A')$. If $x\in A''_o[l_o,l'-2]$, then $y=g-(l'-x)\in A'_o$
by the minimality of $g$ and Claim \ref{middleu2}.2.2.1. Hence
$x+g=l'+y\in (2A')$. If $x\in A''_o[u'+2,H-2]$, then
$x+g=H+(g-(H-x))\in (2A')$. If $x\in A_e[0,u_e-2]$, then
$x+g=u_e+(g-(u_e-x))\in (2A')$. From the arguments above, we have
$(2A'')\smallsetminus (2A')\subseteq \{H+g,u_e+g\}$, which
contradicts the maximality of $|A'|$ by Lemma \ref{addone}. \quad
$\Box$(Claim \ref{middleu2}.2.2.2)

\medskip

{\bf Claim \ref{middleu2}.2.2.3}\quad $A'_o[l_o,t-2]=O[l_o,t-2]$.

Proof of Claim \ref{middleu2}.2.2.3:\quad Suppose the claim is not
true and let $g=\max (O[l_o,t-2]\smallsetminus A'_o)$. Let
$A''=A'\cup\{g\}$. Then $(2A'')\smallsetminus (2A')
\subseteq\{0+g,l_o+g\}$, which contradicts the maximality of
$|A'|$. \quad $\Box$(Claim \ref{middleu2}.2.2.3)

\medskip

By the three claims above we have $A'_o=O[l_o,u_o]$.

Without loss of generality, we can assume that our original
counterexample $A$ satisfies $A_o=O[l_o,u_o]$. Hence
$|A_o|=\frac{H-l_o}{2}+1\mbox{ or }2|A_o|-1=H-l_o+1$.
Note also that since $H+1=H-l_o+1+l_o=2|A|-1+l_o-2|A_e|$, we can
assume that $l_o\geqslant 2|A_e|+1$.

\medskip

{\bf Claim \ref{middleu2}.2.2.4}\quad $2l_o\geqslant u_e+4$.

Proof of Claim \ref{middleu2}.2.2.4:\quad Suppose $2l_o\leqslant
u_e+2$. For each even number $z\in E[u_e,H-1]$ let $A_z=A\cup
E[u_e,z]$. Let
\begin{eqnarray*}
\lefteqn{S=\{z\in E[u_e,H-1]:|2A_z|=3|A_z|-3+b_z\geqslant 3|A_z|-3,}\\
& & b_z\leqslant b,\mbox{ and }H+1>2|A_z|-1+2b_z\}.
\end{eqnarray*}
Clearly, $u_e\in S$. For each $z\succ 0$,
$|A_z|=|A|+E(u_e+2,z)\succ\frac{1}{2}H$,
which implies $z\not\in S$. Let $z_0=\max S\sim 0$. We now derive a
contradiction.

By Theorem \ref{3k-3} we can assume $b_{z_0}>0$.
Note that since $|A_{z_0+2}|=|A_{z_0}|+1$, we have
\[H+1>2|A_{z_0}|-1+2b_{z_0}=2|A_{z_0+2}|-1+2(b_{z_0}-1).\]
Since for each $x\in A_{z_0+2}\cap E$, $x+z_0+2\in E[2l_o,2H]$.
Hence $(2A_{z_0+2})\smallsetminus (2A_{z_0})\subseteq\{z_0+2+H\}$,
which contradicts the maximality of $z_0$ by Lemma \ref{addone}. 
\quad $\Box$(Claim \ref{middleu2}.2.2.4)

\medskip

Let $d'=\gcd(A_e)$. Then $d'$ must be an even number. Let $q=\min
(A_e[l_o+1,u_e])$ and $q'=\max (A_e[0,l_o-1])$.

\medskip

{\bf Subsubcase \ref{middleu2}.2.2.1}\quad $d'\geqslant 4$.

First we can assume that $u_e=q$ by the following argument.

Let $A'=A\smallsetminus A_e[q+2,u_e]$. Then
$|A'|=|A|-A_e(q+2,u_e)$. Note that $A_e(q+2,u_e)\leqslant
\frac{u_e-q}{d'}\leqslant\frac{u_e-q}{4}$. Since $(2A')\subseteq
(2A)\smallsetminus O[q+2+H,u_e+H]$ and $O[q+2+H,u_e+H]\subseteq
A_e+A_o\subseteq (2A)$, then there is a $b'\geqslant 0$ such that
\begin{eqnarray*}
\lefteqn{3|A'|-3+b'=|2A'|}\\
& &\leqslant |2A|-\frac{u_e-q}{2}\\
& &=3|A|-3+b-\frac{u_e-q}{2}\\
& &=3|A'|-3+b-\frac{u_e-q}{2}+3A_e(q+2,u_e)\\
& &\leqslant 3|A'|-3+b+A_e(q+2,u_e).
\end{eqnarray*}
This shows $b'\leqslant b+A_e(q+2,u_e)$. So
\begin{eqnarray*}
\lefteqn{H+1>2|A|-1+2b}\\
& &=2|A'|+2A_e(q+2,u_e)-1+2b\\
& &=2|A'|-1+2(b+A_e(q+2,u_e))\\
& &\geqslant 2|A'|-1+2b'.
\end{eqnarray*}
Hence $A'$ is the desired counterexample. Now identify $A$ with $A'$.

\medskip

{\bf Subsubsubcase \ref{middleu2}.2.2.1.1}\quad $u_e=l_o+1$.

Since $q'\leqslant u_e-d'<u_e-2$, then $q'+u_e\leqslant 2l_o-2$
and $q'+u_e\not\in A_e[0,q']+A_e[0,q']$. Hence we have
$(2A_e)(0,2l_o-2)\geqslant 2A_e(0,q')$.
So we have
\begin{eqnarray*}
\lefteqn{|2A|=|2A_o|+|A_o+A_e|+(2A_e)(0,2l_o-2)}\\
& &\geqslant 2|A_o|-1+|A_o|+\frac{u_e}{2}+2A_e(0,q')\\
& &\geqslant 3|A|-1-3|A_e|+\frac{u_e}{2}+2|A_e|-2\\
& &\geqslant 3|A|-3-|A_e|+\frac{u_e}{2}.
\end{eqnarray*}
This shows $-|A_e|+\frac{u_e}{2}\leqslant b$. On the other hand,
\begin{eqnarray*}
\lefteqn{H+1=2|A_o|-1+l_o}\\
& &=2|A|-1+l_o-2|A_e|\\
& &\leqslant 2|A|-1+l_o+2(b-\frac{u_e}{2})\\
& &=2|A|-1+2b+(l_o-u_e)<2|A|-1+2b.
\end{eqnarray*}
This contradicts the assumption that $H+1>2|A|-1+2b$.
\quad $\Box$(Subsubsubcase \ref{middleu2}.2.2.1.1)

\medskip

{\bf Subsubsubcase \ref{middleu2}.2.2.1.2}\quad $u_e>l_o+1$.

Then we have
\begin{eqnarray*}
\lefteqn{|2A|=|2A_o|+|A_o+A_e|+(2A_e)(0,2l_o-2)}\\
& &\geqslant 2|A_o|-1+|A_o|+\frac{u_e}{2}+2A_e(0,q')-1\\
& &\geqslant 3|A|-1-3|A_e|+\frac{u_e}{2}+2|A_e|-3\\
& &\geqslant 3|A|-3-|A_e|+\frac{u_e}{2}-1.
\end{eqnarray*}
This shows $-|A_e|+\frac{u_e}{2}-1\leqslant b$. On the other hand,
\begin{eqnarray*}
\lefteqn{H+1=2|A_0|-1+l_o}\\
& &=2|A|-1+l_o-2|A_e|\\
& &\leqslant 2|A|-1+l_o+2(b-\frac{u_e}{2}+1)\\
& &=2|A|-1+2b+(l_o-u_e+2)\\
& &\leqslant 2|A|-1+2b,
\end{eqnarray*}
by the assumption $u_e>l_o+1$. This again contradicts
$H+1>2|A|-1+2b$. \quad $\Box$(Subsubcase \ref{middleu2}.2.2.1)

\medskip

{\bf Subsubcase \ref{middleu2}.2.2.2}\quad $d'=2$.

\medskip

{\bf Subsubsubcase \ref{middleu2}.2.2.2.1}\quad $|2A_e|=2|A_e|-1+b_e$ for
some $b_e<|A_e|-2$.

Since $A_e\subseteq E$, then by Theorem \ref{2k-1+b}
$\frac{u_e}{2}+1\leqslant |A_e|+b_e$. We also have
\begin{eqnarray*}
\lefteqn{|2A|=|2A_o|+|2A_e|+|A_o+A_e|-(2A_e)(2l_o,2u_e)}\\
& &\geqslant 2|A_o|-1+2|A_e|-1+b_e+|A_o|+\frac{u_e}{2}-\frac{2u_e-2l_o}{2}-1\\
& &\geqslant 3|A|-3-|A_e|+b_e+\frac{u_e}{2}+l_o-u_e,
\end{eqnarray*}
which implies $-|A_e|+b_e-\frac{u_e}{2}+l_o\leqslant b$. Hence
\begin{eqnarray*}
\lefteqn{\frac{u_e}{2}+1+\frac{H-l_o}{2}+1}\\
& &\leqslant |A_e|+b_e+|A_o|=|A|+b_e\\
& &\leqslant |A|+b+|A_e|+\frac{u_e}{2}-l_o.
\end{eqnarray*}
This implies $u_e+2+H-l_o+2\leqslant 2|A|+2b+2|A_e|+u_e-2l_o$.
Hence \[H+1\leqslant 2|A|-1+2b+(2|A_e|-l_o-2)<2|A|-1+2b\] because
\[|A|=|A_e|+|A_o|=|A_e|+\frac{H-l_o}{2}+1\leqslant\frac{1}{2}(H+1)\]
implies $l_o> 2|A_e|$.
\quad $\Box$(Subsubsubcase \ref{middleu2}.2.2.2.1)

\medskip

{\bf Subsubsubcase \ref{middleu2}.2.2.2.2}\quad $|2A_e|\geqslant 3|A_e|-3$
and $u_e-l_o\leqslant 2|A_e|-4$.

Since
\begin{eqnarray*}
\lefteqn{|2A|=2|A_o|-1+|A_o|+\frac{u_e}{2}+(2A_e)(0,2l_o-2)}\\
& &=3|A|-3+\frac{u_e}{2}-3|A_e|+(2A_e)(0,2l_o-2)+2,
\end{eqnarray*}
then $\frac{u_e}{2}-3|A_e|+(2A_e)(0,2l_o-2)+2\leqslant b$. Hence
$\frac{H-l_o}{2}+1=|A_o|$ implies
\begin{eqnarray*}
\lefteqn{H+1=2|A|-1+l_o-2|A_e|}\\
& &\leqslant 2|A|-1+l_o+2(b-\frac{u_e}{2}+2|A_e|-(2A_e)(0,2l_o-2)-2)\\
& &\leqslant 2|A|-1+2b+(4|A_e|-4-(u_e-l_o)-2(2A_e)(0,2l_o-2)).
\end{eqnarray*}
Since
\begin{eqnarray*}
\lefteqn{(2A_e)(0,2l_o-2)=(2A_e)(0,2u_e)-(2A_e)(2l_o,2u_e)}\\
& &\geqslant 3|A_e|-3-(u_e-l_o+1)=3|A_e|-4-u_e+l_o,
\end{eqnarray*}
then
\begin{eqnarray*}
\lefteqn{4|A_e|-4-(u_e-l_o)-2(2A_e)(0,2l_o-2)}\\
& &\leqslant 4|A_e|-4-(u_e-l_o)-2(3|A_e|-4-u_e+l_o)\\
& &\leqslant -2|A_e|+4+u_e-l_o\leqslant 0.
\end{eqnarray*}
Hence $H+1\leqslant 2|A|-1+2b$, a contradiction.
\quad $\Box$(Subsubsubcase \ref{middleu2}.2.2.2.2)

\medskip

{\bf Subsubsubcase \ref{middleu2}.2.2.2.3}\quad $|2A_e|\geqslant 3|A_e|-3$
and $u_e-l_o> 2|A_e|-4$.

This time we use the fact that $(2A_e)(0,2l_o-2)\geqslant |0+A_e|=|A_e|$ 
implied by Claim \ref{middleu2}.2.2.4. Since
\begin{eqnarray*}
\lefteqn{|2A|\geqslant 3|A|-3+\frac{u_e}{2}-3|A_e|+(2A_e)(0,2l_o-2)+2}\\
& &\geqslant 3|A|-3+\frac{u_e}{2}-2|A_e|+2,
\end{eqnarray*}
then $\frac{u_e}{2}-2|A_e|+2\leqslant b$. Hence
\begin{eqnarray*}
\lefteqn{H+1=2|A|-1+l_o-2|A_e|}\\
& &\leqslant 2|A|-1+l_o+2(b-\frac{u_e}{2}+|A_e|-2)\\
& &=2|A|-1+2b+(l_o-u_e+2|A_e|-4)\\
& &<2|A|-1+2b.
\end{eqnarray*}
This ends the proof of the lemma. \quad
$\Box$(Lemma \ref{middleu2})

\begin{lemma}\label{twoaps}

Suppose $A=A_e\cup A_o$, where $A_e=A\cap E$ is the set of all
even numbers in $A$ and $A_o=A\cap O$ is the set of all odd
numbers in $A$. Let $u_i=\max A_i$ for $i=e,o$ and $l_o=\min A_o$. 
If (a) $u_e=H$ and $u_o-l_o\sim 0$ or (b) $H=u_o$, $0\prec l_o<u_e\prec
H$, and $u_e-l_o\sim 0$, then $H+1\leqslant 2|A|-1+2b$.

\end{lemma}

\noindent {\bf Proof}:\quad The proof of (a) of the lemma is identical
to the proof of Case \ref{middleu2}.1. We sketch the proof of (b) using
Lemma \ref{addone}. It is easy to see that $A_e$ is full in $E[0,u_e]$
and $A_o$ is full in $O[l_o,H]$. Suppose the lemma is not true.
Following the same ideas as in the proof of Claim \ref{bplittle}.1,
Claim \ref{bplittle}.2, and Claim \ref{bplittle}.3, we 
can assume that $A=E[0,u_e]\cup O[l_o,H]$. However, this implies
\[|A|=\frac{u_e}{2}+1+\frac{H-l_o}{2}+1=\frac{H+(u_e-l_o)}{2}+2>\frac{H+1}{2},\]
which contradicts (\ref{smallA}).
\quad $\Box$(Lemma \ref{twoaps})

\begin{lemma}\label{ftlittle}

If $A$ is a forward triangle from 
$0$ to $H$, then $H+1\leqslant 2|A|-1+2b$.

\end{lemma}

\noindent {\bf Proof}:\quad Assume the lemma is not true and we
need to derive a contradiction. Clearly we have
$\underline{d}_{U}(A)\geqslant\frac{1}{2}$. Furthermore we can
assume $\underline{d}_{U}(A)\geqslant\frac{2}{3}$ by the
following argument: If $\underline{d}_{U}(A)<\frac{2}{3}$,
then there exists $y'\succ 0$ in $A$ such that 
for all $0\prec y\leqslant y'$
\[(2A)(0,2y)\sim 2y=3\cdot\frac{2}{3}y\succ 3A(0,y).\]
If for every
$x\succ 0$ in $A$ we have $\gcd(A[x,H]-x)>1$, then there is a $u\sim 0$
in $A$ such that $\gcd(A[u,H]-u)>1$ by Lemma \ref{overspill}. This 
implies that for each $x\succ 0$ we have
$A(0,x)\preceq\frac{1}{2}x$, which contradicts
that $A[0,H]$ is a forward triangle. Hence we can choose $y$ with
$0\prec y\leqslant y'$ in $A$ such that $\gcd(A[y,H]-y)=1$.
By Lemma \ref{smallbig} we have $|2A|\succ 3|A|$. Let
\[z=\max\{x\in U:A(0,x-1)\leqslant\frac{1}{2}x\}.\]
Note that the smallest possible value of $z$ is $0$.
Note also that $z$ is well defined because
$\underline{d}_{U}(A)\geqslant\frac{2}{3}$. It is easy to check
that $z\in A$, $A(0,z-1)=\frac{1}{2}z$, and for every $x\geqslant z$
in $U$ we have $A(z,x)>\frac{1}{2}(x-z+1)$ by the maximality of $z$.

Define $a$ by
\[a=\min\{x\in [z,H]:A(z,x)\leqslant\frac{1}{2}(x-z+1)\}.\]
The number $a$ is well defined by the fact that
$A(0,z-1)=\frac{1}{2}z$, $H\in A$, and
$|A|\leqslant\frac{1}{2}(H+1)$. It is also easy to check that
$a\sim H$, $A[z,a]$ is a forward triangle from $z$ to $a$,
$a\not\in A$, and $A(z,a)=\frac{1}{2}(a-z+1)$. If $z\leqslant x<a$,
then $A(z,x)>\frac{1}{2}(x-z+1)$ by the minimality of $a$.
Let $a'=\max (A[z,a])$.

\medskip

{\bf Claim \ref{ftlittle}.1}:\quad $[z,a+z-1]\subseteq (2A)$.

Proof of Claim \ref{ftlittle}.1:\quad Let $x\in [z,a+z-1]$.
If $x<2z$, then $x\prec H\sim a$. Hence
\[A(0,x)=A(0,z-1)+A(z,x)>\frac{1}{2}z+\frac{1}{2}(x-z+1)
=\frac{1}{2}(x+1).\]
This implies $A[0,x]\cap (x-A[0,x])\not=\emptyset$.
Hence $x\in (2A)$. If $x\geqslant 2z$, then
$A(z,x-z)>\frac{1}{2}(x-2z+1)$. Hence
$A[z,x-z]\cap (x-A[z,x-z])\not=\emptyset$. This again implies
$x\in (2A)$.\quad $\Box$(Claim \ref{ftlittle}.1)

\medskip

Let $c=\min (A[a+1,H])$. If $2a'<c$, then $a'+H\sim a'+c\prec 2c$.
Hence $A$ is a subset of the \bp\ $[0,a']\cup [c,H]$. So from now on
in this lemma we can assume $2a'\geqslant c$.

\medskip

{\bf Claim \ref{ftlittle}.2}:\quad Suppose $2a'\geqslant a+z$ and
$a+z<x\leqslant\min\{2a',c+z\}$. If $(2A)(a+z,x-1)<\frac{1}{2}
(x-a-z)$, then $x\in (2A)$.

Proof of Claim \ref{ftlittle}.2:\quad Since
$A(z,a)=\frac{1}{2}(a-z+1)$ and $A(a'+1,a)=\emptyset$, then
\begin{eqnarray*}
\lefteqn{A(z+a-a',a')=A(z,a)-A(z,z+a-a'-1)}\\
& &\geqslant\frac{1}{2}(a-z+1)-(a-a')=\frac{1}{2}(2a'-z-a+1).
\end{eqnarray*}
Note that $x-a'\leqslant a'$. Since $a'+A[a+z-a',x-1-a']\subseteq
(2A)[a+z,x-1]$, then $A(a+z-a',x-1-a')<\frac{1}{2}(x-a-z)$.
Hence
\begin{eqnarray*}
\lefteqn{A(x-a',a')=A(a+z-a',a')-A(a+z-a',x-a'-1)}\\
& &>\frac{1}{2}(2a'-z-a+1)-\frac{1}{2}(x-a-z)=
\frac{1}{2}(2a'-x+1).
\end{eqnarray*}
This shows that $A[x-a',a']\cap(x-A[x-a',a'])\not=\emptyset$,
which implies $x\in (2A)$.\quad $\Box$(Claim \ref{ftlittle}.2)

\medskip

Let $S=(2A)[0,z-1]\smallsetminus A[0,z-1]$.

\medskip

{\bf Claim \ref{ftlittle}.3}:\quad Suppose $2a'<c+z$ and
$\max\{a+z,2a'\}\leqslant x<c+z-1$. If
$(2A)(x+1,c+z-1)<\frac{1}{2}(c+z-x-1)$, then $x\in (2A)$ or $x-c\in S$. 

Proof of Claim \ref{ftlittle}.3:\quad Assume $(2A)(x+1,c+z-1)
<\frac{1}{2}(c+z-x-1)$. Since
$c+A[x+1-c,z-1]\subseteq(2A)[x+1,c+z-1]$, then
$A(x+1-c,z-1)<\frac{1}{2}(c+z-x-1)$. Hence
\begin{eqnarray*}
\lefteqn{A(0,x-c)=A(0,z-1)-A(x+1-c,z-1)}\\
& &>\frac{1}{2}z-\frac{1}{2}(c+z-x-1)=\frac{1}{2}(x-c+1).
\end{eqnarray*}
Hence $A[0,x-c]\cap (x-c-A[0,x-c])\not=\emptyset$.
This shows $x-c\in (2A)[0,z-1]$. If $x-c\in A[0,z-1]$, then
$x\in c+A[0,z-1]\subseteq (2A)$. If $x-c\not\in A[0,z-1]$,
then $x-c\in S$.\quad $\Box$(Claim \ref{ftlittle}.3)

\medskip

{\bf Claim \ref{ftlittle}.4}:\quad
$(2A)(0,c+z)\geqslant 3A(0,z-1)+2A(z,a)-1+\frac{1}{2}(c-a+1)$.

Proof of Claim \ref{ftlittle}.4:\quad The proof is divided into
three cases for $2a'\geqslant c+z$, $2a'\leqslant a+z$,
and $a+z<2a'<c+z$.

\medskip

{\bf Case \ref{ftlittle}.4.1}:\quad $2a'\geqslant c+z$.

By Claim \ref{ftlittle}.2
we have $(2A)(a+z,c+z-1)\geqslant\frac{1}{2}(c-a-1)$ (this can be proven
by induction on $x\in [a+z,c+z-1]$). Hence
$(2A)(a+z,c+z)\geqslant\frac{1}{2}(c-a+1)$ because $c+z\in (2A)$.
So we have
\begin{eqnarray*}
\lefteqn{(2A)(0,c+z)}\\
& &=(2A)(0,z-1)+(2A)(z,a+z-1)+(2A)(a+z,c+z)\\
& &\geqslant A(0,z-1)+a+\frac{1}{2}(c-a+1)\\
& &=3A(0,z-1)+2A(z,a)-1+\frac{1}{2}(c-a+1).
\end{eqnarray*}
\quad $\Box$(Case \ref{ftlittle}.4.1)

\medskip

{\bf Case \ref{ftlittle}.4.2}:\quad $2a'\leqslant a+z$.

By Claim \ref{ftlittle}.3 we have $(2A)(a+z,c+z)\geqslant
\frac{1}{2}(c-a+1)-|S|$. Hence
\begin{eqnarray*}
\lefteqn{(2A)(0,c+z)=(2A)(0,z-1)+a+(2A)(a+z,c+z)}\\
& &\geqslant A(0,z-1)+|S|+2A(0,z-1)\\
& &\quad +2A(z,a)-1+\frac{1}{2}(c-a+1)-|S|\\
& &=3A(0,z-1)+2A(z,a)-1+\frac{1}{2}(c-a+1).
\end{eqnarray*}
\quad $\Box$(Case \ref{ftlittle}.4.2)

\medskip

{\bf Case \ref{ftlittle}.4.3}:\quad $a+z<2a'<c+z$.

By Claim \ref{ftlittle}.2 we have $(2A)(a+z,2a')
\geqslant\frac{1}{2}(2a'-a-z+1)$ and by
Claim \ref{ftlittle}.3 we have $(2A)(2a'+1,c+z)
\geqslant\frac{1}{2}(c+z-2a')-|S|$. Hence
\begin{eqnarray*}
\lefteqn{(2A)(0,c+z)=A(0,z-1)+|S|+2A(0,z-1)}\\
& &\quad +2A(z,a)-1+(2A)(a+z,2a')+(2A)(2a'+1,c+z)\\
& &\geqslant 3A(0,z-1)+2A(z,a)-1+\frac{1}{2}(2a'-a-z+1)
        +\frac{1}{2}(c+z-2a')\\
& &\geqslant 3A(0,z-1)+2A(z,a)-1+\frac{1}{2}(c-a+1).
\end{eqnarray*}
\quad $\Box$(Claim \ref{ftlittle}.4)

\medskip

We now prove the lemma. The proof is divided into two cases.
The first case is easy and the second case is hard.

\medskip

{\bf Case \ref{ftlittle}.1}\quad $H-c\leqslant 2A(c+1,H)=2A(c,H)-2$.

Since $c-z+A(c,H)\succ A(z,c)+2A(c,H)-2$, then by Theorem \ref{A+B}
we have $|A[z,c]+A[c,H]|\geqslant A(z,c)+2A(c,H)-2$. Hence
\begin{eqnarray*}
\lefteqn{3|A|-3+b=|2A|}\\
& &\geqslant (2A)(0,c+z)-1+|A[z,c]+A[c,H]|+|H+A[c+1,H]|\\
& &\geqslant 3A(0,z-1)+2A(z,a)-1+\frac{1}{2}(c-a+1)-1\\
& &\quad +A(z,c)+2A(c,H)-2+A(c+1,H)\\
& &=3A(0,z-1)+3A(z,a)+3A(c,H)-4+\frac{1}{2}(c-a+1)\\
& &=3|A|-3+\frac{1}{2}(c-a-1).
\end{eqnarray*}
This shows $\frac{1}{2}(c-a-1)\leqslant b$. Hence
\begin{eqnarray*}
\lefteqn{H+1=H-c+c-a-1+a-z+z+2}\\
& &\leqslant 2A(c,H)-2+2b+2A(z,a)-1+2A(0,z-1)+2\\
& &=2|A|-1+2b.
\end{eqnarray*}
\quad $\Box$(Case \ref{ftlittle}.1)

\medskip

{\bf Case \ref{ftlittle}.2}\quad $H-c\geqslant 2A(c+1,H)+1=2A(c,H)-1$.

Note that Case \ref{ftlittle}.1 covers the case for $c=H$.
So we can assume $c<H$. First we prove a claim.

\medskip

{\bf Claim \ref{ftlittle}.2.1}\quad If
$(2A)(c+z,2H)\geqslant\frac{1}{2}(H-c)+A(c,H)+A(z,H)-1$, then
$H+1\leqslant 2|A|-1+2b$.

Proof of Claim \ref{ftlittle}.2.1\quad By the assumption we have
\begin{eqnarray*}
\lefteqn{3|A|-3+b=|2A|}\\
& &\geqslant 3A(0,z-1)+2A(z,a)-1+\frac{1}{2}(c-a+1)+(2A)(c+z,2H)-1\\
& &\geqslant 3A(0,z-1)+2A(z,a)-1+\frac{1}{2}(c-a+1)\\
& &\quad +\frac{1}{2}(H-c)+A(c,H)+A(z,H)-2\\
& &\geqslant 3|A|-3+\frac{1}{2}(H-c)-A(c,H)+\frac{1}{2}(c-a+1).
\end{eqnarray*}
Hence $\frac{1}{2}(H-c)-A(c,H)+\frac{1}{2}(c-a+1)\leqslant b$.
This implies
\begin{eqnarray*}
\lefteqn{H+1=H-c+c-a+1+a-z+z}\\
& &\leqslant 2(b+A(c,H))+2A(z,a)-1+2A(0,z-1)\\
& &=2|A|-1+2b.
\end{eqnarray*}
\quad $\Box$(Claim \ref{ftlittle}.2.1)

\medskip

By Claim \ref{ftlittle}.2.1 we need only to show that
$(2A)(c+z,2H)\geqslant\frac{1}{2}(H-c)+A(c,H)+A(z,H)-1$
is true. We divide the proof into cases according to the
structural properties of $A[c,H]$.

\medskip

{\bf Subcase \ref{ftlittle}.2.1}\quad $\gcd(A[c,H]-c)=1$.

Note that $c\not=H$ and $c\not= H-1$ by the condition of
Case \ref{ftlittle}.2. Since $\gcd(A[c,H]-c)=1$, then
$A(c,H)\geqslant 3$ and
$H-c\geqslant 2A(c+1,H)+1\geqslant 5$. Since
$\underline{d}_{U}(A)\geqslant\frac{2}{3}$, there is a $t\in U$
such that for all $u\geqslant t$ in $U$ we have
$\frac{A(t,u)}{u-t+1}\ggs\frac{2}{3}$. Let $u=t+H-c-1$.
Then there is a non-negative infinitesimal $r$ such that
$\frac{A(t,u)}{u-t+1}=\frac{A(t,u)}{H-c}\geqslant \frac{2}{3}-r$.
By Theorem \ref{A+B} we have
\begin{eqnarray*}
\lefteqn{(2A)(c+z,2H)\geqslant |A[c,H]+A[z,H]|}\\
& &\geqslant |c+A[z,t-1]|+|A[c,H]+A[t,u]|+|H+A[u+1,H]|\\
& &\geqslant A(z,t-1)+\min\{H-c+A(t,u),A(c,H)+2A(t,u)-2\}+A(u+1,H).
\end{eqnarray*}
Since
\begin{eqnarray*}
\lefteqn{H-c=\frac{1}{2}(H-c)+\frac{1}{2}(H-c)}\\
& &\geqslant \frac{1}{2}(H-c)+\frac{1}{2}(2A(c,H)-1)\\
& &=\frac{1}{2}(H-c)+A(c,H)-\frac{1}{2}
\end{eqnarray*}
and
\begin{eqnarray*}
\lefteqn{A(t,u)\geqslant(\frac{2}{3}-r)(H-c)}\\
& &\geqslant \frac{1}{2}(H-c)+(\frac{1}{6}-r)(H-c)\geqslant
   \frac{1}{2}(H-c)+\frac{5}{6}-5r,
\end{eqnarray*}
by the fact that $H-c\geqslant 5$, then we have
\begin{eqnarray*}
\lefteqn{\min\{H-c+A(t,u),A(c,H)+2A(t,u)-2\}}\\
& &\geqslant \frac{1}{2}(H-c)+A(c,H)+A(t,u)-\frac{7}{6}-5r.
\end{eqnarray*}
Hence
\begin{eqnarray*}
\lefteqn{(2A)(c+z,2H)\geqslant |A[c,H]+A[z,H]|}\\
& &\geqslant A(z,t-1)+\frac{1}{2}(H-c)+A(c,H)+A(t,u)\\
& &\quad -\frac{7}{6}-5r+A(u+1,H)\\
& &=A(z,H)+A(c,H)+\frac{1}{2}(H-c)-1-(\frac{1}{6}+5r).
\end{eqnarray*}
Since $(2A)(c+z,2H)$ is an integer, then we have
\[(2A)(c+z,2H)\geqslant A(z,H)+A(c,H)+\frac{1}{2}(H-c)-1.\]
Now the lemma follows from Claim \ref{ftlittle}.2.1.
\quad $\Box$(Subcase \ref{ftlittle}.2.1)

\medskip

{\bf Subcase \ref{ftlittle}.2.2}\quad $\gcd(A[c,H]-c)=d>1$ but
$d\not=3$.

Again let $t\in U\cap A$ be such that
$\frac{A(t,u)}{u-t+1}\ggs\frac{3}{2}$ 
for all $u\geqslant t$ in $U$.

\medskip

{\bf Claim \ref{ftlittle}.2.2.1}\quad For each $x\in A[c,H]$,
$(2A)(t+c,t+x-1)\geqslant A(c,x-1)+\frac{1}{2}(x-c)$.

Proof of Claim \ref{ftlittle}.2.2.1:\quad We prove the claim by
induction on $x\geqslant c$.

The case of $x=c$ is trivially true.

Suppose the claim is true for $y\in A[c,H-1]$. Let $x=\min
A[y+1,H]$. Since $[y+1,x-1]\cap A=\emptyset$, and $x-y=nd$ for some
$n>0$ implies $x-y=2$ or $x-y\geqslant 4$, then
\begin{eqnarray*}
\lefteqn{(2A)(t+y,t+x-1)\geqslant |y+A[t,t+(x-y)-1]|}\\
& &=A(t,t+(x-y)-1)\geqslant (\frac{2}{3}-r)(x-y)\\
& &=\frac{1}{2}(x-y)+(\frac{1}{6}-r)(x-y),
\end{eqnarray*}
for some non-negative infinitesimal $r$. Since either $x-y=2$
or $x-y>3$, and $(2A)(t+y,t+x-1)$ is an integer, then
we have
\[(2A)(t+y,t+x-1)\geqslant\frac{1}{2}(x-y)+1
=\frac{1}{2}(x-y)+A(y,x-1).\] Hence
\begin{eqnarray*}
\lefteqn{(2A)(t+c,t+x-1)}\\
& &=(2A)(t+c,t+y-1)+(2A)(t+y,t+x-1)\\
& &\geqslant A(c,y-1)+\frac{1}{2}(y-c)+\frac{1}{2}(x-y)+A(y,x-1)\\
& &=A(c,x-1)+\frac{1}{2}(x-c).
\end{eqnarray*}
\quad $\Box$(Claim \ref{ftlittle}.2.2.1)

\medskip

Following Claim \ref{ftlittle}.2.2.1 we now have
\begin{eqnarray*}
\lefteqn{(2A)(t+c,t+H)=(2A)(t+c,t+H-1)+1}\\
& &\geqslant A(c,H-1)+\frac{1}{2}(H-c)+1=A(c,H)+\frac{1}{2}(H-c).
\end{eqnarray*}
This implies
\begin{eqnarray*}
\lefteqn{(2A)(c+z,2H)\geqslant |A[c,H]+A[z,H]|}\\
& &\geqslant |c+A[z,t-1]|+(2A)(t+c,t+H)+|H+A[t+1,H]|\\
& &\geqslant A(z,t-1)+\frac{1}{2}(H-c)+A(c,H)+A(t+1,H)\\
& &=\frac{1}{2}(H-c)+A(c,H)+A(z,H)-1.
\end{eqnarray*}
Now the lemma follows from Claim \ref{ftlittle}.2.1.
\quad $\Box$(Subcase \ref{ftlittle}.2.2)

\medskip

{\bf Subcase \ref{ftlittle}.2.3}\quad $\gcd(A[c,H]-c)=3$.

Note that $t,t+1\in A$. Suppose $\{x\in A[z,a]:x-t\equiv 2\,(\mod 3)\}
=\emptyset$. If $c\in\{t,t+1\}\,(\mod 3)$, then $A[z,H]$ is a subset of
the \bp\ $(t\pm (3\ast\,^*\!\nat))\cup(t+1\pm (3\ast\,^*\!\nat))$.
Hence the lemma follows from Lemma \ref{bplittle}.
So we can assume that $c\not\in\{t,t+1\}\,(\mod 3)$. This implies
$(A[c,H]+A[c,H])\cap (A[c,H]+A[z,a])=\emptyset$ and hence
\begin{eqnarray*}
\lefteqn{|A[c,H]+A[z,H]|}\\
& &\geqslant |c+A[z,z+H-c-1]|+|H+A[z,a]|+|A[c,H]+A[c,H]|\\
& &\geqslant A(z,z+H-c-1)+A(z,a)+2A(c,H)-1\\
& &\geqslant\frac{1}{2}(H-c)+A(c,H)+A(z,H)-1.
\end{eqnarray*}
Now the lemma follows from Claim \ref{ftlittle}.2.1. So we can
assume \[\{x\in A[z,a]:x-t\equiv 2\,(\mod 3)\}\not =\emptyset.\]

Suppose $\{x\in A[t,u]:x-t\equiv 2\,(\mod 3)\}=\emptyset$, where
$u=t+H-c-1$.

If $\{x\in A[z,t-1]:x-t\equiv 2\,(\mod 3)\}\not =\emptyset$, let
\[k=\max\{x\in A[z,t-1]:x-t\equiv 2\,(\mod 3)\}.\] Then
$(k+A[c+1,H])\cap (c+A[z,t-1])=\emptyset$ and $(k+A[c+1,H])\cap
(A[c,H]+A[t,u])=\emptyset$. Hence
\begin{eqnarray*}
\lefteqn{(2A)(c+z,2H)}\\
& &\geqslant |c+A[z,t-1]|+|A[c,H]+A[t,u]|+|k+A[c+1,H]|+|H+A[u+1,H]|\\
& &\geqslant A(z,t-1)+|c+A[t,u]|+|H+A[t,u]|+A(c+1,H)+A(u+1,H)\\
& &\geqslant A(z,t-1)+2A(t,u)+A(c,H)-1+A(u+1,H)\\
& &\geqslant A(z,H)+A(c,H)+\frac{1}{2}(u-t+1)-1\\
& &=\frac{1}{2}(H-c)+A(c,H)+A(z,H)-1.
\end{eqnarray*}
Now the lemma follows from Claim \ref{ftlittle}.2.1.

If $\{x\in A[z,u]:x-t\equiv 2\,(\mod 3)\}=\emptyset$, let
$k=\min\{x\in A[u+1,a]:x-t\equiv 2\,(\mod 3)\}$. Then
$(k+A[c,H])\cap (A[c,H]+A[z,k-1])=\emptyset$. Hence
\begin{eqnarray*}
\lefteqn{(2A)(c+z,2H)\geqslant |A[c,H]+A[z,H]|}\\
& &\geqslant |c+A[z,t-1]|+|A[c,H]+A[t,u]|+|H+A[u+1,k-1]|\\
& &\quad +|k+A[c,H]|+|H+A[k+1,H]|\\
& &\geqslant A(z,t-1)+|c+A[t,u]|+|H+A[t,u]|\\
& &\quad +A(u+1,k-1)+A(c,H)+A(k+1,H)\\
& &\geqslant A(z,t-1)+2A(t,u)+A(u+1,k-1)+A(c,H)+A(k+1,H)\\
& &\geqslant A(z,H)-1+A(c,H)+\frac{1}{2}(u-t+1)\\
& &=\frac{1}{2}(H-c)+A(c,H)+A(z,H)-1.
\end{eqnarray*}
This implies the lemma.

Now we can assume that $\{x\in A[t,u]:x-t\equiv 2\,(\mod
3)\}\not=\emptyset$. For $i=0,1,2$ let
\[A_i=\{x\in A[t,u]:x-t\equiv i\,(\mod 3)\}.\]
Clearly $A_i\not=\emptyset$ for $i=0,1,2$. By Theorem \ref{A+B},
\[|A[c,H]+A_i|\geqslant\min\{\frac{1}{3}(H-c)+|A_i|,A(c,H)+2|A_i|-2\}.\]
Let $Q=|A[c,H]+A[t,u]|=\sum_{i=0}^2|A[c,H]+A_i|$. Note that each
term $|A[c,H]+A_i|$ in the sum has two possible lower bounds
$\frac{1}{3}(H-c)+|A_i|$ or $A(c,H)+2|A_i|-2$. We divide the proof
into the cases according to the different combinations of these
lower bounds of $|A[c,H]+A_i|$ for $i=0,1,2$.

\medskip

{\bf Subsubcase \ref{ftlittle}.2.3.1}\quad
$Q\geqslant\sum^2_{i=0}(\frac{1}{3}(H-c)+|A_i|)$.

Together with the assumption of Case \ref{ftlittle}.2, this subsubcase implies
\[Q\geqslant H-c+A(t,u)\geqslant\frac{1}{2}(H-c)+
A(c,H)-\frac{1}{2}+A(t,u).\] Hence
\begin{eqnarray*}
\lefteqn{(2A)(c+z,2H)\geqslant |A[c,H]+A[z,H]|}\\
& &\geqslant |c+A[z,t-1]|+Q+|H+A[u+1,H]|\\
& &\geqslant A(z,t-1)+\frac{1}{2}(H-c)+A(c,H)-\frac{1}{2}+A(t,u)+
   A(u+1,H)\\
& &\geqslant\frac{1}{2}(H-c)+A(c,H)+A(z,H)-1.
\end{eqnarray*}
Now the lemma follows from Claim \ref{ftlittle}.2.1.
\quad $\Box$(Subsubcase \ref{ftlittle}.2.3.1)

\medskip

{\bf Subsubcase \ref{ftlittle}.2.3.2}\quad
$Q\geqslant \sum_{i=0}^1(\frac{1}{3}(H-c)+|A_i|)+A(c,H)+2|A_2|-2$.

Then $Q\geqslant \frac{2}{3}(H-c)+A(t,u)+A(c,H)+|A_2|-2$. Hence
\begin{eqnarray*}
\lefteqn{(2A)(c+z,2H)\geqslant |A[c,H]+A[z,H]|}\\
& &\geqslant |c+A[z,t-1]|+Q+|H+A[u+1,H]|\\
& &\geqslant A(z,t-1)+\frac{2}{3}(H-c)+A(t,u)+A(c,H)+|A_2|-2+A(u+1,H)\\
& &>\frac{1}{2}(H-c)+A(c,H)+A(z,H)-1.
\end{eqnarray*}
\quad $\Box$(Subsubcase \ref{ftlittle}.2.3.2)

\medskip

{\bf Subsubcase \ref{ftlittle}.2.3.3}\quad
$Q\geqslant (\frac{1}{3}(H-c)+|A_0|)+\sum_{i=1}^2(A(c,H)+2|A_i|-2)$.

If $H-c=3$, then $A(c,H)=2$ and $A[t,u]=\{t,t+1,t+2\}$. Hence
$|A[c,H]+A_i|=2=\frac{1}{3}(H-c)+|A_i|$, which implies the
assumption of Subsubcase \ref{ftlittle}.2.3.1. So we can assume
$H-c\geqslant 6$ and $A(c,H)\geqslant 3$. Since $A(t,u)\geqslant 
(\frac{2}{3}-r)(H-c)$ and $|A_0|\leqslant\frac{1}{3}(H-c)$, then
$|A_1|+|A_2|\geqslant(\frac{1}{3}-r)(H-c)$ for some non-negative
infinitesimal $r$. Hence
\begin{eqnarray*}
\lefteqn{Q\geqslant\frac{1}{3}(H-c)+A(t,u)+2A(c,H)+|A_1|+|A_2|-4}\\
& &\geqslant\frac{1}{2}(H-c)+(\frac{1}{6}-r)(H-c)
   +A(t,u)+2A(c,H)-4\\
& &\geqslant\frac{1}{2}(H-c)+(\frac{1}{6}-r)(H-c)
   +A(t,u)+A(c,H)-1\\
& &>\frac{1}{2}(H-c)+A(t,u)+A(c,H)-1.
\end{eqnarray*}
Hence again
\[(2A)(c+z,2H)\geqslant |A[c,H]+A[z,H]|
\geqslant\frac{1}{2}(H-c)+A(c,H)+A(z,H)-1.\]
\quad $\Box$(Subsubcase \ref{ftlittle}.2.3.3)

\medskip

{\bf Subsubcase \ref{ftlittle}.2.3.4}\quad
$Q\geqslant\sum_{i=0}^2(A(c,H)+2|A_i|-2)$.

Again we can assume $H-c\geqslant 6$ and $A(c,H)\geqslant 3$. Then
\begin{eqnarray*}
\lefteqn{Q\geqslant 3A(c,H)+2A(t,u)-6}\\
& &\geqslant A(t,u)+\frac{1}{2}(H-c)+A(c,H),
\end{eqnarray*}
for some non-negative infinitesimal $r$.
Hence $Q>A(t,u)+\frac{1}{2}(H-c)+A(c,H)-1$. This implies
\[(2A)(c+z,2H)\geqslant |A[c,H]+A[0,H]|
\geqslant\frac{1}{2}(H-c)+A(c,H)+A(z,H)-1,\]
which again implies the lemma. The rest of the cases can be proven by
symmetry of the proofs above. 
\quad $\Box$(Lemma \ref{ftlittle})

\begin{lemma}\label{btft}

Suppose $0\prec s\prec H$ such that
$A[0,s]$ is a backward triangle from $0$ to $s$ and $A[s+1,H]$
is a forward triangle from $s+1$ to $H$. Then
$H+1\leqslant 2|A|-1+2b$.

\end{lemma}

\noindent {\bf Proof}:\quad Let
\[u=\min\{x:\frac{s}{2}<x<\frac{s+H}{2}\mbox{ and }
A(0,x)\geqslant\frac{1}{2}(x+1)\}.\]
Clearly $u\succeq s$ because $A[0,s]$ is a backward triangle.
Also $u\preceq s$ because otherwise 
\[A(0,u)\sim A(0,s)+A(s,u)\succ\frac{1}{2}s+\frac{1}{2}(u-s)
\sim\frac{1}{2}u.\] Hence we have $u\sim s$.
It is easy to see that $u\in A$ and $A(0,u)=\frac{1}{2}(u+1)$
by the minimality of $u$. Also by the minimality of $u$
we have that for any $0\prec x\leqslant u$,
$A(x,u)>\frac{1}{2}(u-x+1)$. Let
\[X=\{x:u+1\leqslant x<\frac{u+H}{2}\mbox{ and }
A(u+1,x)\leqslant\frac{1}{2}(x-u)\}.\]
If $X\not=\emptyset$, let $z=1+\max X$. Otherwise let $z=u+1$.
It is also easy to see that $z\in A$, $z\sim u$, and
$A(u+1,z-1)=\frac{1}{2}(z-u-1)$. Since
$A(0,z-1)=\frac{1}{2}z$, $A(0,H)\leqslant\frac{1}{2}(H+1)$,
and $H\in A$, then the number $a$ below is well defined.
\[a=\min\{x:z\leqslant x<H\mbox{ and }
A(z,x)\leqslant\frac{1}{2}(x-z+1)\}.\]
Clearly $a<H$, $a\sim H$, $a\not\in A$, and
$A(z,a)=\frac{1}{2}(a-z+1)$. By the minimality we have that
for any $z\leqslant x<a$, $A(z,x)>\frac{1}{2}(x-z+1)$.
Now let $a'=\max (A[z,a-1])$
and $c=\min (A[a,H])$. Since $a'\succeq\frac{z+H}{2}$ and
$z\succ 0$, then $2a'\succ c$. Let
$S=(2A)[0,z-1]\smallsetminus A[0,z-1]$.

Without loss of generality we can assume that $A[z,H]$ is not
a subset of a \bp\ of difference $3$ by the following reason:

Suppose not. By symmetry, we can also assume that there is a
$z'\sim z$ such that $A[0,z']$ is a subset of a \bp\ of difference
$3$. Now the lemma follows from Lemma \ref{bplittle3}.

The rest of the proofs are almost identical to the proofs of 
Lemma \ref{ftlittle}. We will refer to the proofs of Lemma \ref{ftlittle}
when the steps are the same and add more proofs when
the steps are not the same.

\medskip

{\bf Claim \ref{btft}.1}:\quad $[z,a+z-1]\subseteq (2A)$.

Proof of Claim \ref{btft}.1:\quad The proof here is slightly
different from the proof of Claim \ref{ftlittle}.1.

If $2z\leqslant x<a+z$, then
$z\leqslant x-z<a$. Hence $A(z,x-z)>\frac{1}{2}(x-2z+1)$.
This implies $A[z,x-z]\cap (x-A[z,x-z])\not=\emptyset$. Hence
$x\in (2A)$. Suppose $z\leqslant x<2z$. Then $0\leqslant x-z<z$.
If $0\prec x-z\leqslant u$, then $A(x-z,z)=A(x-z,u)+A(u+1,z)
>\frac{1}{2}(u-x+z+1)+\frac{1}{2}(z-u-1)+1>\frac{1}{2}(2z-x+1)$.
Hence $A[x-z,z]\cap (x-A[x-z,z])\not=\emptyset$, which implies
$x\in (2A)$. If $u<x-z<z$, then by choosing a $y\succ 0$ with
$y\prec\min\{a-z,u\}$ we have $A(x-z-y,x-z)+A(z,z+y)\succ y+1$.
Hence $A[x-z-y,x-z]\cap (x-A[z,z+y])\not=\emptyset$, which
implies $x\in (2A)$. Suppose $x-z\sim 0$. If $A(0,x-z)\geqslant
\frac{1}{2}(x-z+1)$, then by $A(z,x)>\frac{1}{2}(x-z+1)$ we have
$A[0,x-z]\cap (x-A[z,x])\not=\emptyset$,
which implies $x\in (2A)$. If $A(0,x-z)<\frac{1}{2}(x-z+1)$, then
$A(x-z+1,z-1)>\frac{1}{2}(2z-x-1)$. Hence $A[x-z+1,z-1]\cap 
(x-A[x-z+1,z-1])\not=\emptyset$, which again implies $x\in (2A)$.
\quad $\Box$(Claim \ref{btft}.1)

\medskip

{\bf Claim \ref{btft}.2}:\quad Suppose $2a'>a+z$ and
$a+z<x\leqslant\min\{2a',c+z\}$. If $(2A)(a+z,x-1)<\frac{1}{2}(x-a-z)$,
then $x\in 2A$.

Proof of Claim \ref{btft}.2:\quad The proof is identical
to the proof of Claim \ref{ftlittle}.2.
\quad $\Box$(Claim \ref{btft}.2)

\medskip

{\bf Claim \ref{btft}.3}:\quad Suppose $2a'<c+z$ and
$\max\{a+z,2a'\}\leqslant x< c+z-1$. If $(2A)(x+1,c+z-1)
<\frac{1}{2}(c+z-x-1)$, then $x\in (2A)$ or $x-c\in S$.

Proof of Claim \ref{btft}.3:\quad Identical to the proof
of Claim \ref{ftlittle}.3.
\quad $\Box$(Claim \ref{btft}.3)

\medskip

{\bf Claim \ref{btft}.4}:\quad $(2A)(0,c+z)\geqslant
3A(0,z-1)+2A(z,a)-1+\frac{1}{2}(c-a+1)$.

Proof of Claim \ref{btft}.4:\quad The proof is divided into
three cases for $2a'\geqslant c+z$, $2a'\leqslant a+z$, and
$a+z<2a'<c+z$.

\medskip

{\bf Case \ref{btft}.4.1}:\quad $2a'\geqslant c+z$.

Identical to the proof of Case \ref{ftlittle}.4.1.
\quad $\Box$(Case \ref{btft}.4.1)

\medskip

{\bf Case \ref{btft}.4.2}:\quad $2a'\leqslant a+z$.

Identical to the proof of Case \ref{ftlittle}.4.2.
\quad $\Box$(Case \ref{btft}.4.2)

\medskip

{\bf Case \ref{btft}.4.3}:\quad $a+z<2a'<c+z$.

Identical to the proof of Case \ref{ftlittle}.4.3.
\quad $\Box$(Case \ref{btft}.4.3)

\medskip

Now we prove the lemma. The proof is divided into two cases.

\medskip

{\bf Case \ref{btft}.1}\quad $H-c\leqslant 2A(c+1,H)=2A(c,H)-2$.

Identical to the proof of Case \ref{ftlittle}.1.
\quad $\Box$(Case \ref{btft}.1)

\medskip

{\bf Case \ref{btft}.2}\quad $H-c\geqslant 2A(c+1,H)+1=2A(c,H)-1$.

\medskip

{\bf Claim \ref{btft}.2.1}\quad If
$(2A)(c+z,2H)\geqslant\frac{1}{2}(H-c)+A(c,H)+A(z,H)-1$, then
$H+1\leqslant 2|A|-1+2b$.

Proof of Claim \ref{btft}.2.1\quad Identical to 
the proof of Claim \ref{ftlittle}.2.1.
\quad $\Box$(Claim \ref{btft}.2.1)

\medskip

By Claim \ref{btft}.2.1 we need only to show that
$(2A)(c+z,2H)\geqslant\frac{1}{2}(H-c)+A(c,H)+A(z,H)-1$
is true. We divide the proof into cases according to the
structural properties of $A[c,H]$.

\medskip

{\bf Subcase \ref{btft}.2.1}\quad $\gcd(A[c,H]-c)=1$.

The proof is the same as the proof of 
Subcase \ref{ftlittle}.2.1 except that the term $U$ needs to
be replaced by the term $z+U$ throughout the remaining of the 
proof of the lemma. Note that we can also
assume that $\underline{d}_{z+U}\geqslant\frac{2}{3}$
by the same reason as stated at the beginning of the proof
of Lemma \ref{ftlittle}.
\quad $\Box$(Subcase \ref{btft}.2.1)

\medskip

{\bf Subcase \ref{btft}.2.2}\quad $\gcd(A[c,H]-c)=d>1$ but
$d\not=3$.

\medskip

{\bf Claim \ref{btft}.2.2.1}\quad For each $x\in A[c,H]$,
$(2A)(t+c,t+x-1)\geqslant A(c,x-1)+\frac{1}{2}(x-c)$.

Proof of Claim \ref{btft}.2.2.1:\quad 
Identical to the proof of Claim \ref{ftlittle}.2.2.1.
\quad $\Box$(Claim \ref{btft}.2.2.1)

\medskip

Following Claim \ref{btft}.2.2.1 we now have
\begin{eqnarray*}
\lefteqn{(2A)(t+c,t+H)=(2A)(t+c,t+H-1)+1}\\
& &\geqslant A(c,H-1)+\frac{1}{2}(H-c)+1=A(c,H)+\frac{1}{2}(H-c).
\end{eqnarray*}
This implies
\begin{eqnarray*}
\lefteqn{(2A)(c+z,2H)\geqslant |A[c,H]+A[z,H]|}\\
& &\geqslant |c+A[z,t-1]|+(2A)(t+c,t+H)+|H+A[t+1,H]|\\
& &\geqslant A(z,t-1)+\frac{1}{2}(H-c)+A(c,H)+A(t+1,H)\\
& &=\frac{1}{2}(H-c)+A(c,H)+A(z,H)-1.
\end{eqnarray*}
Now the lemma follows from Claim \ref{btft}.2.1.
\quad $\Box$(Subcase \ref{btft}.2.2)

\medskip

{\bf Subcase \ref{btft}.2.3}\quad $\gcd(A[c,H]-c)=3$.

the proof is identical to the proof of Subcase \ref{ftlittle}.2.3.
Note that we assume $\{x\in A[z,a]:x-t\equiv 2\,(\mod 3)\}\not
=\emptyset$ in the beginning of the proof of this lemma. 
\quad $\Box$(Lemma \ref{btft})

\begin{lemma}\label{btap}

Suppose $A=A[0,s]\cup A[s+1,H]$ with $0\prec s\prec H$ such that $A[0,s]$ is a
backward triangle and $A[s+1,H]$ is a subset of an \ap\ of difference $d>1$. Then
$H+1\leqslant 2|A|-1+2b$.

\end{lemma}

\noindent {\bf Proof}:\quad Since (\ref{half}), we have 
$A(s+1,H)\sim\frac{1}{2}(H-s)$, which implies $d=2$. Let
$E$ be the set of all even numbers and let $c=\min A[s+1,H]$. Then $c\sim s$ and
$A[c,H]$ is full. Without loss of generality we can assume that $s\in A$ and
$c-s$ is odd. By the pigeonhole principle and Lemma \ref{overspill}
we can find $e\sim 2c$ and $e'\sim 2H$ such that 
$E[e,e']=(A[c,H]+A[c,H])[e,e']$.

\medskip

{\bf Claim \ref{btap}.1}:\quad If $s\prec x\prec s+H$, then $x\in 2A$.

Proof of Claim \ref{btap}.1:\quad If $s\prec x\prec 2s$, then $0\prec x-s\prec s$.
Hence $A(x-s,s)\succ\frac{1}{2}(2s-x+1)$. So $A[x-s,s]\cap (x-A[x-s,s])\not=\emptyset$,
which implies $x\in 2A$. Now we assume $2s\preceq x\prec s+H$. Let $z_1=[\frac{x}{2}]$. 
Then $2z_1\sim x\prec s+H$ implies $z_1-s\prec H-z_1$.
Choose a $y$ with $z_1-s\prec y\prec\min\{H-z_1,z_1\}$. Then $0\prec z_1-y\prec s$.
Hence 
\begin{eqnarray*}
\lefteqn{A(z_1-y,z_1)\sim A(z_1-y,s)+A(s+1,z_1)}\\
& &\succ\frac{1}{2}(s-z_1+y+1)+\frac{1}{2}(z_1-s)=\frac{1}{2}(y+1)
\end{eqnarray*}
and $A(x-z_1,x-z_1+y)\sim\frac{1}{2}(y+1)$. Hence
$A[z_1-y,z_1]\cap (x-A[x-z_1,x-z_1+y])\not=\emptyset$, which implies $x\in 2A$.
\quad $\Box$(Claim \ref{btap}.1)

\medskip

Now let's assume that the lemma is not true. Let $A'\subseteq [0,H]$ be the set
with the largest cardinality $|A'|$ such that  
$A[s+1,H]\subseteq A'[s+1,H]\subseteq (s+1+E)$, $A'[0,s]=A[0,s]$,
and satisfying (\ref{cond3k-3+b}), (\ref{smallA}), and (\ref{longap})
with $A$ replaced by $A'$ and $b$ replaced by $b'$. We will derive a contradiction.

\medskip

{\bf Claim \ref{btap}.2}:\quad $A'[s+1,H]=(s+1+E)[s+1,H]$.

Proof of Claim \ref{btap}.2:\quad The proof is divided into three cases. Let
$x\in (s+1+E)$ be such that $s+1\leqslant x\leqslant H$. We want to show that
$x\in A'$.

\medskip

{\bf Case \ref{btap}.2.1}\quad $s\prec x\prec\frac{s+H}{2}$.

Suppose $x\not\in A'$. Let $A''=A'\cup\{x\}$. For each $y\in A'[s+1,H]\cup\{x\}$,
$x+y\in E[e,e']\subseteq 2A\subseteq 2A'$, and for each $y\in A'[0,s]$ we have
$s\prec x+y\prec s+\frac{s+H}{2}$, which implies $s\prec x+y\prec s+H$. Hence
$x+y\in 2A\subseteq 2A'$ by Claim \ref{btap}.1. So $2A''=2A'$, which contradict
the maximality of $|A'|$ by Lemma \ref{addone}.
\quad $\Box$(Case \ref{btap}.2.1)

\medskip

{\bf Case \ref{btap}.2.2}:\quad $x\geqslant s+1$ and $x\sim s+1$.

Suppose $x\not\in A'$. Without loss of generality we can, 
by Case \ref{btap}.2.1, assume
\[x=\max((s+1+E)[s+1,\frac{3s+H}{4}]\smallsetminus A').\] 
Let $y\in A''=A'\cup\{x\}$.

If $0\prec y\prec s$, then $A(y+1,s)\succ\frac{1}{2}(s-y)$ and $A(x-(s-y),x-1)
\succ\frac{1}{2}(s-y)$. Hence $A[y+1,s]\cap (x+y-A[x-(s-y),x-1])\not=\emptyset$,
which implies $x+y\in A[y+1,s]+A[x-(s-y),x-1]\subseteq (2A')$.

If $y\sim s$ and $y\leqslant x$, then choose a $z\prec y$ such that
$y-z\prec H-s$. Hence $A(z,y-1)\succ\frac{1}{2}(y-z)$ and
$A(x+1,x+(y-z))\sim\frac{1}{2}(y-z)$ imply $A[z,y-1]\cap (x+y-A[x+1,x+(y-z)])\not=
\emptyset$, which implies $x+y\in (2A)\subseteq (2A')$.

If $x<y\prec H$, then choose a $z$ with $H\geqslant z\succ y$ 
such that $z-y\prec s$. Hence $A(y+1,z)\sim\frac{1}{2}(z-y)$ 
and $A(x-(z-y),x-1)\succ\frac{1}{2}(z-y)$ imply
$A[y+1,z]\cap (x+y-A[x-(z-y),x-1])\not=\emptyset$, which implies
$x+y\in (2A')$.

If $y\sim H$, then $y\in (s+1+E)$. Hence $x+y$ is even and $2s\prec x+y\prec 2H$.
Now we have $x+y\in E[e,e']\subseteq (2A')$. 

If $y\sim 0$, $y>0$, and $y$ is even, then $x+y\in A'$ by the maximality of $x$.
Hence $x+y=0+(x+y)\in (2A')$.

If $y\sim 0$, $y$ is odd, and $y>l=\min\{z\in A':z\mbox{ is odd }\}$, then
$x+(y-l)\in A'$. Hence $x+y=l+(x+y-l)\in (2A')$.

By the arguments above we conclude that $(2A'')\smallsetminus (2A')\subseteq\{x,x+l\}$.
Hence $|2A''|\leqslant |2A'|+2$, which contradicts
the maximality of $|A'|$ by Lemma \ref{addone}. So we conclude that $x\in A'$.
\quad $\Box$(Case \ref{btap}.2.2)

\medskip

{\bf Case \ref{btap}.2.3}:\quad $\frac{s+H}{2}\preceq x\leqslant H$.

Suppose again $x\not\in A'$. Let $x=\min(A'[\frac{3s+H}{4},H])$.
By Case \ref{btap}.2.1 we have $x\succeq\frac{s+H}{2}$. Let $y\in A''$.

If $s+1\leqslant y\prec H$, then $x+y\in E[e,e']\subseteq (2A')$.

If $0\leqslant y\prec s$, then $s\prec x+y\prec s+H$, which implies, by
Claim \ref{btap}.1, $x+y\in 2A'$.

If $y<H$ and $y\sim H$, then $x-(H-y)\in A'$ by the minimality of $x$. Hence
$x+y=H+(x-H+y)\in 2A'$.

If $y<s$, $y\sim s$, and $s-y$ is odd, then $s+1-y$ is even and
$x-(s+1-y)\in A'$ by the minimality of $x$. Hence $x+y=s+1+(x-s-1+y)\in
2A'$.

If $y<s$, $y\sim s$, and $s-y$ is even, then $x-(s-y)\in A'$. Hence $x+y=s+(x-s+y)\in 2A'$.

By the arguments above we conclude that $(2A'')\smallsetminus (2A')\subseteq\{s+x,x+H\}$,
which contradicts the maximality of $|A'|$ by Lemma \ref{addone}.
\quad $\Box$(Claim \ref{btap}.2)

\medskip

Now we are ready to prove the lemma. Without loss of generality we can assume that the
set $A$ is already in the form of the set $A'$ in Claim \ref{btap}.2.
Let $u=\min(\{z\in A[s+1,H]:A[0,z]\mbox{ is not a subset of a \bp }\})$.
Note that if there is a $v\succ s$ such that $A[0,v]$ is a subset of a \bp\ of
difference $d$, then $|2A|\sim 3|A|$ implies that
$A[0,v]$ is full in the \bp\, Hence $d=1$ or $d=3$ because $A[0,s]$ 
is a backward triangle. However, $A[s+1,v]$ is a full subset of an \ap\ of
difference $2$, which contradicts $d=1$ or $d=3$. Hence we have $u\sim s$.
By (\ref{smallA}) and $A(u,H)=\frac{1}{2}(H-u)+1$ we have
$A(0,u)=|A|-A(u,H)+1\leqslant\frac{H+1}{2}-\frac{H-u}{2}=\frac{1}{2}(u+1)$.
Let $|A[0,u]+A[0,u]|=3A(0,u)-3+\bar{b}$. If $\bar{b}<0$, then by Theorem \ref{2k-1+b}
we have $u+1\leqslant 2A(0,u)-2+\bar{b}<2A(0,u)-2$, which contradicts
$A(0,u)\leqslant\frac{1}{2}(u+1)$. Hence we can assume $\bar{b}\geqslant 0$.
Clearly $\bar{b}\sim 0$ because otherwise we would have $|2A|\succ 3|A|$.
By Lemma \ref{ftlittle}, we have $u+1\leqslant 2A(0,u)-1+2\bar{b}$. Hence
\begin{eqnarray*}
\lefteqn{3|A|-3+b=|2A|}\\
& &\geqslant |A[0,u]+A[0,u]|+|s+A[u+2,H]|+|A[u+2,H]+A[u,H]|\\
& &\geqslant 3A(0,u)-3+\bar{b}+A(u+2,H)+2A(u+2,H)\\
& &=3|A|-3+\bar{b}.
\end{eqnarray*}
Above shows $\bar{b}\leqslant b$. Hence
\begin{eqnarray*}
\lefteqn{H+1=H-u+u+1}\\
& &\leqslant 2A(u+2,H)+2A(0,u)-1+2\bar{b}\\
& &=2|A|-1+2\bar{b}\leqslant 2|A|-1+2b,
\end{eqnarray*}
which contradicts $H+1>2|A|-1+2b$.\quad $\Box$(Lemma \ref{btap})

\begin{lemma}\label{largeu}

If $\underline{d}_U(A)>\frac{1}{2}$, then $H+1\leqslant 2|A|-1+2b$.

\end{lemma}

\noindent {\bf Proof}:\quad If $\underline{d}_U(A)>\frac{1}{2}$,
then there exists a $x\succ 0$ such that $A[0,x]$ is a forward
triangle. If $x\sim H$, then the lemma now follows from Lemma
\ref{ftlittle}. If $x\prec H$, then the lemma follows from
Lemma \ref{ftthm}. \quad $\Box$(Lemma \ref{largeu})

\begin{lemma}\label{middleu1}

Let $\underline{d}_U(A)=\frac{1}{2}$. If there is a $x\succ 0$ in
$A$ such that $\gcd(A[x,H]-x)=1$, then $A\cap U$ is either a
subset of an \ap\ of difference $>1$ or a subset of a
$U$--unbounded \bp

\end{lemma}

\noindent {\bf Proof}:\quad The lemma follows from Lemma \ref{smallbig2}
and Lemma \ref{kneser}. \quad $\Box$(Lemma \ref{middleu1})

\begin{lemma}\label{middleu4}

Let $\underline{d}_U(A)=\frac{1}{2}$. If $A\cap U$ is a subset of
a $U$--unbounded \bp, then $H+1\leqslant 2|A|-1+2b$.

\end{lemma}

\noindent {\bf Proof}:\quad Suppose $A\cap U$ is a subset of a
$U$--unbounded \bp\ of difference $d$. Since $\underline{d}_U(A)=\frac{1}{2}$,
then $d=3$ or $d=4$. Let $I_0=d\ast\,^*\!\nat$ and $I_1=c+(d\ast\,^*\!\nat)$
where $c=\min\{z\in A:z\not\equiv 0\,(\mod d)\}$. Then $A\cap U\subseteq 
I_0\cup I_1$. Since $d=3$ or $d=4$, then $\gcd(c,d)=1$.
By Lemma \ref{twoaps} we can assume that there is an
$a\succ 0$ in $A$ such that $\gcd(A[a,H]-a)=1$.

\medskip

{\bf Case \ref{middleu4}.1}:\quad $d=3$.

We want to show this case implies $|2A|\succ 3|A|$, hence $d=3$ is
impossible.

By Lemma \ref{smallbig2} it suffices to show that for 
$\gamma=\frac{1}{25}$ and for every $N\succ 0$ there is a $K\in A$
with $0\prec K\leqslant N$ such that (\ref{3timesmore}) is true.
Suppose $N\succ 0$ is given. Let $0<\epsilon<\frac{1}{12}$.
Without loss of generality we can re-choose
$a$ so that $a\leqslant N$, $A[0,a]\subseteq I_0\cup I_1$, and
for every $0\prec y\leqslant a$ we have $A(0,y)\succeq 
(\frac{1}{2}-\epsilon)y$. Choose a $x$ with $0\prec x\prec\frac{1}{2}a$ 
such that $A(0,x)\preceq (\frac{1}{2}+\epsilon)x$.
Let $K=\min\{z\geqslant x:z,z-1\in A\}$. It is easy to see that
$A(x,K)\preceq\frac{1}{3}(K-x)$ because for any two consecutive numbers
in $(I_0\cup I_1)[x,K-1]$, at least one of them is not in $A$.
Hence we have that $K\preceq 2x\prec a$ and
$A(0,K)\preceq (\frac{1}{2}+\epsilon)K$. Note that
$K,K-1\in A$. So the shortest \bp\ containing $A[0,K]$ must have
length $L\sim\frac{2}{3}K$. Let $A_i=A[0,K]\cap I_i$ for $i=0,1$.
Since $|A_0|\preceq\frac{1}{3}K$, then $|A_1|=A(0,K)-|A_0|\succeq
(\frac{1}{6}-\epsilon)K$. By the same reason we have $|A_0|\succeq
(\frac{1}{6}-\epsilon)K$. Let $|A[0,K]+A[0,K]|=3A(0,K)-3+b$ for some
integer $b$. Since $A[0,K]$ is a subset of a \bp, then $b\geqslant 0$.

If $b\geqslant\frac{1}{3}A(0,K)-3$, then 
\[(2A)(0,2K)\succeq 3A(0,K)+\frac{1}{3}A(0,K)\succeq 3A(0,K)+\frac{1}{3}
(\frac{1}{2}-\epsilon)K,\]
which implies (\ref{3timesmore}).

If $b<\frac{1}{3}A(0,K)-3$, then by Theorem \ref{twolines} we have that
\[\frac{2}{3}K\sim L\leqslant A(0,K)+b\preceq (\frac{1}{2}+\epsilon)K+b.\]
Hence $b\succeq (\frac{1}{6}-\epsilon)K$. So we have
\[(2A)(0,2K)\succeq 3A(0,K)+b\succ 3A(0,K)+\frac{1}{12}K,\]
which again implies (\ref{3timesmore}).
This ends the proof.\quad $\Box$(Case \ref{middleu4}.1)

\medskip

{\bf Case \ref{middleu4}.2}:\quad $d=4$.

Without loss of generality we assume $0\in I_0$ and $1\in I_1$. Suppose 
$A$ is not a subset of a \bp\, We want to derive a contradiction. Let
\[c=\min\{z\in [0,H]:A[0,z]\mbox{ is not a subset of a \bp\ of difference }4\}.\]
Let $A[0,c-1]=A_0\cup A_1$ where $A_i=A[0,c-1]\cap I_i$ for $i=0,1$. 
Then $|A_i|\succ 0$ for $i=0,1$ because otherwise 
$\underline{d}_U(A)\leqslant\frac{1}{4}$. Note that since $d=4$, then
there is an $i=0$ or $i=1$ such that $(c+A_i)\cap |A[0,c-1]+A[0,c-1]|
=\emptyset$. Hence we can assume $c\prec H$ because otherwise
$|2A|\succeq 3|A|+|c+A_i|\succ 3|A|$.

\medskip

{\bf Subcase \ref{middleu4}.2.1}:\quad $A(c,H)\succ\frac{1}{2}(H-c)$.

By Lemma \ref{triangle}
there exist $x\prec c\prec y\leqslant H$ such that $A[x,y]$ is a
backward triangle and $A(y,H)\sim\frac{1}{2}(H-y+1)$. 

If $x\succ 0$, then the lemma follows from either Lemma \ref{ftthm} 
or Lemma \ref{moreftmore2}. so we can assume $x\sim 0$.

If $y\sim H$, then the lemma follows from Lemma \ref{ftlittle}.
So we can assume $y\prec H$. If for any $y\prec y'\prec H$ in $A$, $\gcd(A[y',H]-y')>1$,
then the lemma follows from Lemma \ref{btap}. So we can assume that there is a $y'\in A$,
$y\prec y'\prec H$ such that $\gcd(A[y',H]-y')=1$.

If $\underline{d}_{y+U}(A)>\frac{1}{2}$, then by Lemma \ref{triangle} there is a 
$y\prec z\leqslant H$ such that $A[y,z]$ is a forward triangle. Now the lemma follows
from Lemma \ref{btft} if $z\sim H$. If $z\prec H$, then by Lemma \ref{ftthm} 
$|2A|\sim 3|A|$ implies that $A[y,H]$
is a full subset of a \bp\ $[y,z']\cup [z'',H]$. Hence $A[y,H]$ is the union of a forward
triangle $A[y,2z'-y]$ and a backward triangle $A[2z'-y+1,H]$. Now the lemma follows from
Lemma \ref{ftbtmore}.

If $\underline{d}_{y+U}<\frac{1}{2}$, then by Lemma \ref{triangle} there are 
$y\leqslant z\prec z'\leqslant H$ such that $A[z,z']$ is a backward triangle
and $A(z',H)\sim\frac{1}{2}(H-z')$. Now the lemma follows from Lemma \ref{ftthm}
because $A[0,z']$ cannot be a full subset of a \bp\ of difference $1$ or $3$.

Assume $\underline{d}_{y+U}(A)=\frac{1}{2}$. Since $A[0,y]$ is a backward triangle,
then we can assume $A(y,y+3)\geqslant 3$. Hence $(A-y)\cap U$ is neither a subset of
an \ap\ of difference $>1$ nor a subset of a $U$--unbounded \bp\ of difference $d\not=3$.
If $(A-y)\cap U$ is a subset of a $U$--unbounded \bp\ of difference $3$, then by the proof
of Case \ref{middleu4}.1 we have $|A[y,H]+A[y,H]|\succ 3A(y,H)$, which implies
$|2A|\succ 3|A|$. Note that if there is an $y'\sim y$ in $A$ such that
$\gcd(A[y',H]-y')=d'>1$, then $d'=2$ and the lemma follows from Lemma \ref{btap}.
So we can assume that $(A-y)\cap U$ is neither a subset of an \ap\ of difference $>1$
nor a subset of a $U$--unbounded \bp\ and there is a $y'\succ y$ in $A$ such that
$\gcd(A[y',H]-y')=1$. Now the lemma follows from Lemma \ref{kneser}, 
Lemma \ref{smallbig2}, and Lemma \ref{smallbig}.\quad $\Box$(Subcase \ref{middleu4}.2.1)

\medskip

{\bf Subcase \ref{middleu4}.2.2}:\quad $A(c,H)\preceq\frac{1}{2}(H-c)$. 

By (\ref{half}) we have $A(c,H)\sim\frac{1}{2}(H-c)$.
Since $|A_0\cup A_1|=A(0,c-1)\sim\frac{1}{4}c+\frac{1}{4}c$ and
$\gcd(A_0)=\gcd(A_1-1)=4$, then $A_i$ is full for $i=0,1$. 
Hence we can find a $c'\sim c$ in $A$ such that
$\gcd(A[c',H]-c')=1$. Since there is an $i\in\{0,1\}$ such that
$(2A)(0,2c)\succeq 3A(0,c-1)+|c+A_i|\succ 3A(0,c)$,
then by Lemma \ref{smallbig} we have $|2A|\succ 3|A|$. 
\quad $\Box$(Lemma \ref{middleu4})

\begin{lemma}\label{middleu3}

Let $\underline{d}_U(A)=\frac{1}{2}$. If $A\cap U$ is a subset of
an \ap\ of difference $>1$, then $H+1\leqslant 2|A|-1+2b$.

\end{lemma}

\noindent {\bf Proof}:\quad Let $l_o=\min\{x\in
A:\gcd(A[0,x])=1\}$. Since $\underline{d}_U(A)=\frac{1}{2}$, then
$\gcd(A[0,l_o-1])=2$ and $l_o$ is odd.

\medskip

{\bf Case \ref{middleu3}.1}:\quad There are $0\prec
a\prec c\leqslant H$ such that $A[a,c]$ is a backward triangle and
$A(c,H)\sim\frac{1}{2}(H-c)$.

Note that $l_o\prec c$. By Lemma \ref{ftthm} we have either
$|A[0,c]+A[0,c]|\succ 3A(0,c)$, which is impossible because it
implies $|2A|\succ 3|A|$ by Lemma \ref{smallbig}, or $A[0,c]$ is a
full subset of a \bp\ $[0,x]\cup [x',c]$, which contradicts 
$\underline{d}_U(A)=\frac{1}{2}$, or $A[0,c]$ is a full subset of
a \bp\ of difference $3$, which again contradicts 
$\underline{d}_U(A)=\frac{1}{2}$.\quad $\Box$(Case \ref{middleu3}.1)

\medskip

{\bf Case \ref{middleu3}.2}:\quad
$A(l_o,H)\succ\frac{1}{2}(H-l_o)$.

By lemma \ref{triangle} there are $0\leqslant y\prec l_o\prec
y'\leqslant H$ such that $A(y',H)\sim\frac{1}{2}(H-y')$ and
$A[y,y']$ is a backward triangle. Without loss of generality
we can assume that $A(y',y'+3)\geqslant 3$. By Case \ref{middleu3}.1 
we can assume $y\sim 0$ and by Lemma \ref{ftlittle} we can assume
$y'\prec H$.

\medskip

{\bf Subcase \ref{middleu3}.2.1}:\quad
$\underline{d}_{y'+U}(A)>\frac{1}{2}$.

By Lemma \ref{triangle} there is a $z\succ y'$ such that $A[y',z]$
is a forward triangle. If $z\sim H$, then the lemma follows from
Lemma \ref{btft}. So we can assume $z\prec H$. By Lemma
\ref{smallbig} we can assume $|A[y',H]+A[y',H]|\sim 3A(y',H)$.
By Lemma \ref{ftthm}
this implies that $A[y',H]$ is either a full subset of a \bp\ 
of difference $1$ or a full subset of a \bp\ of difference $3$. 
Note that $A(z,H)\sim\frac{1}{2}(H-z)$. So $A[y',H]$ cannot be
a full subset of a \bp\ of difference $3$ because that would imply
$A[z,H]\sim\frac{1}{3}(H-z)\prec\frac{1}{2}(H-z)$.
If $A[y',H]$ is a full subset of the \bp\ $[y',c']\cup [c,H]$,
then $c\prec H$ because otherwise $A[y',H]$ is a forward triangle, 
which contradicts $z\prec H$. However, $c\prec H$ implies that $A[y',H]$ 
is the union of a forward triangle $A[y',2c'-y]$ and a backward triangle
$A[2c'-y+1, H]$, which implies $|2A|\succ 3|A|$ by Lemma \ref{ftbtmore}.
\quad $\Box$(Subcase \ref{middleu3}.2.1)

\medskip

{\bf Subcase \ref{middleu3}.2.2}:\quad
$\underline{d}_{y'+U}(A)<\frac{1}{2}$.

By Lemma \ref{triangle} there are $y'\leqslant z\prec z'\leqslant
H$ such that $A[z,z']$ is a backward triangle and
$A(z',H)\sim\frac{1}{2}(H-z'+1)$. By Lemma \ref{ftthm} we 
have $|A[0,z']+A[0,z']|\succ 3A(0,z')$ because $A[0,z']$
cannot be a full subset of a \bp\ of difference $1$ or $3$. 
Hence $|2A|\succ 3|A|$ by Lemma \ref{smallbig}. 
\quad $\Box$(Subcase \ref{middleu3}.2.2)

\medskip

{\bf Subcase \ref{middleu3}.2.3}:\quad
$\underline{d}_{y'+U}(A)=\frac{1}{2}$.

If for every $x\succ y'$ in $A$, $\gcd(A[x,H]-x)>1$, then there is
an $x\sim y'$ such that $A[x,H]$ is a subset of an \ap\ of
difference $2$. Hence the lemma follows from Lemma \ref{btap}. So
we can assume that there is a $x\succ H$ in $A$ such that
$\gcd(A[x,H]-x)=1$.

Note that $A(y',y'+3)\geqslant 3$. Hence $(A-y')\cap U$ is neither a
subset of an \ap\ of difference $>1$ nor a subset of a
$U$--unbounded \bp\ of difference $d\not=3$. 
If $(A-y')\cap U$ is not a subset of a $U$--unbounded \bp\ of difference
$3$, then the lemma follows from Lemma \ref{smallbig2} and 
Lemma \ref{smallbig}. If $(A-y')\cap U$ is a subset of a \bp\ of
difference $3$, then again $|A[y',H]+A[y',H]|\succ 3A(y',H)$ by the 
proof of Case \ref{middleu4}.1, which again implies $|2A|\succ 3|A|$ by
Lemma \ref{smallbig}. \quad $\Box$(Case \ref{middleu3}.2)

\medskip

{\bf Case \ref{middleu3}.3}:\quad
$A(l_o,H)\preceq\frac{1}{2}(H-l_o)$.

Since $A(0,l_o)\preceq\frac{1}{2}l_o$, then we have
$A(0,l_o)\sim\frac{1}{2}l_o$ and
$A(l_o,H)\sim\frac{1}{2}(H-l_o)$. Let $u_e=\max (A[0,l_o-1])$.
Then $u_e\sim l_o$ and $A[0,u_e]$ is full.

\medskip

{\bf Subcase \ref{middleu3}.3.1}:\quad
$\underline{d}_{l_o+U}(A)>\frac{1}{2}$.

By Lemma \ref{triangle} there is a $y\succ l_o$ such that
$A[l_o,y]$ is a forward triangle. If $y\sim H$, then the lemma
follows from Lemma \ref{btap}. So we can assume $y\prec H$, which
implies that $A[l_o,H]$ cannot be a full subset of a \bp\ of difference $3$. 
By Lemma \ref{smallbig} we can assume that $|A[l_o,H]+A[l_o,H]|\sim
3A(l_o,H)$. Hence by Lemma \ref{ftthm} $A[l_o,H]$ is a full subset
of a \bp\ $[l_o,z]\cup [z',H]$ for some $l_o\prec z\prec
z'\leqslant H$. If $z'\sim H$, then $A[l_o,H]$ is a forward
triangle, which contradicts $y\prec H$. If $z'\prec H$, then
$A[2z'-H,H]$ is a backward triangle. Now the lemma follows
from Lemma \ref{ftbtmore}. 
\quad $\Box$(Subcase \ref{middleu3}.3.1)

\medskip

{\bf Subcase \ref{middleu3}.3.2}:\quad
$\underline{d}_{l_o+U}(A)<\frac{1}{2}$.

By Lemma \ref{triangle} there are $l_o\leqslant z\prec z'\leqslant
H$ such that $A[z,z']$ is a backward triangle and
$A(z',H)\sim\frac{1}{2}(H-z'+1)$. Now the lemma follows from 
Lemma \ref{ftthm} or Lemma \ref{moreftmore2}.
\quad $\Box$(Subcase \ref{middleu3}.3.2)

\medskip

{\bf Subcase \ref{middleu3}.3.3}:\quad
$\underline{d}_{l_o+U}(A)=\frac{1}{2}$.

If there exists a $x\sim l_o$ in $A$ such that
$\gcd(A[x,H]-x)=d>1$, then $d=2$ and the lemma follows from
Lemma \ref{twoaps}. So we can assume that there is a $x\succ l_o$
in $A$ such that $\gcd(A[x,H]-x)=1$. Since $|2A|\sim 3|A|$, then
$|A[0,l_o]+A[0,l_o]| \sim |A[0,u_e]+A[0,u_e]|+|l_o+A[0,u_e]|
\sim 3A(0,l_o)$ implies that $A[0,u_e]$ is full. Without loss
of generality we can assume $u_e,u_e-2,u_e-4\in A$. Hence
$(A-(u_e-4))\cap U$ is neither a subset of an \ap\ of difference
$>1$ nor a subset of a $U$--unbounded \bp\, Now the lemma follows
from Lemma \ref{kneser}, Lemma \ref{smallbig2}, and Lemma \ref{smallbig}. 
\quad $\Box$(Lemma \ref{middleu3})

\begin{lemma}\label{smallu}

Suppose $\underline{d}_U(A)<\frac{1}{2}$. Then $H+1\leqslant 2|A|-1+2b$.

\end{lemma}

\noindent {\bf Proof}:\quad Since
$\underline{d}_U(A)<\frac{1}{2}$, by (4) of Lemma \ref{triangle} 
there are $0\leqslant x\prec a\leqslant H$ such that $A[x,a]$
is a backward triangle from $x$ to $a$ and $A(a,H)\sim\frac{1}{2}(H-a)$. 
If $x\sim 0$ and $a\sim H$, then the lemma follows from 
Lemma \ref{ftlittle}. If $a\sim H$
and $x\succ 0$, then the lemma follows from Lemma \ref{ftthm}.
If $x\succ 0$ and $a\prec H$, then $A(0,x)\sim\frac{1}{2}x$. Hence 
the lemma follows from Lemma \ref{moreftmore} because $A[0,a]$ cannot
be a full subset of a \bp\ of difference $3$.

So we can now assume $x\sim 0$ and $a\prec H$.

If $\underline{d}_{a+U}(A)>\frac{1}{2}$, then there is a $z\succ
a$ such that $A[a,z]$ is a forward triangle. If $z\sim H$, then
the lemma follows from Lemma \ref{btft}. If $z\prec H$, then
$A[a,H]$ cannot be a full subset of a \bp\ of difference $3$
because $A(z,H)\sim\frac{1}{2}(H-z)$. Hence the lemma 
follows from Lemma \ref{moreftmore}.

If $\underline{d}_{a+U}(A)<\frac{1}{2}$, then there are $a<z\prec
z'\leqslant H$ such that $A(z',H)\sim\frac{1}{2}(H-z')$ and
$A[z,z']$ is a backward triangle. Note that $A[0,z']$ contains two
backward triangle, hence is not a full subset of a \bp\,
By Lemma \ref{ftthm} we have
$|A[0,z']+A[0,z']|\succ 3A(0,z')$. Hence $|2A|\succ 3|A|$ by Lemma
\ref{smallbig}.

Assume $\underline{d}_{a+U}(A)=\frac{1}{2}$. Since $A[0,a]$ is a
backward triangle, we can assume there is a $c\sim a$ in $A$ such
that $A(c,c+3)\geqslant 3$. This implies (1) $(A-c)\cap U$ is not
a subset of an \ap\ of difference $>1$ and (2) if $(A-c)\cap U$ is
a subset of a $U$--unbounded \bp\ of difference $d$, then $d=3$.

Suppose $(A-c)\cap U$ is a subset of a $U$--unbounded \bp\ of
difference $3$. By Case \ref{middleu4}.1 we have
$|A[c,H]+A[c,H]|\succ 3A(c,H)$, which imply $|2A|\succ 3|A|$.

Suppose $(A-c)\cap U$ is not a subset of a $U$--unbounded \bp\ of
difference $3$. If for each $x'\succ c$ in $A$, $\gcd(A[x',H]-x')>1$,
then the lemma follows from Lemma \ref{btap}. So we can assume that
there is a $x'\succ c$ such that $\gcd(A[x',H]-x')=1$.
Now the lemma follows from Lemma \ref{kneser}, Lemma \ref{smallbig2}, 
and Lemma \ref{smallbig}. \quad $\Box$(Lemma \ref{smallu})

\begin{theorem}\label{halfthm}

Let $A\subseteq [0,H]$ be such that $0,H\in A$, $\gcd(A)=1$, $A$
is not a subset of a \bp, $|A|\leqslant\frac{1}{2}(H+1)$,
$|A|\sim\frac{1}{2}H$, and $|2A|=3|A|-3+b$ for
$0\leqslant b\sim 0$. Then $H+1\leqslant 2|A|-1+2b$.

\end{theorem}

\noindent {\bf Proof}:\quad If $\underline{d}_U(A)>\frac{1}{2}$,
then the theorem follows from Lemma \ref{largeu}. If
$\underline{d}_U(A)<\frac{1}{2}$, then the theorem follows from
Lemma \ref{smallu}. Suppose $\underline{d}_U(A)=\frac{1}{2}$.

If $A\cap U$ is a subset of an \ap\ of difference $>1$, then the
theorem follows from Lemma \ref{middleu3}. If $A\cap U$ is a
subset of a $U$--unbounded \bp, then the theorem follows from
Lemma \ref{middleu4}. So we can assume that $A\cap U$ is neither a
subset of an \ap\ of difference $>1$ nor a subset of a
$U$--unbounded \bp\, Hence by Lemma \ref{kneser} for every $x\succ 0$
there exists a $y\in A$ with $0\prec y\prec x$ such that
$(2A)(0,2y)\succ 3A(0,y)$.

If for any $x\succ 0$ in $A$ we have $\gcd(A[x,H]-x)>1$, then the
theorem follows from Lemma \ref{middleu2}. So we can assume that
there is a $x\succ 0$ in $A$ such that $\gcd(A[x,H]-x)=1$. But by
Lemma \ref{smallbig2}, we have $|2A|\succ 3|A|$, which contradicts
the assumption of the theorem. Now we finish the proof. \quad
$\Box$(Theorem \ref{halfthm}).

\begin{remark}

We already proved that Theorem \ref{main} follows from Theorem
\ref{nsamain}. Now Theorem \ref{nsamain} follows from 
Theorem \ref{subsetofbp}, Theorem \ref{lesshalf}, and 
Theorem \ref{halfthm}.

\end{remark}

\section{A Corollary for Upper Asymptotic Density}

In this section we slightly improve the most important part of
the main theorem in \cite{jin3} using Theorem \ref{main}.

For an infinite set $A\subseteq\nat$ the upper asymptotic density
of $A$ is defined by
\[\bar{d}(A)=\limsup_{n\rightarrow\infty}\frac{A(0,n-1)}{n}.\]
In this section we assume that $0\in A$ and $\gcd(A)=1$.
By Theorem \ref{2k-1+b} it is not hard to prove that
$0<\bar{d}(A)<\frac{1}{2}$ implies $\bar{d}(2A)\geqslant\frac{3}{2}\bar{d}(A)$.
In \cite{jin3} the structure of $A$ was characterized when 
$0<\bar{d}(A)<\frac{1}{2}$ and $\bar{d}(2A)=\frac{3}{2}\bar{d}(A)$. 
Next we improve this result by substituting the condition
$\bar{d}(2A)=\frac{3}{2}\bar{d}(A)$ with the condition
$\frac{3}{2}\bar{d}(A)\leqslant\bar{d}(2A)<
\frac{3+\epsilon}{2}\bar{d}(A)$ for some positive real number 
$\epsilon$.

\begin{corollary}\label{up}

Let $0<\epsilon\leqslant\frac{1}{3}$ be the real number in Theorem \ref{main}.
For every real number $\delta$ with $0\leqslant\delta<\epsilon$,
if $0<\bar{d}(A)=\alpha<\frac{1}{2(1+\delta)}$ and 
$\bar{d}(2A)=\frac{3+\delta}{2}\alpha$, then either

(a) there exist $d\geqslant 4$ and $c\in [1,d-1]$ such that
$A\subseteq (d\ast\nat)\cup (c+(d\ast\nat))$ and 
$\frac{6}{(2\delta+3)d}\leqslant\alpha\leqslant\frac{2}{d}$,
or 

(b) for every increasing sequence $\langle h_n:n\in\nat\rangle$ with
$\lim_{n\rightarrow\infty} \frac{A(0,h_n)}{h_n+1}=\alpha$, there
exist two sequences $0\leqslant c_n\leqslant b_n\leqslant h_n$
such that $A\cap [c_n+1,b_n-1]=\emptyset$ for each $n\in\nat$,
\[\limsup_{n\rightarrow\infty}\frac{c_n+h_n-b_n}{h_n}\leqslant
\alpha (1+\delta),\]
\[\mbox{ and }\limsup_{n\rightarrow\infty}
\frac{c_n}{h_n-b_n}\leqslant\frac{\delta}{1-\alpha (1+\delta)}.\]

\end{corollary}

\noindent {\bf Proof}:\quad Let $N$ be any hyperfinite integer
and $H=h_N$ be the term in the internal sequence 
$\langle h_n:n\in\,^*\!\nat\rangle$ from (b). Without loss
of generality we can assume $H\in\,^*\!A$. Let $B=\,^*\!A[0,H]$.
Then $|B|\sim\alpha H$ and 
$|2B|\preceq (3+\delta)|B|\sim 3|B|-3+\delta |B|$. If $B$ is a 
subset of an \ap\ of length $l\sim 2|B|-1+2\delta |B|$,
then $H+1\leqslant l\sim 2(1+\delta)\alpha H\prec H$, which is
absurd. So by Theorem \ref{main} we conclude that $B$ is a 
subset of a \bp\ of difference $d$ of length at most $L\preceq
(1+\delta)|B|$. Now the proofs can be found in \cite{bordes} 
(with $\delta=2\sigma-3$) that
if $d>1$, then $B\subseteq (d\ast\,^*\!\nat)\cup (c+(d\ast\,^*\!\nat))$ 
and $\frac{6}{(2\delta+3)d}\leqslant\alpha\leqslant\frac{2}{d}$, and
if $d=1$, then there are $0\leqslant c\leqslant b\leqslant H$ such that
$A\subseteq [0,c]\cup [b,H]$, \[c+H-b\preceq \alpha (1+\delta)H,\]
\[\mbox{ and } c\preceq\frac{\delta}{1-\alpha (1+\delta)}(H-b).\]
Note that the first inequality is a trivial consequence of 
Theorem \ref{main} and the second inequality indicates that
the interval $[0,c]$ is much shorter than $[b,H]$. 
Since the arguments above are true for every hyperfinite
integer $N$, then the corollary follows from the transfer principle.
\quad $\Box$(Corollary \ref{up})

\begin{remark}

(1) The result in \cite{jin3} mentioned above is a special case
of Corollary \ref{up} with $\delta=0$.

(2) Corollary \ref{up} is very similar to the main theorem in 
\cite{bordes}. The main theorem in \cite{bordes} allows
all $\delta<\frac{1}{3}$ instead of $\delta<\epsilon$ in 
Corollary \ref{up}. However, Corollary \ref{up} allows, for example
$\bar{d}(A)=\alpha\leqslant\frac{3}{8}$ (note that 
$\frac{3}{8}\leqslant\frac{1}{2(1+\epsilon)}<\frac{1}{2(1+\delta)}$), 
instead of $\alpha<\alpha_0$ for a small $\alpha_0>0$ in the
main theorem in \cite{bordes}. The reason for the difference is
that the main theorem in \cite{bordes} is a corollary of
Theorem \ref{bilufreiman} while Corollary \ref{up} is a corollary of
Theorem \ref{main}. It should be interesting to see whether one can
prove Corollary \ref{up} with the condition 
$\bar{d}(2A)<\frac{3+\epsilon}{2}\bar{d}(A)$ replaced by
$\bar{d}(2A)<\frac{5}{3}\bar{d}(A)$. In fact, this is
a corollary of Conjecture \ref{weakconj}.

\end{remark}

\section{Comments and a Conjecture}

The reader might notice that the proof of the case when 
$|A|\prec\frac{1}{2}H$ is much more ``nonstandard'' than the proof of the 
case when $|A|\sim\frac{1}{2}H$, which is combinatorial. However,
the proof of the latter is significantly simplified after the
possibility of $|A|\prec\frac{1}{2}H$ is eliminated.

After reading all the proofs above, the reader should be able to
see the crucial role that Lemma \ref{kneser} plays. In order
to violate the condition $|2A|\sim 3|A|$ one needs only to find a
small segment $A[a,b]$ of the set $A$, which already violates
$(2A)(2a,2b)\sim 3A(a,b)$, as long as the rest of the $A$ at each
side of the segment is not too dense and is not a subset of an
\ap\ of difference $>1$ 
(see the condition of Lemma \ref{smallbig}). So if $A\cap U$
does not have expected structural properties such as
$\underline{d}_U(A)>\frac{1}{2}$, $\underline{d}_U(A)=0$,
$A\cap U$ is a subset of an \ap\
of difference $>1$, or $A\cap U$ is a subset of a $U$--unbounded
\bp, then the segment mentioned above is guaranteed by Lemma
\ref{kneser}. Otherwise $A$ must have one of some 
desired structural properties 
in an interval $[0,x]$ for some $x\succ 0$, which gives us a 
high standing ground to reach our final goal. When $A[0,x]$ has
these structural properties, the proof of the main theorems can be
clearly divided into a few possible cases.

Lemma \ref{kneser} is inspired by Kneser's Theorem (cf.\ \cite{HR}) 
and uses the fact that $U$ is an additive semigroup. This
tool is not available in the standard setting, i.e. an initial segment of 
a finite interval cannot be closed under usual addition. This indicates
that the use of nonstandard analysis in this paper is non-trivial.

Although Theorem \ref{main} is a significant advancement of the current
results, it confirms only a weak version of Conjecture \ref{conjecture}.
It is interesting to see whether the ideas from
nonstandard analysis can play a major role in the ultimate
solution of Conjecture \ref{conjecture}. Many lemmas including Lemma
\ref{kneser} in this paper may be generalized. If these
generalizations are achieved, then one can generalize 
Theorem \ref{main} by allowing $|2A|\leqslant\alpha<\frac{10}{3}|A|$. 
I would like to state that as a conjecture. 
The conjecture stated below should be much easier to prove than
proving Conjecture \ref{conjecture}. However, the
solution of the following conjecture is useful for improving
Corollary \ref{up} and could be the last stepping
stone to the ultimate solution of Conjecture \ref{conjecture}.

\begin{conjecture}\label{weakconj}

For any real number $\alpha$ with $3<\alpha<3+\frac{1}{3}$ there
exists a $K\in\nat$ such that for every finite set $A\subseteq\nat$
with $|A|>K$, if $3|A|-3\leqslant
|2A|=3|A|-3+b\leqslant\alpha|A|$, then $A$ is either a subset of
an \ap\ of length at most $2|A|-1+2b$ or a subset of a \bp\ of
length at most $|A|+b$.

\end{conjecture}

\appendix\section{Appendix}

The following theorem is in \cite{freiman2} and
in \cite[p.28]{nathanson}.

\begin{theorem}[G. A. Freiman]\label{2k-1+b}

Let $A$ be a finite set of integers and $|A|=k$. If
$|2A|=2k-1+b<3k-3$, then $A$ is a subset of an \ap\ of length
at most $k+b$.

\end{theorem}

The following theorem is in \cite{freiman2,bilu2}

\begin{theorem}[G. A. Freiman]\label{3k-3}

Let $A$ be a finite set of integers and $|A|=k$. If $|2A|=3k-3$,
then $A$ is either a subset of an \ap\ of length at most $2k-1$
or a \bp

\end{theorem}

The following theorem is in \cite{freiman2}.

\begin{theorem}[G. A. Freiman]\label{twolines}

Let $A\subseteq {\mathbb Z}^2$ be such that $|A|=k>10$. If
$|2A|=3k-3+b$ for $0\leqslant b<\frac{1}{3}k-2$ and $A$ is not a subset 
of a straight line, then $A$ is $F_2$-isomorphic to a subset of 
$\{(0,0),(1,0),\ldots,(l_1-1,0)\}\cup\{(0,1),(1,1),\ldots,(l_2-1,1)\}$
where $l_1+l_2\leqslant k+b$.

\end{theorem}

The following theorem is in \cite[p.118]{nathanson} and in \cite{LS}.

\begin{theorem}[V. Lev \& P. Y. Smeliansky]\label{A+B}

Let $A$ and $B$ be two finite set of non-negative integers such that $0\in
A\cap B$, $|A|,|B|>1$, $\gcd(A)=1$, $m=\max A$, and $n=\max
B\leqslant m$. If $m=n$, then
$|A+B|\geqslant\min\{m+|B|,|A|+2|B|-3\}$. If $m>n$, then
$|A+B|\geqslant\min\{m+|B|,|A|+2|B|-2\}$.

\end{theorem}

Note that Theorem \ref{A+B} is trivially true when $|B|=1$.
In this paper Theorem \ref{A+B} is used in different variations. 
For example, we can replace the conditions of the theorem by 
$|A|,|B|>1$, $\gcd(A-\min A)=1$,  and $m=\max A-\min A\geqslant 
n=\max B-\min B$. We can also consider that both $A$ and $B$ are subsets of
\ap\ 's of difference $d$ with the conditions $|A|,|B|>1$, 
$\gcd(A-\min A)=d$,  and $m=\frac{1}{d}(\max A-\min A)\geqslant 
n=\frac{1}{d}(\max B-\min B)$.

The next theorem is Bilu's version of Freiman's
famous theorem for the inverse problems about the addition of
finite sets \cite[Theorem 1.2 and Theorem 1.3]{bilu2}. we state
only this weak version in order to make the paper a little shorter. 

\begin{theorem}[Y. Bilu \& G. A. Freiman]\label{bilufreiman}

Let $\sigma<4$, $A$ be a finite subset of integers such that
$k=|A|>6$, and $|2A|\leqslant\sigma k$. Then $A$ is a subset of an
$F_2$--progression $P=P(x_0;x_1,x_2;b_1,b_2)$ such that
$|P|\leqslant c_1k$ for some constant $c_1$ and $b_2<c_2$ 
for some constant $c_2$. The constants $c_1,c_2$ are independent 
of $k$.

\end{theorem}

\end{document}